\documentclass[a4paper,12pt]{amsart}
\usepackage[]{amscd,amsfonts,amsmath,amsxtra,amssymb}
\usepackage{graphics}
\usepackage{extarrows}
\usepackage[all]{xy}
\usepackage{hyperref}
\usepackage{eucal}
\usepackage{rotating}
\usepackage[T1]{fontenc}
\newtheorem{sub}{}[section]
\newtheorem{subsub}{}[sub]
\usepackage{rotating}
\usepackage[active]{srcltx}
\def\ov#1{\overline{#1}}

\def\Hom{\mathop{\rm Hom}\nolimits}

\def\HHom{\mathop{\mathcal Hom}\nolimits}

\def\Ext{\mathop{\rm Ext}\nolimits}
\def\EExt{\mathop{\mathcal Ext}\nolimits}

\def\Pic{\mathop{\rm Pic}\nolimits}

\def\Aut{\mathop{\rm Aut}\nolimits}
\def\AAut{\mathop{\mathcal Aut}\nolimits}

\def\End{\mathop{\rm End}\nolimits}

\def\EEnd{\mathop{\mathcal End}\nolimits}
\def\AAd{\mathop{\mathcal Ad}\nolimits}
\def\GL{\mathop{\rm GL}\nolimits}

\def\imm{\mathop{\rm im}\nolimits}
\def\deg{\mathop{\rm deg}\nolimits}

\def\rk{\mathop{\rm rk}\nolimits}

\def\det{\mathop{\rm det}\nolimits}
\def\spec{\mathop{\rm spec}\nolimits}

\def\lra{\longrightarrow}

\def\sigg{\mathop{\hbox{$\displaystyle\sum$}}\limits}

\def\hfl#1#2{\smash{\mathop{\ \hbox to 12mm{\rightarrowfill}}
\limits^{\scriptstyle#1}_{\scriptstyle#2} \ }}
\def\hflb#1#2{\smash{\mathop{\hbox to 12mm{\leftarrowfill}}
\limits^{\scriptstyle#1}_{\scriptstyle#2}}}

\def\pline#1{<\hskip-3.5pt#1\hskip-3.5pt>}
\def\m#1{{\hbox{$#1$}}}

\def\ot{\otimes}
\def\og{\leavevmode\raise.3ex\hbox{$\scriptscriptstyle\langle\!\langle$}}
\def\fg{\leavevmode\raise.3ex\hbox{$\scriptscriptstyle\,\rangle\!\rangle$}}

\def\nsp{\lbrace 0\rbrace}

\def\dsp{\displaystyle}
\def\Ssect#1#2{\pagebreak[3]\begin{sub}\label{#2}{\sc\small\small
#1}\rm\medskip}

\def\sepsec{\vskip 1.4cm}
\def\sepsub{\vskip 0.7cm}
\def\sepsubsub{\vskip 0.5cm}
\def\sepprop{\vskip 0.5cm}

\def\xmat#1{\[\xymatrix{#1}\]}
\def\flinc{\ar@{^{(}->}}
\def\fleq{\ar@{=}}
\def\flon{\ar@{->>}}
\def\fmaps{\ar@{|-{>}}}
\def\fflat{\ar@{-}}
\def\fimpl{\ar@{=>}}
\def\Nligne{\hfil\break}

\def\ED{\vskip 1cm\end{document}}

\def\Inv{\text{Inv}}

\newcommand{\Z}{{\mathbb Z}}

\newcommand{\C}{{\mathbb C}}

\renewcommand{\P}{{\mathbb P}}
\newcommand{\D}{{\mathbb D}}

\newcommand{\E}{{\mathbb E}}

\renewcommand{\L}{{\mathbb L}}

\newcommand{\ka}{{\mathcal A}}

\newcommand{\kc}{{\mathcal C}}

\newcommand{\ke}{{\mathcal E}}
\newcommand{\kf}{{\mathcal F}}
\newcommand{\kg}{{\mathcal G}}
\newcommand{\kh}{{\mathcal H}}
\newcommand{\ki}{{\mathcal I}}

\newcommand{\kl}{{\mathcal L}}

\newcommand{\ko}{{\mathcal O}}

\newcommand{\kt}{{\mathcal T}}

\newcommand{\kx}{{\mathcal X}}

\begin{document}

\def\refname{References}
\def\contentsname{Summary}
\def\proofname{Proof}
\def\abstractname{Resume}

\author{Jean--Marc Dr\'{e}zet}
\address{
Institut de Math\'ematiques de Jussieu - Paris Rive Gauche\\
Case 247\\
4 place Jussieu\\
F-75252 Paris, France}
\email{jean-marc.drezet@imj-prg.fr}
\title[{Primitive multiple schemes}] {Primitive multiple schemes}

\begin{abstract}
A primitive multiple scheme is a Cohen-Macaulay scheme $Y$ such that the 
associated reduced scheme $X=Y_{red}$ is smooth, irreducible, and that $Y$ can 
be locally embedded in a smooth variety of dimension \m{\dim(X)+1}. If $n$ is 
the multiplicity of $Y$, there is a canonical filtration \m{ X=X_1\subset
X_2\subset\cdots\subset X_n=Y}, such that \m{X_i} is a primitive multiple 
scheme of multiplicity $i$. The simplest example is the trivial primitive 
multiple scheme of multiplicity $n$ associated to a line bundle $L$ on $X$: it 
is the $n$-th infinitesimal neighborhood of $X$, embedded in the line bundle 
$L^*$ by the zero section.

Let ${\bf Z}_n={\rm spec}(\C[t]/(t^n))$. The 
primitive multiple schemes of multiplicity $n$ are obtained by taking an open 
cover $(U_i)$ of a smooth variety $X$ and by gluing the schemes
$U_i\times{\bf Z}_n$ using automorphisms of $U_{ij}\times {\bf Z}_n$ that leave 
$U_{ij}$ invariant. This leads to the study of the sheaf of nonabelian groups 
$\kg_n$ of automorphisms of $X\times {\bf Z}_n$ that leave the $X$ invariant, 
and to the study of its first cohomology set. If \m{n\geq 2} there is an 
obstruction to the extension of \m{X_n} to a primitive multiple scheme of 
multiplicity $n+1$, which lies in the second cohomology group \m{H^2(X,E)} of a 
suitable vector bundle $E$ on $X$. 

In this paper we study these obstructions and the parametrization of primitive 
multiple schemes. As an example we show that if \m{X=\P_m} with \m{m>=3} all 
the primitive multiple schemes are trivial. If \m{X=\P_2}, there are only two 
non trivial primitive multiple schemes, of multiplicities $2$ and $4$, which 
are not quasi-projective. On the other hand, if $X$ is a projective bundle over 
a curve, we show that there are infinite sequences \ \m{ X=X_1\subset
X_2\subset\cdots\subset X_n\subset X_{n+1}\subset\cdots} \ of non trivial 
primitive multiple schemes.
\end{abstract}

\maketitle
\tableofcontents


\vskip 1cm

\section{Introduction}\label{intro}

A {\em multiple scheme} is a Cohen-Macaulay scheme $Y$ over $\C$ such that 
\m{Y_{red}=X} is a smooth connected variety. We call $X$ the {\em support} of 
$Y$. If $Y$ is quasi-projective we say that it is a {\em multiple variety with 
support $X$}. In this case $Y$ is projective if $X$ is.

Let \m{\ki_X} be the ideal sheaf of $X$ in $Y$.
Let $n$ be the smallest integer such that \m{Y=X^{(n-1)}}, \m{X^{(k-1)}}
being the $k$-th infinitesimal neighborhood of $X$, i.e. \
\m{\ki_{X^{(k-1)}}=\ki_X^{k}} . We have a filtration \ \m{X=X_1\subset
X_2\subset\cdots\subset X_{n}=Y} \ where $X_i$ is the biggest Cohen-Macaulay
subscheme contained in \m{Y\cap X^{(i-1)}}. We call $n$ the {\em multiplicity}
of $Y$.

Let \m{d=\dim(X)}.
We say that $Y$ is {\em primitive} if, for every closed point $x$ of $X$,
there exists a smooth variety $S$ of dimension \m{d+1}, containing a 
neighborhood of $x$ in $Y$ as a locally closed subvariety. In this case, 
\m{L=\ki_X/\ki_{X_2}} is a line bundle on $X$, \m{X_j} is a primitive multiple 
scheme of multiplicity $j$ and we have \ 
\m{\ki_{X_j}=\ki_X^j}, \m{\ki_{X_{j}}/\ki_{X_{j+1}}=L^j} \ for \m{1\leq j<n}. 
We call $L$ the line bundle on $X$ {\em associated} to $Y$. The ideal sheaf 
\m{\ki_{X}} can be viewed as a line bundle on \m{X_{n-1}}. If \m{n=2}, $Y$ is
called a {\em primitive double scheme}.

The simplest case is when $Y$ is contained in a smooth variety $S$ of dimension 
\m{d+1}. Suppose that $Y$ has multiplicity $n$. Let \m{P\in X} and 
\m{f\in\ko_{S,P}}  a local equation of $X$. Then we have \ 
\m{\ki_{X_i,P}=(f^{i})} \ for \m{1<j\leq n} in $S$, in particular 
\m{\ki_{Y,P}=(f^n)}, and \ \m{L=\ko_X(-X)} .

For any \m{L\in\Pic(X)}, the {\em trivial primitive variety} of multiplicity 
$n$, with induced smooth variety $X$ and associated line bundle $L$ on $X$ is 
the $n$-th infinitesimal neighborhood of $X$, embedded by the zero section in 
the dual bundle $L^*$, seen as a smooth variety.

The case \m{\dim(X)=1} is well known. The primitive multiple curves have been 
defined and studied in \cite{ba_fo}. Double curves had previously been used in 
\cite{fe}. The primitive double schemes have then been constructed and 
parametrized in \cite{ba_ei}. Further work on primitive multiple curves has 
been done in \cite{ch-ka}, \cite{dr1b}, \cite{dr2}, \cite{dr1}, \cite{dr4},  
\cite{dr5}, \cite{dr6}, \cite{dr7}, \cite{dr8} and \cite{dr9}, where higher 
multiplicity is treated, as well as coherent sheaves on primitive curves and 
their deformations. See also \cite{gonz1}, \cite{ga_go_pu0}, \cite{sa1}, 
\cite{sa2} and \cite{sa3}. Primitive double schemes of greater dimension are 
also studied in\cite{b_g_g}, \cite{b_m_r}, \cite{ga_go_pu1} and 
\cite{ga-go-pu}.

If \m{\dim(X)\geq 2}, the main difference with the case of curves is the 
presence of {\em obstructions} : given \m{X_n} they are
\begin{enumerate}
\item[--] Obstruction to the extension of a vector bundle on $X_{n-1}$ to a 
vector bundle on $X_n$.
\item[--] Obstruction to the extension of $X_n$ to a primitive multiple scheme 
of multiplicity $n+1$.
\end{enumerate}
We will see that these obstructions depend on the vanishing of elements in 
cohomology groups \m{H^2(X,E)}, where $E$ is a suitable vector bundle on $X$. 
Hence if \m{\dim(X)=1} the obstructions disappear, and it is always possible to 
extend vector bundles or primitive multiple schemes. This is why one can 
obtain many primitive multiple curves. 

On the other hand in some cases, when \m{\dim(X)\geq 2}, it may happen that non 
trivial primitive multiple schemes are very rare. For example, we will see that 
if \m{X=\P_m}, \m{m\geq 2}, then there are exactly two non trivial primitive 
multiple schemes (for \m{m=2}, in multiplicities 2 and 4), and these two 
schemes are not quasi-projective (the scheme of multiplicity 2 appears already 
in \cite{b_m_r}, see also \cite{ga-go-pu}, 2-). But in some other cases (for 
example if $X$ is of dimension 2 and is a projective bundle over a curve) there 
are big families and infinite sequences
\[X=X_1\subset X_2\subset\cdots\subset X_{n}\subset X_{n+1}\subset\cdots\]
of non trivial primitive multiple schemes.

If \m{\dim(X)\geq 2}, it may also happen that primitive multiple schemes are 
not quasi-projective (this has also been observed in \cite{b_m_r}, 
\cite{ga-go-pu}). In fact, \m{X_n} is quasi-projective (and even 
projective) if and only if it is possible to extend an ample line bundle on $X$ 
to a line bundle on \m{X_n}. For example the two non trivial primitive multiple 
schemes for \m{X=\P_2} are not quasi-projective.

\sepprop

I am grateful to the anonymous referee for giving some useful remarks. 

\sepprop

{\em Notations:} -- In this paper, an {\em algebraic variety} is a 
quasi-projective scheme over $\C$. A {\em vector bundle} on a scheme is an 
algebraic vector bundle.

-- If $X$ is a scheme and \m{P\in X} is a closed point, we denote 
by \m{m_{X,P}} (or \m{m_P}) the maximal ideal of $P$ in \m{\ko_{X,P}}.

-- If $X$ is an scheme and \m{Z\subset X} is a closed subscheme, 
\m{\ki_Z} (or \m{\ki_{Z,X}}) denotes the ideal sheaf of $Z$ in $X$.

-- If \m{V} is a finite dimensional complex vector space, \m{\P(V)} denotes the 
projective space of lines in $V$, and we use a similar notation for projective 
bundles.

\sepprop

Let \m{X_n} be a primitive multiple scheme of multiplicity $n$, 
\m{X=(X_n)_{red}}, 
and \m{L\in\Pic(X)} the associated line bundle. 

\sepprop

\Ssect{Construction and parametrization of primitive multiple schemes}{intro-1}

\begin{subsub} The sheaves of groups associated to primitive multiple schemes 
-- \rm
For every open subset \m{U\subset X}, let \m{U^{(n)}} denote the corresponding 
open subset of \m{X_n}. Let \ \m{{\bf Z}_n=\spec(\C[t]/(t^n))}. It is proved in 
\cite{dr1}, th. 5.2.1,  that \m{X_n} is locally trivial, i.e. there exists an 
open cover 
\m{(U_i)_{i\in I}} of $X$ such that for evey \m{i\in I} there is an isomorphism
\[\delta_i:U_i^{(n)}\lra U_i\times{\bf Z}_n\]
inducing the identity on \m{U_i}. For every open subset $U$ of $X$, we have
\ \m{\ko_{X\times{\bf Z}_n}(U)=\ko_X(U)[t]/(t^n)}. It is then natural to 
introduce the sheaf of (non abelian) groups \m{\kg_n} on $X$, where 
\m{\kg_n(U)} is the group of automorphisms $\theta$ of the $\C$-algebra 
\m{\ko_X(U)[t]/(t^n)} such that for every \m{\alpha\in\ko_X(U)[t]/(t^n)}, 
\m{\theta(\alpha)_{|U}=\alpha_{|U}}. We have \m{\kg_1=\ko_X^*}.

Let
\xmat{\delta_{ij}=\delta_j\delta_i^{-1}:U_{ij}\times{\bf Z}_n\ar[r]^-\simeq & 
U_{ij}\times{\bf Z}_n , }
and \ 
\m{\delta_{ij}^*=\delta_i^{*-1}\delta_j^*\in\ko_X(U_{ij})[t]/(t^n)}.
Then \m{(\delta_{ij}^*)} is a cocycle of \m{\kg_n}, which describes completely 
\m{X_n}. 

In this way we see that there is a canonical bijection between the cohomology 
set \m{H^1(X,\kg_n)} ans the set of isomorphism classes of primitive multiple 
schemes \m{X_n} such that \m{X=(X_n)_{red}}. Most results of this paper come 
from computations in {\em \v Cech cohomology}.

There is an obvious morphism \ 
\m{\rho_{n+1}:\kg_{n+1}\to\kg_n}, such that
\[H^1(\rho_{n+1}):H^1(X,\kg_{n+1})\lra H^1(X,\kg_n)\]
sends \m{X_{n+1}} to \m{X_n} \ if \m{n\geq 2}, whereas 
\[H^1(\rho_2):H^1(X,\kg_2)\lra H^1(X,\ko_X^*)=\Pic(X)\]
sends \m{X_2} to $L$.

We prove that \m{\ker(\rho_2)\simeq T_X} (proposition \ref{g2_lem2}) and
\ \m{\ker(\rho_{n+1})\simeq T_X\oplus\ko_X} \ if \m{n\geq 2} (proposition 
\ref{gn_lem1}), where \m{T_X} is the tangent bundle of $X$. The fact that they 
are sheaves of {\em abelian} groups allows to compute {\em obstructions}. Let 
\m{g_n\in H^1(X,\kg_n)}, corresponding to the primitive multiple scheme \m{X_n},
and for \m{1\leq i<n}, \m{g_i} the image of \m{g_n} in \m{H^1(X,\kg_i)}. Let 
\m{\ker(\rho_{n+1})^{g_n}} be the associated sheaf of groups (cf. 
\ref{coh_gr}). We prove that
\[\ker(\rho_2)^{g_1} \ \simeq \ T_X\ot L  \] 
(proposition \ref{g2_lem1}). By cohomology theory, we find that there is a 
canonical surjective map
\[H^1(X,T_X\ot L)\lra H^1(\rho_2)^{-1}(L)\]
sending 0 to the trivial primitive scheme, and whose fibers are the orbits of 
the action of $\C^*$ on \m{H^1(X,T_X\ot L)} by multiplication. Hence there is a 
bijection between the set of non trivial double schemes with associated line 
bundle $L$, and \m{\P(H^1(X,T_X\ot L))}. For a non trivial double scheme 
\m{X_2}, we have an exact sequence 
\[0\lra L\lra\Omega_{X_2|X}\lra\Omega_X\lra 0 \ , \]
corresponding to \ \m{\sigma\in\Ext^1_{\ko_X}(\Omega_X,L)=H^1(T_X\ot L)}, and 
\m{\C\sigma} is the element of \m{\P(H^1(T_X\ot L))} corresponding to \m{X_2}.
This result has already been proved in \cite{ba_ei} using another method, and in 
\cite{dr1} when $X$ is a curve, in the same way as here.

If \m{n>2} we have
\[\ker(\rho_{n+1})^{g_n} \ \simeq \ (\Omega_{X_2|X})^*\ot L^n\]
(proposition \ref{g2_lem2}). We have then (by the theory of cohomology of 
sheaves of groups) an {\em obstruction map}
\[\Delta_n:H^1(X,\kg_n)\lra H^2((\Omega_{X_2|X})^*\ot L^n)\]
such that \ \m{g_n\in\imm(H^1(\rho_{n+1}))} \ if and only \ \m{\Delta_n(g_n)=0}.

Suppose that \m{n>2} and \ \m{\Delta_n(g_n)=0}, and let \m{g_{n+1}\in 
H^1(X,\kg_{n+1})} be such that \Nligne \m{H^1(\rho_{n+1})(g_{n+1})=g_n}. Then 
\m{H^1(\rho_{n+1})^{-1}(g_n)} is 
the set of extensions of \m{X_n} to a primitive multiple scheme of multiplicity 
\m{n+1}. Let \m{\Aut(X_n)} be the set of automorphisms of \m{X_n} inducing the 
identity on $X$. There is a canonical surjective map
\[\lambda_{g_{n+1}}:H^1(X,(\Omega_{X_2|X})^*\ot L^n)\lra 
H^1(\rho_{n+1})^{-1}(g_n)\]
which sends 0 to \m{g_n}, and whose fibers are the orbits of an action of 
\m{\Aut(X_n)} on\Nligne \m{H^1(X,(\Omega_{X_2|X})^*\ot L^n)}.

\end{subsub}

\sepprop

\begin{subsub} Extensions of the ideal sheaf of $X$ -- \rm The ideal sheaf 
\m{\ki_{X,X_n}} is a line bundle on \m{X_{n-1}}. A necessary condition to 
extend \m{X_n} to a primitive multiple scheme \m{X_{n+1}} of multiplicity 
\m{n+1} is that \m{\ki_{X,X_n}} can be extended to a line bundle on \m{X_n} 
(namely \m{\ki_{X,X_{n+1}}}). This is why we can consider pairs \m{(X_n,\L)}, 
where $\L$ is a line bundle on \m{X_n} such that \ 
\m{\L_{|X_{n-1}}\simeq\ki_{X,X_n}}. 

The corresponding sheaf of groups on $X$ is defined as follows: for every open 
subset \m{U\subset X}, \m{\kh_n(U)} is the set of pairs \m{(\phi,u)}, where 
\m{\phi\in\kg_n(U)}, and \m{u\in\ko_X(U)[t]/(t^n)} is such that \ 
\m{\phi(t)=ut} (with an adequate definition of the group law, cf. 
\ref{prim_can}). The set of isomorphism classes of the above pairs \m{(X_n,\L)} 
can then be identified with the cohomology set \m{H^1(X,\kh_n)}.

There is an obvious morphism \ \m{\tau_n:\kg_{n+1}\to\kh_n}, such that
\[H^1(\tau_n):H^1(X,\kg_{n+1})\lra H^1(X,\kh_n)\]
sends \m{X_{n+1}} to \m{(X_n,\ki_{X,X_{n+1}})}.
Let \m{g\in H^1(X,\kh_n)}. We prove in proposition \ref{prop5} that \ 
\m{\ker(\tau_n)\simeq T_X} \ and
\[\ker(\tau_n)^g \ \simeq \ T_X\ot L^n \ . \]
Consequently there is again, by cohomology theory, an {\em obstruction map}
\[\Delta''_n:H^1(X,\kh_n)\lra H^2(T_X\ot L^n)\]
such that, if \m{(X_n,\L)} corresponds to $g$, there is an extension of \m{X_n} 
to a primitive multiple scheme \m{X_{n+1}} of multiplicity \m{n+1} with \ 
\m{\ki_{X,X_{n+1}}\simeq\L} \ if and only if \ \m{\Delta''_n(g)=0}.

Suppose \ \m{\Delta''_n(g)=0}, and let \m{g_{n+1}\in H^1(X,\kg_{n+1})} be such 
that \ \m{H^1(\tau_n)(g_{n+1})=g}. Then \m{H^1(\tau_n)^{-1}(g)} is 
the set of extensions of \m{X_n} to a primitive multiple scheme \m{X_{n+1}} of 
multiplicity \m{n+1}, such that \ \m{\ki_{X,X_{n+1}}\simeq\L}. Let 
\m{\Aut^\L(X_n)} be the set of automorphisms $\sigma$ of \m{X_n} inducing the 
identity on $X$, and such that \m{\sigma^*(\L)\simeq\L}. There is a canonical 
surjective map
\[\lambda''_{g_{n+1}}:H^1(X,T_X\ot L^n)\lra H^1(\rho_{n+1})^{-1}(g_n)\]
which sends 0 to \m{g_{n+1}}, and whose fibers are the orbits of an action of 
\m{\Aut^\L(X_n)} on\Nligne \m{H^1(X,T_X\ot L^n)}.

\end{subsub}

\end{sub}

\sepsub

\Ssect{Automorphisms of primitive multiple schemes}{intro-2}
\end{sub}

The group \m{\Aut(X_n)} appears in the parametrization of the extensions of 
\m{X_n} to multiplicity \m{n+1}. Let \m{\Aut_0(X_n)} be group of 
automorphisms of \m{X_n} such that the induced automorphism of $L$ is the 
identity. Let
\[\kt_n \ = \ (\Omega_{X_n})^* \ . \]
We prove that there is a natural bijection
\[\Aut_0(X_n)\lra H^0(X_n,\ki_X\kt_n)\]
(theorem \ref{prop22}, unless \m{n=2} this is not a morphism of groups). 
Locally this corresponds to the following result: let \m{U\subset X} be 
an open subset such that \m{\Omega_U} is trivial. For every derivation \m{D_0} 
of \m{\ko_X(U)[t]/(t^n)}, if \m{D=tD_0}, then
\[\chi_D \ = \ \sigg_{k\geq 0}\frac{1}{k!}D^k:\ko_X(U)[t]/(t^n)\lra
\ko_X(U)[t]/(t^n)\]
is an element of \m{\kg_n(U)} such that \m{\chi_D(t)} is of the form \ 
\m{\chi_D(t)=(1+\beta t)t} (for some \m{\beta\in\ko_X(U)[t]/(t^n)}). Moreover 
for every \m{\chi\in\kg_n(U)} with this property, there exists a unique 
derivation $D$ multiple of $t$ such that \m{\chi=\chi_D} (theorem \ref{prop12}).

We show that if \m{X_n} is not trivial, and \m{\Aut_0(X_n)} is trivial, then 
\m{\Aut(X_n)} is finite (theorem \ref{prop21} and corollary \ref{coro3}). 

\sepsub

\Ssect{Extensions and obstructions}{intro-3}

\begin{subsub} Extension of vector bundles -- \rm Let $r$ be a positive 
integer. If $Y$ is a primitive multiple scheme, let \m{\GL(r,\ko_Y)} denote the 
sheaf of groups on $Y$ of invertible \m{r\times r}-matrices with coefficients 
in $Y$ (the group law being the multiplication of matrices). If $\ke$ is a 
coherent sheaf on $Y$, let \m{M(r,\ke)} denote the sheaf of groups on $Y$ of 
matrices with coefficients in $\ke$ (the group law being the addition of 
matrices). There is a canonical bijection between \m{H^1(Y,\GL(r,\ko_Y))} and 
the set of isomorphism classes of rank $r$ vector bundles on $Y$.

Suppose that \m{X_n} can be extended to a primitive multiple scheme \m{X_{n+1}} 
of multiplicity \m{n+1}. The kernel of the canonical morphism \ 
\m{p_n:\GL(r,\ko_{X_{n+1}})\to\GL(r,\ko_{X_n})} \ is isomorphic to 
\m{M(r,L^n)}. Let \ \m{g\in H^1(X_n,\GL(r,\ko_{X_n}))}, corresponding to the 
vector bundle $\E$ on \m{X_n}. Let \m{E=\E_{|X}}. We prove that the associated 
sheaf \m{M(r,L^n)^g} is isomorphic to \m{E^*\ot E\ot L^n}. It follows that we 
have a canonical obstruction map
\[\Delta:H^1(X_{n},\GL(r,\ko_{X_n}))\lra H^2(X,E^*\ot E\ot L^n)\]
such that \ \m{g\in\imm(H^1(p_n))} \ if and only \ \m{\Delta(g)=0}. Hence 
\m{\Delta(g)} is the obstruction to the extension of $\E$ to \m{X_{n+1}}. If 
$\E$ is a line bundle, this obstruction lies in \m{H^2(X,L^n)}.

We have a canonical exact sequence of sheaves on \m{X_n}
\[0\lra L^n\lra\Omega_{X_{n+1}|X_n}\lra\Omega_{X_n}\lra 0 \ , \]
corresponding to \ \m{\sigma_n\in\Ext^1_{\ko_{X_n}}(\Omega_{X_n},L^n)}, 
inducing \ \m{\ov{\sigma}_n\in\Ext^1_{\ko_{X_n}}(\E\ot\Omega_{X_n},\E\ot L^n)}. 
Let \ \m{\nabla_0(E)\in\Ext^1(\E,\E\ot\Omega_{X_n})} \ be the canonical class 
of $\E$ (cf. \ref{can_cl}). We 
have a canonical product
\[\Ext^1_{\ko_{X_n}}(\E\ot\Omega_{X_n},\E\ot L^n)\times
\Ext^1_{\ko_{X_n}}(\E,\E\ot\Omega_{X_n})\lra
\Ext^2_{\ko_{X_n}}(\E,\E\ot L^n)=H^2(X,E^*\ot E\ot L^n) \ . \]
We prove that  \ \m{\Delta(g)=\ov{\sigma}_n\nabla_0(\E)} (theorem \ref{prop4}). 
\end{subsub}

\sepprop

\begin{subsub} Extension of primitive multiple schemes -- \rm Suppose that
\m{n=2}, and that \m{\ki_{X,X_2}} can be extended to a line bundle $\L$ on 
\m{X_2}. Let \m{g_2\in H^1(X,\kh_2)} be associated to \m{X_2}. We give in 
theorem \ref{prop8} a description of the obstruction map
\[\Delta''_2:H^1(X,\kh_2)\lra H^2(X,T_X\ot L^2)\]
corresponding to \m{(X,\L)}, used to decide if \m{X_2} can be extended to 
multiplicity 3. This description is made with \v Cech cohomology. 

We could not find a simpler expression for \m{\Delta''_2(g_2)}, except in the 
case \m{L=\ko_X}. The parametrization of extensions of \m{\ko_X} to a line 
bundle on \m{X_2} exhibits an element \m{\theta_\L\in H^1(\ko_X)} associated to 
$\L$. On the other hand, we have \m{\sigma\in H^1(T_X)} associated to \m{X_2} 
(from the construction of double primitive schemes, or by the exact sequence \ 
\m{0\to L\to\Omega_{X_2|X}\to\Omega_X\to 0}). We find that
\[\Delta''_2(g_2) \ = \ \theta_\L\sigma+\frac{1}{2}[\sigma,\sigma] \ . \]
(where \m{[-,-]} is the Poisson bracket).

Suppose that \m{X_2} is trivial. We give in \ref{hi_mu} a complete description 
of the action of \m{\Aut(X_2)} on \m{H^1(X,(\Omega_{X_2|X})^*\ot L^2)} (cf 
.\ref{intro-1}).

For any \m{n\geq 2}, if \m{X_n} is trivial, we describe in \ref{hi_mu2} the 
action of \m{\C^*\subset\Aut(X_n)} on\Nligne \m{H^0(X,(\Omega_{X_2|X})^*\ot 
L^n)}. It follows that if \m{h^0(X_n,\ki_X\kt_n)=0}, then 
\m{\Aut(X_n)\simeq\C^*}, and  the extensions of \m{X_n} to multiplicity \m{n+1} 
are parametrized by a {\em weighted projective space} (cf. \ref{weigh}).

\end{subsub}

\end{sub}

\sepsub

\Ssect{Examples}{intro-4}

\begin{subsub} The case of surfaces and canonical associated line bundle -- \rm
Suppose that $X$ is a surface and \m{L=\omega_X}. Let \m{X_2} be a non trivial 
primitive double scheme with associated line bundle \m{\omega_X}. As indicated 
in \ref{intro-1}, such schemes are parametrized by \m{\P(H^1(\Omega_X))}. Let 
\m{\sigma\in H^1(\Omega_X)} such that \m{\C\sigma} corresponds to \m{X_2}. We 
prove in proposition \ref{prop28} that if \m{H^0(T_X)=\nsp}, then \m{X_2} can 
be extended to a primitive multiple scheme of multiplicity 3 if and only if \ 
\m{\sigma.\nabla_0(\omega_X)=0} \ in \m{H^2(\omega_X)}.

\end{subsub}

\sepprop

\begin{subsub} The case of projective spaces -- \rm We suppose that \m{X=\P_m}, 
\m{m\geq 2}. We prove that if \m{m>2}, then there are only trivial primitive 
multiple schemes, and if \m{m=2}, there are exactly two non trivial primitive 
multiple schemes (theorem \ref{theo4}). The first, \m{{\bf X}_2}, is of 
multiplicity 2, with \m{L=\ko_{\P_2}(-3)}. The second, \m{{\bf X_4}}, is of 
multiplicity 4, with \m{L=\ko_{\P_2}(-1)}.

We prove that the only line bundles on \m{{\bf X}_2} and \m{{\bf X}_4} are the 
trivial line bundles. Consequently, \m{{\bf X}_2} and \m{{\bf X}_4} are not 
quasi-projective (theorems \ref{theo2} and \ref{theo3}). Since a primitive 
multiple scheme is a locally complete intersection, it has a dualizing 
sheaf which is a line bundle. We must have \ \m{\omega_{{\bf 
X}_2}\simeq\ko_{{\bf X}_2}} \ and \ \m{\omega_{{\bf 
X}_4}\simeq\ko_{{\bf X}_4}}, i.e. \m{{\bf X}_2} \m{{\bf X}_4} are {\em K3 
carpets} according to definition 1.2 of \cite{ga-go-pu}.

\end{subsub}

\sepprop

\begin{subsub} The case of projective bundles over curves -- \rm Let $C$ be a 
smooth irreducible projective curve and $E$ a rank 2 vector bundle on $C$. Let 
\ \m{X=\P(E)} \ be the associated projective bundle and \ \m{\pi:\P(E)\to C} \ 
the canonical projection. We suppose that
\[L \ = \ \pi^*(D)\ot\ko_{\P(E)}(-k) \ , \]
with \m{D\in\Pic(C)} and \m{k\geq 3}. We prove in \ref{pro_c} that if 
\m{\deg(D)} is sufficiently big, then there exist infinite sequences
\[X_1=\P(E)\subset X_2\subset\cdots\subset X_n\subset X_{n+1}\subset\cdots \]
of projective non trivial primitive multiple schemes (\m{X_n} of multiplicity 
$n$). 

For example, if \m{C=\P_1}, \m{E=\ko_{\P_1}\ot\C^2}, we have 
\m{X=\P_1\times\P_1}. Let \m{L=\ko(d,-k)}, with \m{d>0}, \m{k\geq 3}. Then such 
infinite sequences exist, and for every \m{n>1}, the extensions of \m{X_n} to 
primitive multiple schemes of multiplicity \m{n+1} form a family of dimension \ 
\m{d(3kn^2-5n-2kn+k+1)+5kn-7-k}.
\end{subsub}

\end{sub}

\sepsub

\Ssect{Outline of the paper}{intro-5}

In chapter 2, we give several ways to use \v Cech cohomology, and some basic 
results on cohomology of sheaves of groups are recalled. We give also some 
technical results about polynomials that are used later.

In chapter 3 the definition and some properties of the canonical class of a 
vector bundle are given.

In chapter 4 we give the definition and the properties of the sheaves of groups 
that are used to define primitive multiple schemes, as described in 
\ref{intro-1}.

The chapter 5 contains a detailed study of the groups \m{\kg_n(U)} of 
\ref{intro-1}, and of the canonical sheaves of primitive multiple schemes. 
These results are used is to perform computations in \v Cech cohomology. 

In chapter 6 we give some descriptions and properties of the automorphism group 
of a primitive multiple scheme. Some demonstrations are omitted when they are 
similar to the case of primitive multiple curves, treated in \cite{dr1}.

In chapter 7 we study the two extensions problems: extension of vector bundles 
and extension of primitive multiple schemes to higher multiplicity. It contains 
also a comparison of the present work and of the other possible construction of 
primitive multiple schemes given in \cite{ba_ei}.

In Chapter 8 we apply the previous results to the case \m{X=\P_m}, \m{m\geq 2}.

In chapter 9 we apply the previous results to the case \m{X=\P(E)}, where $E$ 
is a rank 2 vector bundle on a smooth projective curve.

\end{sub}

\sepsec

\section{Preliminaries}

\Ssect{\v Cech cohomology and trivializations}{cech}

Let $X$ be a scheme over $\C$ and $E$, $F$ vector bundles on $X$ of ranks 
$r$, $s$ respectively. Let \m{(U_i)_{i\in I}} be an open cover of $X$ such that 
we have trivializations:
\xmat{\alpha_i:E_{|U_i}\ar[r]^-\simeq & \ko_{U_i}\ot\C^r .}

\sepprop

\begin{subsub}\label{cech1}\rm Let
\xmat{\alpha_{ij}=\alpha_i\circ\alpha_j^{-1}:\ko_{U_{ij}}\ot\C^r\ar[r]^-\simeq &
\ko_{U_{ij}}\ot\C^r ,}
so that we have the relation \ \m{\alpha_{ij}\alpha_{jk}=\alpha_{ik}}.
Let $n$ be a positive integer, and for every sequence \m{(i_0,\ldots,i_n)} of 
distinct elements of $I$, \m{\sigma_{i_0\cdots i_n}\in H^0(U_{i_0\cdots 
i_n},\ko_X\ot\C^r)}. Let
\[\theta_{i_0\cdots i_n} \ = \ \alpha_{i_0}^{-1}\sigma_{i_0\cdots i_n} \ \in \
H^0(U_{i_0\cdots i_n},E) \ . \]
The family \m{(\theta_{i_0\cdots i_n})} represents an element of \m{H^n(X,E)} 
if the cocycle relations are satisfied: for every sequence  
\m{(i_0,\ldots,i_{n+1})} of distinct elements of $I$, 
\[\sigg_{k=0}^n(-1)^k\theta_{i_0\cdots\widehat{i_k}\cdots i_{n+1}} \ = \ 0 \ , 
\]
which is equivalent to
\[\alpha_{i_0i_1}\sigma_{i_1\cdots i_{n+1}}+\sigg_{k=1}^{n+1}
(-1)^k\sigma_{i_0\cdots\widehat{i_k}\cdots i_{n+1}} \ = \ 0 \ . \]
For \m{n=1}, this gives that elements of \m{H^1(X,E)} are represented by 
families \m{(\sigma_{ij})},\Nligne \m{\sigma_{ij}\in H^0(U_{ij},\ko_X\ot\C^r)}, 
such that
\[\alpha_{ij}\sigma_{jk}+\sigma_{ij}-\sigma_{ik} \ = \ 0 \ . \]
In \v Cech cohomology, it is generally assumed that 
\m{\theta_{ji}=-\theta_{ij}}. This implies that \ 
\m{\sigma_{ji}=-\alpha_{ji}\sigma_{ij}}.
\end{subsub}

\sepprop

\begin{subsub}\label{cech2}\rm Similarly, an element of \m{H^n(X,E\ot F)} is 
represented by a family  \m{(\mu_{i_0\cdots i_n})}, with \ \m{\mu_{i_0\cdots 
i_n}\in H^0(U_{i_0\cdots i_n},\C^r\ot F)}, satisfying the relations
\[(\alpha_{i_0i_1}\ot I_F)(\mu_{i_1\cdots i_{n+1}})+\sigg_{k=1}^{n+1}
(-1)^k\mu_{i_0\cdots\widehat{i_k}\cdots i_{n+1}} \ = \ 0 \ . \]
The corresponding element of \m{H^0(U_{i_0\cdots i_n},E\ot F)} is \ 
\m{\theta_{i_0\cdots i_n}=(\alpha_{i_0}\ot I_F)^{-1}(\mu_{i_0\cdots i_n})}.

Suppose that \m{n=1}. Let \m{\nu\in H^1(E\ot F)} be the element defined by 
\m{(\theta_{ij})}. For every \m{i\in I}, let \m{\beta_i} be an automorphism 
of \m{\ko_{U_i}\ot\C^r}, and \m{\theta_i\in H^0(U_i,E\ot F)}. The cocycle 
\m{(\theta_{ij}+\theta_i-\theta_j)} represents also $\nu$. Let \ 
\m{\rho_i=(\beta_i\alpha_i\ot I_F)(\theta_i)\in 
H^0(U_i,\C^r\ot F)}. If we replace \m{\alpha_i} with \m{\beta_i\alpha_i}, we 
see that $\nu$ is represented by the family \m{(\mu'_{ij})}, \m{\mu'_{ij}\in 
H^0(U_{ij},\C^r\ot F)}, with
\[\mu'_{ij} \ = \ (\beta_i\ot 
I_F)(\mu_{ij})+\rho_i-(\beta_i\alpha_{ij}\beta_j^{-1}\ot I_F)(\rho_j) \ . \]

\end{subsub}

\sepprop

\begin{subsub}\label{cech2c}Representation of morphisms -- \rm Suppose that we 
have also local trivializations of $F$:
\xmat{\beta_i:F_{|U_i}\ar[r]^-\simeq & \ko_{U_i}\ot\C^s .}
We have then local trivializations of \m{\HHom(E,F)}
\[\xymatrix@R=5pt{\Delta_i:\HHom(E,F)_{|U_i}\ar[r]^-\simeq & 
\ko_{U_i}\ot L(\C^r,\C^s)\\ \phi\fmaps[r] & \beta_i\circ\phi\circ\alpha_i^{-1} 
}\]
such that, for every \ \m{\lambda\in H^0(\ko_{U_{ij}}\ot L(\C^r,\C^s))}, we have
\[\Delta_{ij}(\lambda) \ = \ \Delta_i\Delta_j^{-1}(\lambda) \ = \
\beta_{ij}\lambda\alpha_{ij}^{-1} \ . \]
\end{subsub}

\sepsubsub

\begin{subsub}\label{cech2b}Construction of vector bundles via local 
isomorphisms -- \rm Le $Z$ be a scheme over $\C$, \m{(Z_i)_{i\in I}} an open 
cover of $Z$, and for every \m{i\in I}, a scheme \m{U_i}, with an isomorphism 
\m{\delta_i:Z_i\to U_i}. Let \ 
\m{\delta_{ij}=\delta_j\delta_i^{-1}:\delta_i(Z_{ij})\to\delta_j(Z_{ij})}, 
which is an isomorphism.
Then a vector bundle of $X$ can be constructed (in an obvious way) using vector 
bundles \m{E_i} on \m{U_i} and isomorphisms \ 
\[\Theta_{ij}:\delta_{ij}^*(E_{j|\delta_j(Z_{ij})})\lra 
E_{i|{\delta_i(Z_{ij})}}\]
such that \ \m{\Theta_{ij}\circ\delta_{ij}^*(\Theta_{jk})=\Theta_{ik}} \ on 
\m{\delta_i(Z_{ijk})}.
\end{subsub}

\sepsubsub

\begin{subsub}\label{cech3} Representation of extensions -- \rm Let $E$, $F$ be 
vector bundles on $X$, and\Nligne \m{\sigma\in\Ext^1(F,E)\simeq 
H^1(\Hom(F,E))}. Let
\[0\lra E\lra\ke\lra F\lra 0\]
be the corresponding exact sequence. Suppose that $\sigma$ is represented by a 
cocycle \m{(\sigma_{ij})}, \m{\sigma_{ij}:F_{|U_{ij}}\to E_{|U_{ij}}}. Then 
$\ke$ can be constructed by gluing the vector bundles \m{(E\oplus F)_{|U_i}} 
using the automorphisms of \m{(E\oplus F)_{|U_{ij}}} defined by the matrices
\ \m{\dsp\begin{pmatrix}I_E & \sigma_{ij}\\ 0 & I_F\end{pmatrix}}.
\end{subsub}

\end{sub}

\sepsub

\Ssect{Cohomology of sheaves of non abelian groups}{coh_gr}

(cf. \cite{fr})

Let $X$ be a paracompact topological space and $G$ a sheaf of groups on $X$. We 
denote by \m{H^1(X,G)} the first \v Cech cohomology group of $G$. If $G$ is not 
a sheaf of {\em commutative} groups there is no natural structure of group on 
\m{H^1(X,G)}. Let \m{e\in H^1(X,G)} be defined by the trivial section of $G$ 
(on the trivial open cover \m{(X)} of $X$), we call $e$ the {\em identity 
element} of \m{H^1(X,G)}.

Given an open cover \m{(U_i)_{i\in I}} of $X$, a 1-cocycle of $G$ is a family 
\m{(g_{ij})_{i,j\in I}} such that \m{g_{ij}\in G(U_{ij})} and that the 
relations \m{g_{ij}g_{jk}=g_{ik}} are satisfied.


\sepsubsub

\begin{subsub}\label{coh_gr1} Actions of sheaves of groups -- \rm
Let $H$ be another sheaf of group on $X$, acting on $G$ (the action is 
compatible with the group actions). 
Let \m{z\in H^1(X,H)}, represented by a cocycle \m{(z_{ij})} 
with respect to an open cover \m{(U_i)_{i\in I}} of $X$. Recall that the 
associated sheaf of groups \m{G^z} is defined as follows: we glue the sheaves 
\m{G_{|U_i}} using the automorphisms \m{(z_{ij})} of  \m{G_{|U_{ij}}}. The 
sections of \m{G^z} are represented by families \m{(\gamma_i)_{i\in I}}, 
\m{\gamma_i\in G(U_i)}, such that \m{\gamma_i=z_{ij}\gamma_j} on \m{U_{ij}}. 
The elements of \m{H^1(X,G^z)} are represented by families \m{(\rho_{ij})}, 
\m{\rho_{ij}\in G(U_{ij})}, such that \ 
\m{\rho_{ij}^{-1}\rho_{ik}=z_{ij}\rho_{jk}} \ on \m{U_{ijk}}.
If $G$ is a sheaf of {\em commutative} groups, the elements of the higher 
cohomology groups \m{H^k(X,G^z)}, \m{k\geq 2} can be represented in the same 
way by appropriate families.
\end{subsub}

\sepsubsub

\begin{subsub}\label{coh_gr2}Exact cohomology sequence -- \rm Let 
\m{\Gamma\subset G} be a subsheaf of normal subgroups of $G$, and \m{G/\Gamma} 
the quotient sheaf of groups. Hence we have an exact sequence of sheaves of 
groups
\xmat{0\ar[r] & \Gamma\ar[r]^-i & G\ar[r]^-p & G/\Gamma\ar[r] & 0 . }

Let \m{A_1,\ldots,A_m} be sets, having a particular element \m{e_i\in A_i} for 
\m{1\leq i\leq m}. A sequence of maps
\xmat{A_1\ar[r]^-{f_1} & A_2\ar[r] & \cdots\ar[r] & A_{m-1}
\ar[r]^-{f_{m-1}} & A_m}
is called {\em exact} if for \m{1<i<m} we have \ 
\m{f_i^{-1}(e_i)=f_{i-1}(A_{i-1})} . Then we have: {\em the canonical sequence
\xmat{0\ar[r] & H^0(X,\Gamma)\ar[r]^-{H^0(i)} & H^0(X,G)\ar[r]^-{H^0(p)}
& H^0(X,G/\Gamma)\ar[r]^-\delta & \null\quad\quad\quad\quad}
\xmat{\quad\quad\quad\quad\quad\quad H^1(X,\Gamma)\ar[r]^-{H^1(i)} &
H^1(X,G)\ar[r]^-{H^1(p)} & H^1(X,G/\Gamma)\quad .}
is exact. The subsequence
\xmat{0\ar[r] & H^0(X,\Gamma)\ar[r]^-{H^0(i)} & H^0(X,G)\ar[r]^-{H^0(p)}
& H^0(X,G/\Gamma)}
is an exact sequence of groups. Moreover, if \m{c,c'\in H^0(X,G/\Gamma)}, we 
have \ \m{\delta(c)=\delta(c')} \ if and only if \ 
\m{c{c'}^{-1}\in\imm(H^0(p))}.}

The group \m{H^0(X,G/\Gamma)} acts on \m{H^1(X,\Gamma)} as follows: let \m{a\in 
H^1(X,\Gamma)}, represented by a cocycle \m{(\gamma_{ij})} (with respect to an 
open cover \m{(U_i)} of $X$), and \m{c\in H^0(X,G/\Gamma)}. Suppose also that 
the cover \m{(U_i)} is such that \m{c_{|U_i}} is the image of an element 
\m{c_i} of \m{G(U_i)}. Then \m{(c_i\gamma_{ij}c_j^{-1})} is a cocycle of 
elements of $\Gamma$, and the corresponding element of \m{H^1(X,\Gamma)} 
depends only on $c$ and $a$; we denote it by \m{c.a}. We have then \ 
\m{\delta(c)=c^{-1}.e} \ . Moreover, let \ \m{\gamma,\gamma'\in H^1(X,\Gamma)}. 
Then we have: {\em\m{H^1(i)(\gamma)=H^1(i)(\gamma')} \ if and only if 
\m{\gamma,\gamma'} are in the same \m{H^0(X,G/\Gamma)}-orbit.}
\end{subsub}

\sepsubsub

\begin{subsub}\label{coh_gr3}The fibers of \m{H^1(p)} -- \rm The sheaf of 
groups $G$ acts (by conjugation) on itself and $\Gamma$, \m{G/\Gamma}. Let \ 
\m{\omega\in H^1(X,G/\Gamma)}, \m{g\in H^1(X,G)} be such that \ 
\m{H^1(p)(g)=\omega}. Suppose that $g$ is represented by a cocycle \m{(g_{ij})}.
Let \m{g'\in H^1(p)^{-1}(\omega)}. Then \m{g'} is represented by a cocycle of 
the form \m{(\gamma_{ij}g_{ij})}, with \m{\gamma_{ij}\in\Gamma(U_{ij})}. 
Conversely, a cochain \m{(\gamma_{ij}g_{ij})} is a cocycle if and only if, for 
every \m{i,j,k} we have \ \m{\gamma_{ij}^{-1}\gamma_{ik} = 
g_{ij}\gamma_{jk}g_{ij}^{-1}}, i.e. if and only if \m{(\gamma_{ij})} 
induces a cocycle in \m{\Gamma^g} (cf. \ref{coh_gr1}). In this way we define a 
surjective map
\[\lambda_g:H^1(X,\Gamma^g)\lra H^1(p)^{-1}(\omega)\]
sending the identity element of \m{H^1(X,\Gamma^g)} to $g$.
\end{subsub}
Note that we have an exact sequence \ \m{0\to\Gamma^g\to G^g\to G^g/\Gamma^g\to 
0}.

The group \m{H^0(X,G^g/\Gamma^g)} acts on \m{H^1(X,\Gamma^g)} and the fibers of 
\m{\lambda_g} are the fibers of this action.

\sepsubsub

\begin{subsub}\label{coh_gr4} The case where $\Gamma$ is commutative -- \rm
The sheaf of groups \m{G/\Gamma} acts then by conjugation on $\Gamma$. Let \ 
\m{\omega\in H^1(X,G/\Gamma)}. We define \ \m{\Delta(\omega)\in 
H^2(X,\Gamma^\omega)} \ as follows: we can represent $\omega$ by a cocycle 
\m{(\omega_{ij})} (for a suitable open cover \m{(U_i)} of $X$), such that for 
every \m{i,j} there exists \m{g_{ij}\in G(U_{ij})} over \m{\omega_{ij}}, and 
that \ \m{g_{ji}=g_{ij}^{-1}}. For every indices \m{i,j,k} let \ 
\m{\gamma_{ijk}=g_{ij}g_{jk}g_{ki}}. Then \m{(\gamma_{ijk})} is a family 
representing \m{\Delta(\omega)} (cf. \ref{coh_gr1}). We have then: {\em 
\m{\Delta(\omega)=0} if and only if $\omega$ belongs to the image of 
\m{H^1(p)}}.

The commutativity of $\Gamma$ brings also a supplementary property of the map 
\m{\lambda_g} of \ref{coh_gr3}. The sheaf of groups \m{\Gamma^g} can be seen 
naturally as a subsheaf of commutative groups of of \m{G^g}, and we have a 
canonical isomorphism \ \m{(G/\Gamma)^g\simeq G^g/\Gamma^g}. Hence there is a 
natural action of \m{H^0(X,G^g/\Gamma^g)} on \m{H^1(X,\Gamma^g)}, and we have: 
{\em the fibers of \m{\lambda_g} are the orbits of the action of 
\m{H^0(X,G^g/\Gamma^g)}.} This action is as follows: let \m{h\in 
H^0(G^g/\Gamma^g)}, \m{\gamma\in H^1(X,\Gamma^g)}, represented by 
\m{(\gamma_{ij})}. For a suitable open cover \m{(U_i)}, each \m{h_{|U_i}} can 
be lifted to \m{h_i\in G(U_i)}, and we have \ 
\m{h_ig_{ij}h_jg_{ij}^{-1}\in\Gamma(U_{ij})}. Then \m{h.\gamma} is represented 
by \m{(h_i.\gamma_{ij}.g_{ij}h_j^{-1}g_{ij}^{-1})}.
\end{subsub}

\end{sub}

\sepsub

\Ssect{Extension of polynomials}{ex_pol}

The following results are used in 3.
Let $R$ be a commutative unitary ring and $n$ a positive integer. Let \ 
\m{R_n=R[t]/(t^n)}, which can be seen as the set of polynomials of degree 
\m{<n} in the variable $t$ with coefficients in $R$. In this way we have an 
inclusion \ \m{R_n\subset R_{n+1}}. For every \m{\alpha\in R_n} and every 
integer $p$ such that \m{0\leq p<n}, let \m{\alpha_p} denote the coefficient of 
\m{t^p} in $\alpha$.

Let $r$ be a positive integer, \m{M_{n,r}=M(r,R_n)} \ the $R$-algebra of 
\m{r\times r}-matrices with coefficients if \m{R_n}, and \ 
\m{\GL_{n,r}=\GL(r,R_n)\subset M_{n,r}} \ the group of invertible \m{r\times r} 
matrices. We have inclusions \ \m{M_{n,r}\subset M_{n+1,r}}, 
\m{\GL_{n,r}\subset\GL_{n+1,r}}.

Let \m{A\in M_{n+1,r}}, and $p$ an integer such that \m{0\leq p\leq n}. We 
will denote by \m{A_p\in M(r,R)} the coefficient of the term of degree $p$ of 
$A$, so that \ \m{\dsp A=\sigg_{0\leq p\leq n}A_pt^p}. Let
\[ [A]_n \ = \ \sigg_{0\leq p<n}A_pt^p \ , \]
so that \ \m{A=[A]_n+A_nt^n}.

Let \m{A\in\GL_{n,r}}. An element $B$ of \m{\GL_{n+1,r}} is called an {\em 
extension} of $A$ if \m{[B]_n=A}. Here we will define canonical extensions 
compatible with some actions of automorphisms of \m{R_{n+1}}.

Let $\lambda$ be an 
automorphism of \m{R_{n+1}} such that \ \m{\lambda(\alpha)_0=\alpha} \ for 
every \m{\alpha\in R}. We have \ \m{\lambda(t)=\mu_\lambda t}, with \ 
\m{\mu_\lambda\in R_n}, invertible, and \ 
\m{\mu_{\lambda^{-1}}=\mu_\lambda^{-1}}. We will also denote by $\lambda$ the 
automorphism of \m{R_n} and the bijection \m{M_{n+1,r}\to M_{n+1,r}} induced by 
$\lambda$.

Let \m{A\in\GL_{n+1,r}} (or \m{\GL_{n,r}}). Let \ 
\m{\Inv_\lambda(A)=\lambda^{-1}(A^{-1})}. In this way we define a bijection \ 
\m{\GL_{n+1,r}\to\GL_{n+1,r}} whose inverse is \ \m{A\mapsto 
\Inv_{\lambda^{-1}}(A)}.
Let
\[\Gamma_\lambda(A) \ = \ 
\frac{1}{2}\big(\mu_\lambda^nA_0\Inv_\lambda(A)_nA_0-A_n
\big) \ , \]
so that we have also
\[\Gamma_{\lambda^{-1}}(\Inv_\lambda(A)) \ = \ 
\frac{1}{2}\big(\mu_{\lambda^{-1}}^nA_0^{-1}A_nA_0^{-1}-\Inv_\lambda(A)_n\big) 
\ 
. \]
Let 
\[A_{ext,\lambda} \ = \ A+\Gamma_\lambda(A)t^n \ , \]
which is an extension of $A$ if \m{A\in\GL_{n,r}}.

\sepprop

\begin{subsub}\label{lemX1}{\bf Lemma: }We have \ \rm
\m{\Inv_\lambda(A_{ext,\lambda}) 
\ = \ \Inv_\lambda(A)_{ext,\lambda^{-1}}} .
\end{subsub}

\begin{proof} We have
\begin{eqnarray*}
\Inv_\lambda(A_{ext,\lambda}) & = & 
\lambda^{-1}\big([A+\Gamma_\lambda(A)t^n]^{-1}\big)\\
& = & \lambda^{-1}\big([A(I+A_0^{-1}\Gamma_\lambda(A)t^n)]^{-1})\\
& = & \lambda^{-1}\big((I-A_0^{-1}\Gamma_\lambda(A)t^n))A^{-1})\\
& = & (I-\mu_\lambda^{-n}A_0^{-1}\Gamma_\lambda(A)t^n)\Inv_\lambda(A)\\
& = & \Inv_\lambda(A)-\mu_\lambda^{-n}A_0^{-1}\Gamma_\lambda(A)A_0^{-1}t^n \ .
\end{eqnarray*}
The result follows from the equality \ 
\m{\Gamma_{\lambda^{-1}}(\Inv_\lambda(A))=
-\mu_\lambda^{-n}A_0^{-1}\Gamma_\lambda(A)A_0^{-1}}, which is easily verified.
\end{proof}

\sepprop

\begin{subsub}{\bf Lemma: } Let \m{A\in\GL(n,r)}, viewed as an element of 
\m{\GL(n+1,r)}. Then we have \rm
\[\Inv_\lambda(A_{ext,\lambda}) \ = \ ([\Inv_\lambda(A)]_n)_{ext,\lambda^{-1}} 
\ . 
\]
\end{subsub}

\begin{proof} We have from lemma \ref{lemX1}
\begin{eqnarray*}
\Inv_\lambda(A_{ext,\lambda}) & = & \Inv_\lambda(A)_{ext,\lambda^{-1}}\\
& = & \Inv_\lambda(A)-\frac{1}{2}\Inv_\lambda(A)_nt^n\\
& = & [\Inv_\lambda(A)]_n+\frac{1}{2}\Inv_\lambda(A)_nt^n \ .
\end{eqnarray*}

Let \ \m{B=\Inv_{\lambda^{-1}}([\Inv_\lambda(A)]_n)}, so that \ 
\m{[\Inv_\lambda(A)]_n=B_\lambda}. We have
\[([\Inv_\lambda(A)]_n)_{ext,\lambda^{-1}} \ = \ 
[\Inv_\lambda(A)]_n+\frac{1}{2}\mu_\lambda^{-n}A_0^{-1}B_nA_0^{-1}t^n \ , \]
hence we have to show that \ 
\m{\dsp\Inv_\lambda(A)_n=\mu_\lambda^{-n}A_0^{-1}B_nA_0^{-1}}. We have
\begin{eqnarray*}
B & = & \lambda\big(([\Inv_\lambda(A)]_n)^{-1}\big)\\
& = & \lambda\big((\Inv_\lambda(A)-\Inv_\lambda(A)_nt^n)^{-1}\big)\\
& = & \lambda\big([\Inv_\lambda(A)(1-A_0\Inv_\lambda(A)_nt^n)]^{-1}\big)\\
& = & \lambda\big((1+A_0\Inv_\lambda(A)_nt^n)\Inv_\lambda(A)^{-1}\big)\\
& = & A + \mu_\lambda^nA_0\Inv_\lambda(A)_nA_0t^n \ ,
\end{eqnarray*}
and the result follows immediately since \m{A_n=0}.
\end{proof}

\sepprop

\begin{subsub} Example -- \rm If \m{n=2}, \m{\lambda=I} and \m{A\in\GL_{2,r}}, 
we have
\begin{eqnarray*}
A_{ext,I} & = & A_0+tA_1+\frac{1}{2}t^2A_1A_0^{-1}A_1 \ ,\\ 
(A^{-1})_{ext,I} & = & 
A_0^{-1}-tA_0^{-1}A_1A_0^{-1}+\frac{1}{2}t^2A_0^{-1}A_1A_0^{-1}A_1A_0^{-1} \ . 
\end{eqnarray*}
\end{subsub}

\end{sub}

\sepsec

\section{Canonical class associated to a vector bundle}\label{can_cl}

\Ssect{Definition}{can_cl-def}

Let $X$ be a scheme over $\C$ and $E$ a vector bundle on $X$ of rank $r$. Let 
\m{(U_i)_{i\in I}} be an open cover of $X$ such that we have trivializations:
\xmat{\theta_i:E_{|U_i}\ar[r]^-\simeq & \ko_{U_i}\ot\C^r .}
then 
\xmat{\theta_{ij}=\theta_i\theta_j^{-1}:\ko_{U_{ij}}\ot\C^r\ar[r]^-\simeq &
\ko_{U_{ij}}\ot\C^r}
can be viewed as a $r\times r$ matrix with coefficients in \m{\ko(U_{ij})}. We 
can then define \m{d\theta_{ij}}, which is a $r\times r$ matrix with 
coefficients in \m{\Omega_X(U_{ij})}. Let
\[B_{ij} \ = \ (d\theta_{ij})\theta_{ij}^{-1} \ . \]
We have
\begin{eqnarray*}
B_{ik} & = & (d\theta_{ik})\theta_{ik}^{-1}\\
& = & \big (\theta_{ij}d\theta_{jk}+d\theta_{ij}\theta_{jk}\big)
\theta_{jk}^{-1}\theta_{ij}^{-1}\\
& = & \theta_{ij}(d\theta_{jk})\theta_{jk}^{-1}\theta_{ij}^{-1}+
(d\theta_{ij})\theta_{ij}^{-1}\\
& = & \theta_{ij}B_{jk}\theta_{ij}^{-1}+B_{ij} \ . 
\end{eqnarray*}
From \ref{cech}, it follows that \m{(B_{ij})} represents an element of 
\m{H^1(E\ot E^*\ot\Omega_X)}. It is easy to see that it does not depend on the 
family \m{(\theta_{ij})} defining $E$. We denote by \m{\nabla_0(E)} this 
element, and call it the {\em canonical class of $E$}.

If $L$ is a line bundle on $X$, then we have \ \m{\nabla_0(L)\in H^1(\Omega_X)}.

\end{sub}

\sepsub

\Ssect{Functorial properties}{func_pt}

Let \m{f:X\to Y} be a morphism of schemes over $\C$, and $F$ a vector bundle 
on $Y$. We have a canonical morphism of sheaves \ 
\m{\eta:f^*(\Omega_Y)\to\Omega_X}, hence
\[\eta\ot I:f^*(F\ot F^*\ot\Omega_Y)\lra f^*(F)\ot f^*(F)^*\ot\Omega_X . \]
Let $\iota$ be the canonical map \ \m{H^1(Y, F\ot F^*\ot\Omega_Y)\to 
H^1(X,f^*(F\ot F^*\ot\Omega_Y))}.

\sepprop

\begin{subsub}{\bf Lemma:} We have \ \m{\nabla_0(f^*(F))=H^1(\eta\ot 
I)(\iota(\nabla_0(F)))}.
\end{subsub}

\begin{proof} Let \m{(V_i)_{i\in I}} be an open cover of $Y$ such that for 
every \m{i\in I} there is a trivialization \ 
\m{\theta_i:F_{V_i}\simeq\ko_{V_i}\ot\C^r}, so that, if 
\m{\theta_{ij}=\theta_i\theta_j^{-1}:\ko_{V_{ij}}\ot\C^r\to\ko_{V_{ij}}\ot\C^r},
then \m{(\theta_{ij})} is a cocycle representing $F$. We have then \ 
\m{f^*(\theta_{ij}):\ko_{f^{-1}(V_{ij})}\ot\C^r\to\ko_{f^{-1}(V_{ij})}\ot\C^r}, 
and \m{(f^*(\theta_{ij}))} is a cocycle representing \m{f^*(F)}, with respect 
to the open cover \m{(f^{-1}(V_i)_{i\in I})} of $X$. The result follows 
immediately from the fact that, for every open subset \m{V\subset Y}, and 
\m{\alpha\in\ko_Y(V)} we have \ \m{\eta(d\alpha)=d(f\circ\alpha)}.
\end{proof}

\sepprop

Let \m{E_1,E_2} be vector bundles on $X$. For \m{i=1,2}, let \ 
\m{\Psi_i:\ko_X\to E_i\ot E_i^*} \ be the canonical morphism. If \m{1\leq 
i\not=j\leq 2}, we have then
\[I\ot\Psi_j:E_i\ot E_i^*\ot\Omega_X\lra(E_1\ot E_2)\ot(E_1\ot 
E_2)^*\ot\Omega_X \ . \]

\sepprop

\begin{subsub}{\bf Lemma:} We have \ \m{\nabla_0(E_1\ot 
E_2)=(I\ot\Psi_1)(\nabla_0(E_2))+(I\ot\Psi_2)(\nabla_0(E_1))}.
\end{subsub}

\begin{proof} For \m{k=1,2}, let \m{r_k=\rk(E_k)}. Suppose that \m{E_k} is 
represented by a cocycle \m{(\theta^k_{ij})}, where \m{\theta_{ij}^k} is an 
automorphism of \m{\ko_{U_{ij}}\ot\C^{r_k}}. Then \m{E_1\ot E_2} is 
represented by the cocycle \m{(\theta^1_{ij}\ot\theta^2_{ij})}, 
\m{\theta^1_{ij}\ot\theta^2_{ij}} being an automorphism of \
\m{\ko_{U_{ij}}\ot\C^{r_1}\ot\C^{r_2}}. We have then
\begin{eqnarray*}\big(d(\theta^1_{ij}\ot 
\theta^2_{ij})\big)(\theta^1_{ij}\ot\theta^2_{ij})^{-1} & = & 
\big(d\theta^1_{ij}\ot\theta^2_{ij}+\theta^1_{ij}\ot 
d\theta^2_{ij}\big)(\theta^1_{ij}\ot\theta^2_{ij})^{-1}\\
& = & (d\theta^1_{ij})(\theta^1_{ij})^{-1}\ot 
I+I\ot(d\theta^2_{ij})(\theta^2_{ij})^{-1} \ .
\end{eqnarray*}
The result follows immediately.
\end{proof}

\sepprop

\begin{subsub}{\bf Corollary:} Let \m{L_1,L_2} be line bundles on $X$. Then we 
have \Nligne \m{\nabla_0(L_1\ot L_2)=\nabla_0(L_1)+\nabla_0(L_2)}.
\end{subsub}
\end{sub}

\sepsub

\Ssect{The canonical class of the determinant of a vector bundle}{can_det}

Let $U$ be a scheme over $\C$. Let
\[\xymatrix@R=5pt{\eta_r:\GL(r,H^0(\ko_U))\ar[r] & M(r,H^0(\Omega_U))\\
M\fmaps[r] & M^{-1}dM }\]
We have, for any \ \m{M,N\in\GL(r,H^0(\ko_U))},
\[\eta_r(MN) \ = \ \eta_r(N)+N^{-1}\eta_r(M)N \ . \]
Let
\[\xymatrix@R=5pt{T_r:\GL(r,H^0(\ko_U))\ar[r] & H^0(\Omega_U)\\
M\fmaps[r] & tr(\eta_r(M)) }\]
(where \m{tr} is the trace morphism). We have
\begin{equation}\label{equ1} T_r(MN) \ = \ T_r(M)+T_r(N) \ . \end{equation}

\sepprop

\begin{subsub}\label{prop1}{\bf Proposition: } For every \ 
\m{M\in\GL(r,H^0(\ko_U))} \ we have \ \m{T_r(M)=T_1(\det(M))} .
\end{subsub}

\begin{proof} It is true if \ \m{M\in\GL(r,\C)}, because in this case \ 
\m{T_r(M)=T_1(\det(M)))=0}. 

Let's show that it is also true if $M$ is a triangular matrix (upper or lower). 
Let \Nligne \m{\lambda_1,\ldots,\lambda_r\in\ko^*(U)} be the diagonal elements 
of 
$M$. Then \m{M^{-1}dM} is also triangular, and its diagonal elements are 
\m{\dsp\frac{d\lambda_i}{\lambda_i}}, \m{1\leq i\leq r}. Hence we have
\[T_r(M) \ = \ \sigg_{k=1}^r\frac{d\lambda_k}{\lambda_k} \ . \]
We have \ \m{\det(M)=\prod_{k=1}^r\lambda_k}. Hence
\[d(\det(M)) \ = \ \big(
\prod_{k=1}^r\lambda_k\big)\sigg_{k=1}^r\frac{d\lambda_k}{\lambda_k} \ , \]
and \ \m{T_r(M)=T_1(\det(M)))}.

Note that for any \ \m{A,B\in\ko^*(U)}, we have \ \m{T_1(AB)=T_1(A)+T_1(B)}. 
Hence it follows from ($\ref{equ1}$) that it suffices to show that $M$ is a 
product of matrices that are triangular or in \m{\GL(r,\C)}. We can also 
replace $M$ by the product of $M$ with a triangular matrix, or by the product 
of $M$ with an element of \m{\GL(r,\C)}. Il particular, we can multiply rows or 
columns of $M$ by elements of \m{\ko^*(U)}, add a column (or row) to another 
one, or switch columns or rows.

We now prove proposition \ref{prop1}, more precisely the fact that $M$ is a 
product of matrices that are triangular or in \m{\GL(r,\C)}, by induction on 
$r$. The case \m{r=1} is obvious. Suppose that the theorem is true for $r-1$. 

In the first column of $M$ there is at least one invertible element \m{M_{i1}}. 
By permuting rows of $M$ we can assume that \m{i=1}. For \m{2\leq i\leq r}, by 
subtracting \m{M_{i1}M_{11}^{-1}} times the first row to the $i$-th we can 
assume that \m{M_{i1}=0}. Similarly we can assume that \m{M_{1i}=0}. We have 
then
\[M \ = \ \begin{pmatrix} M_{11} & 0\\ 0 & M'\end{pmatrix} \ , \]
where \ \m{M'\in\GL(r-1,H^0(\ko_U))}. From the induction hypothesis we can write
\[M' \ = \ M'_1\cdots M'_p \ , \]
where, for \m{1\leq i\leq p}, \m{M'_i} is triangular or in \m{\GL(r-1,\C)}. We 
have
\[M \ = \ \begin{pmatrix}M_{11} & 0\\ 0 & I_{r-1}\end{pmatrix}
\begin{pmatrix}1 & 0\\ 0 & M'_1\end{pmatrix}\cdots
\begin{pmatrix}1 & 0\\ 0 & M'_p\end{pmatrix} \ , \]
and \m{\begin{pmatrix}M_{11} & 0\\ 0 & I_{r-1}\end{pmatrix}}, 
\m{\begin{pmatrix}1 & 0\\ 0 & M'_1\end{pmatrix},\ldots,
\begin{pmatrix}1 & 0\\ 0 & M'_p\end{pmatrix}} are triangular or in 
\m{\GL(r,\C)}. 
\end{proof}

\sepprop

Let \ \m{tr_E:H^1(\EEnd(E)\ot\Omega_X)\to H^1(\Omega_X)} \ be the trace 
morphism.

\sepprop

\begin{subsub}{\bf Corollary: } We have \ 
\m{tr_E(\nabla_0(E))=\nabla_0(\det(E))} . \end{subsub}

\sepprop

We have a canonical isomorphism \ \m{\EEnd(E)\simeq\AAd(E)\oplus\ko_X} , where 
\m{\AAd(E)} is the sheaf of endomorphisms of trace 0. Let \m{\nabla(E)} be the 
component of \m{\nabla_0(E)} in \m{H^1(\AAd(E)\ot\Omega_X)}. So we have
\[\nabla_0(E) \ = \ \nabla(E)+\frac{1}{r}\nabla_0(\det(E)) \ . \]
\end{sub}

\sepsec

\section{Primitive multiple schemes}\label{PMV}

\Ssect{Definitions}{PMV_def}

Let $X$ be a smooth connected variety, and \ \m{d=\dim(X)}. A {\em multiple 
scheme with support $X$} is a Cohen-Macaulay scheme $Y$ such that 
\m{Y_{red}=X}. If $Y$ is quasi-projective we say that it is a {\em multiple 
variety with support $X$}. In this case $Y$ is projective if $X$ is.

Let $n$ be the smallest integer such that \m{Y=X^{(n-1)}}, \m{X^{(k-1)}}
being the $k$-th infinitesimal neighborhood of $X$, i.e. \
\m{\ki_{X^{(k-1)}}=\ki_X^{k}} . We have a filtration \ \m{X=X_1\subset
X_2\subset\cdots\subset X_{n}=Y} \ where $X_i$ is the biggest Cohen-Macaulay
subscheme contained in \m{Y\cap X^{(i-1)}}. We call $n$ the {\em multiplicity}
of $Y$.

We say that $Y$ is {\em primitive} if, for every closed point $x$ of $X$,
there exists a smooth variety $S$ of dimension \m{d+1}, containing a 
neighborhood of $x$ in $Y$ as a locally closed subvariety. In this case, 
\m{L=\ki_X/\ki_{X_2}} is a line bundle on $X$, \m{X_j} is a primitive multiple 
scheme of multiplicity $j$ and we have \ 
\m{\ki_{X_j}=\ki_X^j}, \m{\ki_{X_{j}}/\ki_{X_{j+1}}=L^j} \ for \m{1\leq j<n}. 
We call $L$ the line bundle on $X$ {\em associated} to $Y$. The ideal sheaf 
\m{\ki_{X,Y}} can be viewed as a line bundle on \m{X_{n-1}}.

Let \m{P\in X}. 
Then there exist elements \m{y_1,\ldots,y_d}, $t$ of \m{m_{S,P}} whose images 
in \m{m_{S,P}/m_{S,P}^2} form a basis, and such that for \m{1\leq i<n} we have 
\ \m{\ki_{X_i,P}=(t^{i})}. In this case the images of \m{y_1,\ldots,y_d} in 
\m{m_{X,P}/m_{X,P}^2} form a basis of this vector space.

A {\em multiple scheme with support $X$} is primitive if and only if 
\m{\ki_X/\ki_X^2} is zero or a line bundle on $X$ (cf. \cite{dr7}, proposition 
2.3.1).  

Even if $X$ is projective, we do not assume that $Y$ is projective. In fact we 
will give examples of non quasi-projective $Y$.

The simplest case is when $Y$ is contained in a smooth variety $S$ of dimension 
\m{d+1}. Suppose that $Y$ has multiplicity $n$. Let \m{P\in X} and 
\m{f\in\ko_{S,P}}  a local equation of $X$. Then we have \ 
\m{\ki_{X_i,P}=(f^{i})} \ for \m{1<j\leq n} in $S$, in particular 
\m{\ki_{Y,P}=(f^n)}, and \ \m{L=\ko_X(-X)} .


For any \m{L\in\Pic(X)}, the {\em trivial primitive variety} of multiplicity 
$n$, with induced smooth variety $X$ and associated line bundle $L$ on $X$ is 
the $n$-th infinitesimal neighborhood of $X$, embedded by the zero section in 
the dual bundle $L^*$, seen as a smooth variety.

\end{sub}

\sepsub

\Ssect{Construction and parametrization of primitive multiple 
schemes}{PMV_cons}

Let $Y$ be a primitive multiple scheme of multiplicity $n$, \m{X=Y_{red}}.
Let \ \m{{\bf Z}_n=\spec(\C[t]/(t^n))}.
Then for every closed point \m{P\in X}, there exists an open neighborhood $U$ 
of $P$ in $X$, such that if \m{U^{(n)}} is the corresponding neighborhood of 
$P$ in $Y$, there exists a commutative diagram
 \xmat{ & U\flinc[ld]\flinc[rd] \\
U^{(n)}\ar[rr]^-\simeq & & U\times {\mathbf Z}_n}
i.e. $Y$ is locally trivial (\cite{dr1}, th\'eor\`eme 5.2.1, corollaire 5.2.2).

It follows that we can construct a primitive multiple scheme of multiplicity 
$n$ by taking an open cover \m{(U_i)_{i\in I}} of $X$ and gluing the varieties 
\ \m{U_i\times{\bf Z}_n} (with automorphisms of the \ \m{U_{ij}\times{\bf Z}_n} 
\ leaving \m{U_{ij}} invariant).

Let \ \m{\ka_n=\ko_{X\times{\bf Z}_n}}. So we have, for any open subset $U$ of 
$X$
\[\ka_n(U) \ = \ \ko_X(U)\ot_\C\C[t]/(t^n) \ . \]
We have a canonical morphism \ \m{\ka_n\to\ko_X} \ of sheaves on $X$ coming 
from the inclusion \ \m{X\subset X\times{\bf Z}_n}.
Let \m{\kg_n} be the sheaf of groups of automorphisms of \m{\ka_n} leaving $X$ 
invariant, i.e. \m{\kg_n(U)} is the group of automorphisms of $\C$-algebras \ 
\m{\theta:\ka_n(U)\to\ka_n(U)} \ such that the following diagram is commutative
\xmat{
\ka_n(U)\ar[rr]^-\theta\ar[rd] & & \ka_n(U)\ar[ld]\\
& \ko_X(U)
}
Then the set of primitive multiple schemes $Y$ of multiplicity $n$, such that 
\m{X=Y_{red}}, can be identified with the cohomology set \m{H^1(X,\kg_n)}. 

\sepsubsub

\begin{subsub}\label{PMS_coh} Primitive multiple schemes and the associated 
cohomology classes -- \rm For every open subset $U$ of $X$, let \m{U^{(n)}} be 
the corresponding open subset of $Y$.
Let \m{(U_i)_{i\in I}} be an affine open cover of $X$ such that we have 
trivializations
\xmat{\delta_i:U_i^{(n)}\ar[r]^-\simeq & U_i\times{\bf Z}_n , }
and \ \m{\delta_i^*:\ko_{U_i\times{\bf Z}_n}\to\ko_{U_i^{(n)}}} \ the 
corresponding isomorphism. Let
\xmat{\delta_{ij}=\delta_j\delta_i^{-1}:U_{ij}\times{\bf Z}_n\ar[r]^-\simeq & 
U_{ij}\times{\bf Z}_n , }
and \ \m{\delta_{ij}^*=\delta_i^{*-1}\delta_j^*\in\kg_n(U_{ij})}. 
Then \m{(\delta_{ij}^*)} is a cocycle which represents the element \m{g_n} of 
\m{H^1(X,\kg_n)} corresponding to $Y$. 
\end{subsub}

\sepsubsub

\begin{subsub}\label{I_X} The ideal sheaf of $X$ -- \rm
There exists \ \m{a_{ij}\in\ko_X(U_{ij})\ot_\C\C[t]/(t^{n-1})} \ such that \ 
\m{\delta_{ij}^*(t)=a_{ij}t}. Let \ \m{\alpha_{ij}=a_{ij|X}\in\ko_X(U_i)}.
For every \m{i\in I}, \m{\delta_i^*(t)} is a generator of 
\m{\ki_{X,Y|{U^{(n)}}}}. So we have local trivializations
\[\xymatrix@R=5pt{\lambda_i:\ki_{X,Y|{U_i^{(n-1)}}}\ar[r] & 
\ko_{U_i^{(n-1)}}\\
\delta_i^*(t)\fmaps[r] & 1}\]
Hence \ \m{\lambda_{ij}=\lambda_i\lambda_j^{-1}: 
\ko_{U_{ij}^{(n-1)}}\to\ko_{U_{ij}^{(n-1)}}} \ is the multiplication by 
\m{\delta_j^*(a_{ij})}. It follows that 
\m{(\delta_j^*(a_{ij}))} (resp. 
\m{(\alpha_{ij})})  is a cocycle representing the line bundle 
\m{\ki_{X,Y}} (resp. \m{L}) on \m{X_{n-1}} (resp. $X$). 

We have a canonical morphism of sheaves of groups
\[\xi_n:\kg_n\lra\ko_X^*\]
defined as follows: if $U$ is an open subset of $X$ and \m{\phi\in\kg_n(U)}, 
then 
there exists \m{\nu\in\ko_X^*(U)} such that \ \m{\phi(t)=\nu t}. Then \ 
\m{\xi_n(\phi)=\nu}. The map
\[H^1(\xi_n):H^1(X,\kg_n)\lra H^1(X,\ko_X^*)\]
associates to the primitive scheme $Y$ the associated line bundle $L$ on $X$.
This map is surjective: for every line bundle $L$ on $X$, there exists a
primitive variety of multiplicity $n$ and associated line bundle L (the trivial 
primitive variety).
\end{subsub}

\sepprop

\begin{subsub}\label{const_sh} Descriptions using the open cover \m{(U_i)} \rm

{\bf (i)} \ {\em Construction of sheaves -- } Let $\ke$ be a 
coherent sheaf on \m{X_n}. We can define it in the usual way, starting with 
sheaves \m{\kf_i} on the open sets \m{U_i^{(n)}} and gluing them. We take 
these sheaves of the form \ \m{\kf_i=\delta_i^*(\ke_i)}, where \m{\ke_i} is a 
sheaf on \ \m{U_i\times {\bf Z}_n}. To glue the sheaves \m{\kf_i} on the 
intersection 
\m{U_{ij}^{(n)}} we use isomorphisms \ 
\m{\rho_{ij}:\kf_{j|U_{ij}^{(n)}}\to\kf_{i|U_{ij}^{(n)}}}, with the relations \ 
\m{\rho_{ik}=\rho_{ij}\rho_{jk}}. Let
\[\theta_{ij} = 
(\delta_i^*)^{-1}(\rho_{ij}):\delta_{ij}^*(\ke_{j|U_{ij}\times{\bf 
Z}_n})\lra\ke_{i|U_{ij}\times{\bf Z}_n} \ .
\]
We have then the relations \ 
\m{\theta_{ik}=\theta_{ij}\circ\delta_{ij}^*(\theta_{jk})}. Conversely, 
starting with sheaves \m{\ke_i} and isomorphisms \m{\theta_{ij}} satisfying the 
preceding relations, one obtains a coherent sheaf on \m{X_n}.

This applies to trivializations, i.e when \ \m{\ke_i=\ko_{U_i^{(n)}}\ot\C^r}. 
We have then \Nligne \m{\theta_{ij}:\ko_{U_{ij}\times{\bf 
Z}_n}\ot\C^r\to\ko_{U_{ij}\times{\bf Z}_n}\ot\C^r}.

\sepsubsub

{\bf (ii)} \ {\em Morphisms -- }
Suppose that we have another sheaf \m{\ke'} on \m{X_n}, defined by sheaves 
\m{\ke'_i} on \ \m{U_i\times {\bf Z}_n} \ and isomorphisms \m{\theta'_{ij}}. 
One can see easily that a morphism \ \m{\Psi:\ke\to\ke'} \ is defined by 
morphisms \ \m{\Psi_i:\ke_i\to\ke'_i} \ such that \ 
\m{\theta'_{ij}\circ\delta_{ij}^*(\Psi_j)=\Psi_i\circ\theta_{ij}}.

\sepsubsub

{\bf (iii)} \ {\em Inverse images -- } Let $\varphi$ be an automorphism of 
\m{X_n} 
inducing 
the identity on $X$. Let \Nligne 
\m{\varphi_i=\delta_i\circ\varphi\circ\delta_i^{-1}:
U_i\times{\bf Z}_n\to U_i\times{\bf Z}_n}. We have then \ 
\m{\varphi_j\circ\delta_{ij}=\delta_{ij}\circ{\varphi_i}}. And conversely, 
given automorphisms \ \m{\varphi_i:U_i\times{\bf Z}_n\to U_i\times{\bf Z}_n} \ 
inducing the identity on \m{U_i} and satisfying the preceding relations we can 
build a corresponding automorphism $\varphi$ of \m{X_n}.
 
Suppose that $\ke$ is defined by trivializations as in (i), and described with 
the morphisms \m{\theta_{ij}}. In the same way \m{\varphi^*(\ke)} is described 
by the morphisms \m{\varphi_i^*(\theta_{ij})}.
\end{subsub}

\end{sub}

\sepsub

\Ssect{The case of double schemes}{g2}

We suppose that \m{n=2}. Then there exists a map \
\m{D_{ij}:\ko_X(U_{ij})\to\ko_X(U_{ij})} \ such that, for every 
\m{f\in\ko_X(U_{ij})} we have
\[\delta_{ij}^*(f) \ = \ f+tD_{ij}(f) \ , \]
and it is easy to see that \m{D_{ij}} is a derivation of \m{\ko_X(U_{ij})}, 
i.e. a section of \m{T_{X|U_{ij}}}. It follows that \m{\delta_{ij}^*} can be 
represented as the matrix \m{\begin{pmatrix}1 & 0\\ D_{ij} & 
\alpha_{ij}\end{pmatrix}}. The 
formula \ \m{\delta_{ij}^*\delta_{jk}^*=\delta_{ik}^*} \ implies that \ 
\m{D_{ij}+\alpha_{ij}D_{jk}=D_{ik}}. It follows from \ref{cech2} that the 
family \m{(D_{ij})} represents an element of \m{H^1(T_X\ot L)}.
This element does not depend on the choice of the automorphisms \m{\delta_i} 
and the trivializations \m{\lambda_{i|X}} of $L$ : suppose that \m{\tau_i} is 
an automorphism of \m{U_i\times{\bf Z}_n}, represented by the matrix 
\m{\begin{pmatrix}1 & 0\\ D_i & \beta_i\end{pmatrix}}, and that we replace 
\m{\delta_i} with \m{\delta'_i=\tau_i\delta_i}. Then \m{\delta_{ij}^*} is 
replaced with \m{\delta_{ij}^{*'}}, represented by the matrix
\[\begin{pmatrix}1 & 0\\ D'_{ij} & \alpha_{ij}\beta_j\beta_i^{-1}\end{pmatrix}
\ = \ \begin{pmatrix}1 & 0\\ D_i & \beta_i\end{pmatrix}^{-1}\begin{pmatrix}1 & 
0\\ D_{ij} & \alpha_{ij}\end{pmatrix}\begin{pmatrix}1 & 0\\ D_j & 
\beta_j\end{pmatrix} \ , \]
with \ 
\m{D'_{ij}=-\beta_i^{-1}D_i+\beta_i^{-1}D_{ij}+\beta_i^{-1}\alpha_{ij}D_j}. It 
follows from \ref{cech2} that \m{(D'_{ij})} represents the same element of 
\m{H^1(T_X\ot L)} as \m{(D_{ij})}.

We have \ \m{\ker(\xi_2)\simeq T_X}, so we have an exact sequence of sheaves of 
groups on $X$
\xmat{0\ar[r] & T_X\ar[r] & \kg_2\ar[r]^-{\xi_2} & \ko_X^*\ar[r] & 0 .}
Recall that \m{g_2\in H^1(X,\kg_2)} is associated to $Y$. Let 
\m{g_1=H^1(\xi_2)(g_2)}.

\sepprop

\begin{subsub}{\bf Lemma: }\label{g2_lem1} {\bf 1 -- } We have \ 
\m{T_X^{g_1}\simeq T_X\ot L}.

{\bf 2 -- } The action of \m{H^0(X,\ko_X^*)=\C^*} on \m{H^1(X,T_X\ot L)} is the 
multiplication.
\end{subsub}
\begin{proof} Let \m{U\subset X} be an open subset, \m{D,D_0\in H^0(U,T_X)}, 
\m{\alpha\in H^0(U,\ko_X^*)}. The result follows easily from the formula
\[\begin{pmatrix}1 & 0\\D & \alpha\end{pmatrix}
\begin{pmatrix}1 & 0\\D_0 & 1\end{pmatrix}
\begin{pmatrix}1 & 0\\D & \alpha\end{pmatrix}^{-1} \ = \
\begin{pmatrix}1 & 0\\\alpha D_0 & 1\end{pmatrix}\]
and \ref{cech1}.
\end{proof}

\sepprop

\begin{subsub}{\bf Proposition: }\label{g2_lem2} The element of \ 
\m{H^1(X,T_X\ot 
L)=\Ext^1_{\ko_X}(\Omega_X,L)} \ corresponding to the canonical exact sequence
\[0\lra L\lra\Omega_{Y|X}\lra\Omega_X\lra 0\]
is \m{g_2}.
\end{subsub}
\begin{proof} The vector bundle \m{\Omega_{Y|X}} can be constructed, by the 
method of \ref{cech2b}, using the local isomorphisms \m{U_i^{(2)}\simeq 
U_i\times{\bf Z}_2}, the bundles \m{\Omega_{U_i\times{\bf Z}_2|U_i}}, and the 
canonical automorphisms of \m{\Omega_{U_{ij}\times{\bf Z}_2|U_{ii}}} defined by 
\m{\delta_{ij}}. Let
\[\mu_{ij}:\delta_{ij}^*(\Omega_{U_{ij}\times{\bf 
Z}_2})=\Omega_{U_{ij}\times{\bf Z}_2}
\lra\Omega_{U_{ij}\times{\bf Z}_2}\]
be the canonical morphism. Then for every \m{\alpha\in\ko_X(U_{ij})} we have
\begin{eqnarray*}\mu_{ij}(d\alpha) & = & d(\delta_{ij}\circ\alpha)\\
& = & d(\alpha+D_{ij}(\alpha)t)\\
& = & d\alpha+d(D_{ij}(\alpha))t+D_{ij}(\alpha)dt ,\\
\mu_{ij}(dt) & = & d(\delta_{ij}\circ t)\\
& = & d(\alpha_{ij}t)\\ & = & d(\alpha_{ij})t+\alpha_{ij}dt .
\end{eqnarray*}
It follows that \ \m{\mu_{ij|U_{ij}}:\Omega_{U_{ij}\times{\bf Z}_2|U_{ij}}\to
\Omega_{U_{ij}\times{\bf Z}_2|U_{ij}}} \ is defined, with respect to the 
isomorphism \ \m{\Omega_{U_{ij}\times{\bf 
Z}_2|U_{ij}}\simeq\Omega_{U_{ij}}\oplus(\ko_{U_{ij}}\ot\C dt)}, by the matrix
\m{\dsp\begin{pmatrix}I & 0\\ D_{ij} & \alpha_{ij}\end{pmatrix}}. The result 
follows then easily from \ref{cech}.
\end{proof}

\sepprop

\begin{subsub}\label{param_2} Parametrization of primitive double schemes -- \rm
From lemma \ref{g2_lem1} it follows that the non trivial primitive double 
schemes $Y$ such that \m{Y_{red}=X} and with associated line bundle $L$ are 
parametrized by \m{\P(H^1(X,T_X\ot L))}. This result has been proved by another 
method in \cite{ba_ei} (cf. \ref{b_e}).
\end{subsub}

\end{sub}

\sepsub

\Ssect{The case of primitive schemes of multiplicity $n>2$}{gn}

There is a canonical obvious morphism \ \m{\rho_n:\kg_n\to\kg_{n-1}}.

\sepprop

\begin{subsub}{\bf Proposition: }\label{gn_lem1} We have \ 
\m{\ker(\rho_n)\simeq T_X\oplus\ko_X}.
\end{subsub}
\begin{proof} Let \m{U\subset X} be an open subset. Let $\phi$ be an 
automorphism of \ \m{U\times{\bf Z}_n} \ leaving $U$ invariant, and \m{\phi^*} 
the associated automorphism of \m{\ko_X(U)[t]/(t^n)}. It belongs to 
\m{\ker(\rho_n)(U)} if and only if there is a derivation $D$ of \m{\ko_X(U)} 
and \m{\alpha\in\ko_X^*(U)} such that \ 
\m{\phi^*(\lambda)=\lambda+D(\lambda)t^{n-1}} \ for every 
\m{\lambda\in\ko_X(U)}, and \ \m{\phi^*(t)=(1+\alpha t^{t-2})t}. The result 
follows easily.
\end{proof}

\sepprop

Let \m{g_{n}\in H^1(X,\kg_{n})}, and for \m{2\leq k<n}, 
\m{g_k=\rho_{k+1}\circ\cdots\circ\rho_{n}(g_{n})}. We denote by $Y$ the 
double primitive scheme defined by \m{g_2}, and use the notations of 
\ref{PMV_cons} and \ref{g2}. Suppose that \m{g_{n}} is represented by the 
cocycle \m{(\delta_{ij})} with respect to \m{(U_i)}.

\sepprop

\begin{subsub}{\bf Proposition: }\label{gn_lem2} We have \ 
\m{\ker(\rho_n)^{g_{n-1}}\simeq(\Omega_{Y|X})^*\ot L^{n-1}}.
\end{subsub}
\begin{proof} The sheaf \m{\ker(\rho_n)^{g_{n-1}}=(T_X\oplus\ko_X)^{g_{n-1}}} 
is constructed as follows: we glue the sheaves \m{(T_X\oplus\ko_X)_{|U_i}} 
using the automorphisms of \m{\ker(\rho_n)_{|U_{ij}}}: 
\m{\phi\mapsto\delta_{ij}^*\phi(\delta_{ij}^*)^{-1}}. We can write
\[\phi(\alpha) \ = \ \alpha+D(\alpha)t^{n-1} \ , \qquad\qquad \phi(t) \ = \ 
(1+\epsilon t^{n-2})t\]
for every \m{\alpha\in\ko_X(U_{ij})} (for some derivation $D$ of 
\m{\ko_X(U_{ij})} and \m{\epsilon\in\ko_X(U_{ij})}); $D$ and $\epsilon$ are the 
components of $\psi$ if \m{H^0(U_{ij},T_X)} and \m{\ko_X(U_{ij})} respectively.
An easy computation shows then that
\begin{eqnarray*}\delta_{ij}^*\psi(\delta_{ij}^*)^{-1}(\alpha) & = & 
\alpha+\big(D(\alpha)+D_{ij}(\alpha)\epsilon\big)\alpha_{ij}^{n-1}t^{n-1} , 
\\ \delta_{ij}^*\psi(\delta_{ij}^*)^{-1}(t) & = & 
(1+\epsilon\alpha_{ij}^{n-2}t^{n-2})t , \end{eqnarray*}
i.e. we obtain the automorphism of \m{(T_X\oplus\ko_X)_{|U_{ij}}} defined by 
the matrix \Nligne \m{\dsp\begin{pmatrix}\alpha_{ij}^{n-1} & 
\alpha_{ij}^{n-1}D_{ij}\\ 0 & \alpha_{ij}^{n-2}\end{pmatrix}}. The result 
follows then from \ref{cech}.
\end{proof}

\end{sub}

\sepsub

\Ssect{Primitive multiple schemes with extensions of the canonical ideal 
sheaf}{prim_can}

We use the notations of \ref{PMV_def} to \ref{gn}. Suppose that \m{n\geq 2}, 
and let \m{Y=X_n} be a primitive multiple scheme of multiplicity $n$ such that
\m{Y_{red}=X}, with associated line bundle \m{L\in\Pic(X)}. The ideal sheaf 
\m{\ki_{X,Y}} is a line bundle on \m{X_{n-1}}. A necessary condition for the 
possibility to extend $Y$ to a primitive multiple scheme \m{X_{n+1}} of 
multiplicity \m{n+1} is that \m{\ki_{X,Y}} can be extended to a line bundle on 
$Y$, because in this case \m{\ki_{X,X_{n+1}}} is a line bundle on $Y$ and we 
have \ \m{\ki_{X,X_{n+1}|Y}=\ki_{X,Y}}. This is why we will consider pairs 
\m{(Y,\L)}, where $\L$ is a line bundle on $Y$ such that \ 
\m{\L_{|X_{n-1}}\simeq\ki_{Y,X}}. 

\sepprop

\begin{subsub} The corresponding sheaf of groups -- \rm  The sheaf of groups 
\m{\kh_n} on $X$ (as \m{\kg_n} for primitive multiple schemes of multiplicity 
$n$) is defined as follows: for every open subset \m{U\subset X}, \m{\kh_n(U)} 
consists of pairs \m{(\phi,u)}, where \m{\phi\in\kg_n(U)}, and 
\m{u\in\ko_X(U)[t]/(t^n)} is such that \ \m{\phi(t)=ut}. Note that $u$ is then 
defined up to a multiple of \m{t^{n-1}}. The multiplication is defined by:
\[(\phi',u')(\phi,u) \ = \ (\phi'\phi,u'\phi'(u))\]
(cf. the formulas in \ref{str_aut}). Associativity is easily verified. The 
identity element is \m{(I,1)} and the inverse of \m{(\phi,u)} is 
\m{\dsp(\phi^{-1},\phi^{-1}(\frac{1}{u}))}. The set of isomorphism classes of 
the above pairs \m{(Y,\L)} can then be identified with the cohomology set 
\m{H^1(X,\kh_n)}. More precisely, given a cocycle \m{((\phi_{ij},u_{ij}))} of 
\m{\kh_n} associated to \m{(Y,\L)}, \m{(\phi_{ij})} is a cocycle of \m{\kg_n} 
defining $Y$ and \m{(u_{ij})} a family defining \m{\L} (according to 
\ref{I_X} and \ref{const_sh}).

We have canonical obvious morphisms
\[\epsilon_n:\kh_n\lra\kg_n \ , \qquad \tau_n:\kg_{n+1}\lra\kh_n \ . \]

\sepprop

\begin{subsub}\label{prop5}{\bf Proposition: } {\bf 1 -- } We have \ 
\m{\ker(\tau_n)\simeq T_X}.

{\bf 2 -- } Let \m{g\in H^1(X,\kh_{n})}. Then we have \ \m{T_X^g\simeq T_X\ot 
L^n}.
\end{subsub}

The proof is similar to those of propositions \ref{gn_lem1} and \ref{gn_lem2}.

\sepprop

Finally we obtain commutative diagrams with exact rows

\begin{equation}\label{equ2}
\begin{split}
\xymatrix{0\ar[r] & T_X\ar[r]\flinc[d] & \kg_{n+1}\ar[rr]^-{\tau_n}\fleq[d] & &
\kh_n\ar[r]\ar[d]^{\epsilon_n} & 0\\
0\ar[r] & T_X\oplus\ko_X\ar[r] & \kg_{n+1}\ar[rr]^-{\rho_{n+1}} & & \kg_n\ar[r] 
& 0}\\
\xymatrix{0\ar[r] & T_X\ot L^n\ar[r]\ar[d]^\iota & \kg_{n+1}^g\ar[r]\fleq[d] &
\kh_n^g\ar[r]\ar[d] & 0\\
0\ar[r] & (\Omega_{Y|X})^*\ot L^n\ar[r] & \kg_{n+1}^g\ar[r] & \kg_n^g\ar[r] & 0
}
\end{split}
\end{equation}

It is easy to see that $\iota$ is the injective morphism deduced from the exact 
sequence of proposition \ref{gn_lem2}.
\end{subsub}

\sepprop

\begin{subsub} The associated sheaf of groups -- \rm Let \m{g\in 
H^1(X,\kh_{n})}, corresponding to the pair \m{(Y,\L)}. 
Suppose that $g$ is defined by the cocycle \m{((\delta_{ij}^*,u_{ij}))}, 
\m{(\delta_{ij}^*,u_{ij})\in\kh_n(U_{ij})}, so that \m{(\delta_{ij}^*)} is the 
cocycle defining $Y$ (cf. \ref{PMS_coh}).

A global section of \m{\kh_n^g} is defined by a family \m{((\chi_i,v_i))_{i\in 
I}}, \m{(\chi_i,v_i)\in\kh_n(U_i)}, such that
\[(\delta_{ij}^*,u_{ij})(\chi_j,v_j)(\delta_{ij}^*,u_{ij})^{-1} \ = \ 
(\chi_i,v_i) \ , \]
which is equivalent to the two relations
\begin{eqnarray*}
\delta_{ij}^*\chi_j\delta_{ij}^{*-1} & = & \chi_i \ ,\\
u_{ij}\delta_{ij}^*(v_j) & = & v_i\chi_i(u_{ij}) \ .
\end{eqnarray*}
The first one says that \m{(\chi_i)} defines an automorphism $\chi$ of $Y$ 
inducing the identity on $X$. The second that \m{(v_i)} defines an isomorphism 
\ \m{\L\to\chi^*(\L)}. It follows that {\em we can identify 
\m{H^0(X,\kh_n^g)} with the set of pairs \m{(\chi,\eta)}, where $\chi$ is an 
automorphism of $Y$ and \ \m{\eta:\L\to\chi^*(\L)} \ is an isomorphism.}
\end{subsub}

\end{sub}

\sepsub

\Ssect{Ample line bundles and projectivity}{amp_li}

We use the notations of \ref{PMV_def} to \ref{gn}. Suppose that \m{n\geq 2}, 
and let \m{Y=X_n} be a primitive multiple scheme of multiplicity $n$ such that
\m{Y_{red}=X}, with associated line bundle \m{L\in\Pic(X)}. The following 
result follows from prop. 4.2 of \cite{ha2}:

\sepprop

\begin{subsub}\label{prop10}{\bf Proposition: } Let $\D$ be a line bundle on 
\m{X_n}, and \m{D=\D_{|X}}. Then $\D$ is ample if and only $D$ is ample.
\end{subsub}

\sepprop

\begin{subsub}\label{theo1}{\bf Corollary:} The scheme \m{X_n} is 
quasi-projective if and only if there exists an ample line bundle on $X$ that 
can be extended to a line bundle on \m{X_n}.
\end{subsub}

In this case \m{X_n} is even projective.

\end{sub}

\sepsec

\section{Structure of $\kg_n$}\label{struc}

We use the notations of \ref{PMV_def} and \ref{PMV_cons}. 

\sepsub

\Ssect{Description of local automorphisms}{str_aut}

Let \m{U\subset X} be a non empty open subset. Let $\phi$ be an automorphism of 
\ \m{U\times {\bf Z}_n} \ leaving $U$ invariant, that we can view $\phi$ as an 
automorphism of \m{\ko_X(U)[t]/(t^n)} such that for every \m{\alpha\in\ko_X(U)} 
we have \ \m{\phi(\alpha)=\alpha+t\eta(\alpha)}, for some map \ 
\m{\eta:\ko_X(U)\to\ko_X(U)/(t^{n-1})}. We have also \ \m{\phi(t)=\mu t}, 
for some \ \m{\mu\in\ko_X(U)/(t^{n-1})}, invertible, and $\phi$ is completely 
defined by $\eta$ and $\mu$. The map $\eta$ is $\C$-linear and has the 
following property: for every \m{\alpha,\beta\in\ko_X(U)} we have
\[\eta(\alpha\beta) \ = \ 
\eta(\alpha)\beta+\alpha\eta(\beta)+\eta(\alpha)\eta(\beta).t \ . \]
In particular, if \m{n=2} then $\eta$ is a derivation of \m{\ko_X(U)}.
Conversely, given an invertible element \ \m{\mu\in\ko_X(U)/(t^{n-1})} \ and a 
$\C$-linear map $\eta$ satisfying the preceding property, there is a unique 
automorphism $\phi$ of \m{\ko_X(U)[t]/(t^n)} such that for every 
\m{\alpha\in\ko_X(U)} we have \ \m{\phi(\alpha)=\alpha+t\eta(\alpha)}, and \ 
\m{\phi(t)=\mu t}. We will write
\[\phi \ = \ \phi_{\eta,\mu} \ . \]

Let \m{d=\dim(X)}. Suppose that there exist \m{x_1,\ldots,x_d\in\ko_X(U)} such 
that for every \m{P\in U}, \m{x_1-x_1(P),\ldots,x_d-x_d(P)} generate 
\m{m_{X,P}/m_{X,P}^2}. Then $\eta$ is completely determined by $\mu$ and \ 
\m{\eta(x_i)}, \m{1\leq i\leq d}. Conversely, any sequence 
\m{(\lambda_1,\ldots,\lambda_d,\mu)} of elements of 
\m{\mu\in\ko_X(U)/(t^{n-1})} uniquely defines such a map $\eta$ with the 
required properties. For every \m{\alpha\in\ko_X(U)}, 
\m{\phi_{\eta,\mu}(\alpha)} can be computed with the Taylor formula. For 
example, for \m{n=3}
\[\phi_{\eta,\mu}(\alpha) \ = \ 
\alpha+t\sigg_{i=1}^d\lambda_i\frac{\partial\alpha}{\partial x_i}+
\frac{t^2}{2}\sigg_{0\leq i,j\leq 
d}\lambda_i\lambda_j\frac{\partial^2\alpha}{\partial x_i\partial x_j} \ . \]

For every \m{\beta\in\ko_X(U)[t]/(t^n)} and \m{0\leq k<n}, let \m{\beta_k} be 
the coefficient of \m{t^k} in $\beta$, and if \m{k>0}
\[ [\beta]_k \ = \ \sigg_{p=0}^{k-1}\beta_pt^p \ . \]

Let \m{\phi_{\eta,\mu}} and \m{\phi_{\eta',\mu'}} be automorphisms, and
\[\phi_{\eta'',\mu''} \ = \ \phi_{\eta',\mu'}\circ\phi_{\eta,\mu} \ . \]
Then we have, for every \m{\alpha\in\ko_X(U)}
\[\eta''(\alpha) \ = \ 
\eta'(\alpha)+\nu'[\phi_{\eta',\nu'}(\eta(\alpha))]_{n-1} \ , \]
and \ \m{\nu''=[\phi_{\eta',\nu'}(\nu)]_{n-1}\nu'} \ .
\end{sub}

\sepsub

\Ssect{The case $n=2$}{nequ2} 

Let \ \m{\phi_{D,\mu}} \ be an 
automorphism of \m{\ko_X(U)[t]/(t^2)}, where $D$ is a derivation of 
\m{\ko_X(U)} and \m{\mu\in\ko_X(U)} is invertible.

\sepsubsub

\begin{subsub} Composition and inverse -- \rm Let \m{\Phi_{D,\mu}}, 
\m{\Phi_{D',\mu'}}, be automorphisms, and
\[\Phi_{D'',\mu''} \ = \ 
\Phi_{D',\mu'}\circ\Phi_{D,\mu} \ , \quad
\ \Phi_{\widehat{D},\widehat{\mu}} \ = \ 
(\Phi_{D,\mu})^{-1} \ . \]
Using the formulas of \ref{str_aut} or directly, it is easy to prove that
\[D'' \ = \ D'+\mu'D \ , \quad \mu'' \ = \ \mu\mu' \ , \]
\[\widehat{D} \ = \ -\frac{1}{\mu_0}D \ , \ \quad \widehat{\mu} \ = \ 
\frac{1}{\mu} \ . \]
\end{subsub}

\end{sub}

\sepsub

\Ssect{The case $n=3$}{nequ3}

The following computations will be used in \ref{ext_DB} to describe the 
obstruction to extend primitive double schemes to multiplicity 3.

\sepprop

\begin{subsub}\label{nequ3_1}\rm
We can write the automorphism $\Phi$ of 
\m{\ko_X(U)[t]/(t^3)} as follows: for every \m{\alpha\in\ko_X(U)}
\begin{equation}\label{equ_0}\Phi(\alpha) \ = \ \alpha+D(\alpha)t+E(\alpha)t^2 
\ , \end{equation}
where $D$ is a derivation of \m{\ko_X(U)} and $E$ is a {\em $D$-operator}, i.e. 
a $\C$-linear map \ \m{\ko_X(U)\to\ko_X(U)} \ such that \[E(\alpha\beta) \ = \ 
\alpha E(\beta)+\beta E(\alpha)+D(\alpha)D(\beta)\]
for every \ \m{\alpha,\beta\in\ko_X(U)}. We have also
\[\Phi(t) \ = \ \mu t \ , \]
where \ \m{\mu=\mu_0+\mu_1 t\in\ko_X(U)[t]/(t^2)}, with 
\m{\mu_0,\mu_1\in\ko_X(U)}, \m{\mu_0} invertible.
\end{subsub}

\sepsubsub

\begin{subsub}\label{nequ3_2} Properties of $D$-operators - \rm We have
\begin{enumerate}
\item[(i)] If $E$, $E'$ are $D$-operators, then $E-E'$ is a derivation.
\item[(ii)] If $E$ is a $D$-operator and $D_0$ a derivation, then $E+D_0$ is a 
$D$-operator.
\item[(iii)] $\dsp\frac{1}{2}D^2=\frac{1}{2}D\circ D$ \ is a $D$-operator.
\end{enumerate}
It follows that the $D$-operators are the maps of the form
\[E \ = \ \frac{1}{2}D^2 + D^{(1)} \ , \]
where \m{D^{(1)}} is a derivation.
\end{subsub}

\sepsubsub

\begin{subsub} Notations -- \rm We will denote the automorphism $\Phi$ defined 
as in (\ref{equ_0}), with \Nligne \m{\dsp E=\frac{1}{2}D^2+D^{(1)}}, as
\[\Phi \ = \ \Phi_{D,\mu,D^{(1)}} \ . \]
The automorphism of \m{\ko_X(U)[t]/(t^2)} induced by \m{\Phi_{D,\mu,D^{(1)}}} 
is \m{\Phi_{D,\mu_0}}.
\end{subsub}

\sepsubsub

\begin{subsub} Composition and inverse -- \rm Let \m{\Phi_{D,\mu,D^{(1)}}}, 
\m{\Phi_{D',\mu',{D'}^{(1)}}}, be automorphisms, and
\[\Phi_{D'',\mu'',{D''}^{(1)}} \ = \ 
\Phi_{D',\mu',{D'}^{(1)}}\circ\Phi_{D,\mu,D^{(1)}} \ , \quad
\ \Phi_{\widehat{D},\widehat{\mu},\widehat{D}^{(1)}} \ = \ 
(\Phi_{D,\mu,D^{(1)}})^{-1} \ . \]
Using the formulas of \ref{str_aut} or directly, it is easy to prove the
\end{subsub}

\sepprop

\begin{subsub}\label{lem5}{\bf Lemma: } We have
\[D'' \ = \ D'+\mu'_0D \ , \]
\[\mu''_0 \ = \ \mu_0\mu'_0 \ , \qquad \mu''_1 \ = \ 
\mu'_0D'(\mu_0)+{\mu'_0}^2\mu_1+\mu_0\mu'_1 \ , \]
\[{D''}^{(1)} \ = \ {D'}^{(1)}+{\mu'_0}^2D^{(1)}+\mu'_1D
-\frac{1}{2}(D'+\mu'_0D)(\mu'_0)D+\frac{1}{2}\mu'_0(D'D-DD') \ , \]
\[\widehat{D} \ = \ -\frac{1}{\mu_0}D \ , \]
\[\widehat{\mu}_0 \ = \ \frac{1}{\mu_0} \ , \quad
\widehat{\mu}_1 \ = \ \frac{D(\mu_0)-\mu_1}{\mu_0^3} \ , \]
\[\widehat{D}^{(1)} \ = \ 
-\frac{1}{\mu_0^2}D^{(1)}+\frac{1}{\mu_0^3}\big(\mu_1-
\frac{1}{2}D(\mu_0)\big)D \ . \]
\end{subsub}

\sepprop

\begin{subsub}\label{can_ext} Canonical extensions -- \rm Let $D$ be a 
derivation of \m{\ko_X(U)} and \ \m{\mu=\mu_0+\mu_1t\in\ko_X[t]/(t^2)}, with 
\m{\mu_0} invertible. Let
\[D^{(1)} \ = \ \frac{\mu_1-D(\mu_0)}{2\mu_0}D\]
and
\[\Psi(D,\mu) \ = \ \Phi_{D,\mu,D^{(1)}} \ . \]
Let \m{\widehat{D}} and \m{\widehat{\mu}} be defined as in lemma \ref{lem5}.
An easy computation using lemma \ref{lem5} shows that
\end{subsub}

\sepprop

\begin{subsub}\label{lem6}{\bf Lemma: } We have \ \m{\Psi(D,\mu)^{-1}=
\Psi(\widehat{D},\widehat{\mu})} .
\end{subsub}

\sepprop

So we have defined an extension of \m{\Phi_{D,\mu_0}} to an automorphism of 
\m{\ko_X(U)[t]/(t^3)} in such a way that the extension of the inverse is the 
inverse of the extension.

Let \m{\Psi(D,\mu)}, \m{\Psi(D',\mu')}, be automorphisms, 
and
\[\Phi_{D'',\mu'',{D''}^{(1)}} \ = \ \Psi(D',\mu')\circ\Psi(D,\mu)\]
(with \m{D''}, \m{\mu''} as in lemma \ref{lem5}).
In general we don't have \ \m{\Phi_{D'',\mu'',{D''}^{(1)}}=\Psi(D'',\mu'')}.
The following result can be computed using lemma \ref{lem5}, and will be used 
in \ref{ext_MV}.

\sepprop

\begin{subsub}\label{prop0}{\bf Proposition: } We have
\[\Psi(D',\mu')\circ\Psi(D,\mu)\circ\Psi(D'',\mu'')^{-1} \ = \ 
\Phi_{0,1,\widehat{D}^{(1)}} \ , \]
with
\[\widehat{D}^{(1)} \ = \ 
\frac{1}{2}\Big(\mu'_1D+\big(D(\mu'_0)+\frac{\mu'_0}{\mu_0}D(\mu_0)-
\frac{\mu'_0}{\mu_0}\mu_1\big)D'+\mu'_0(D'D-DD')\Big) \ . \]
\end{subsub}

\end{sub}

\sepsub

\Ssect{The canonical sheaf of primitive multiple schemes and extension in 
higher multiplicity}{can_sch}

This section contains technical results that are used in \ref{pro_sp}.
We use the notations of \ref{PMV_def} to \ref{gn}.

\sepprop

\begin{subsub}\label{Can_sh} Canonical sheaves and isomorphisms -- \rm Let $Y$, 
$Z$ be schemes and \ \m{\phi:Y\to Z} \ an isomorphism. Then $\phi$ induces a 
canonical isomorphism
\[\theta:\phi^*(\Omega_Z)\lra\Omega_Y \ . \]
Let \m{U\subset Z} be an open subset and \ \m{f,g\in\ko_Z(U)}. Then \ 
\m{f.dg\in\Omega_Z(U)=\phi^*(\Omega_Z)(\phi^{-1}(U))}, and
\[\theta(f.dg) \ = \ f\circ\phi.d(g\circ\phi) \ = \ \phi^*(f).d(\phi^*(g)) \ .
\]
We will also note \ \m{\theta=\phi^*}.
\end{subsub}

\sepprop

\begin{subsub}\label{can_sh2} Description of canonical sheaves using open 
covers -- \rm Let \m{X_n} be a primitive multiple curve of multiplicity 
\m{n\geq 2}, such that \m{(X_n)_{red}=X}, with associated line bundle $L$ on 
$X$. Let \m{(U_i)} be an open cover of $X$ such that we have trivializations
\[\delta_i^{(n)}:U_i^{(n)}\lra U_i\times{\bf Z}_n\]
as in \ref{PMV_cons}. 

Let \m{(\tau_{ij})} be a cocycle representing \ \m{\tau\in H^1(\Omega_n)} 
(\m{\tau_{ij}\in H^0(U_{ij}^{(n)},\Omega_{X_n})}). Let \ 
\[\theta_i^{(n)}:((\delta_i^{(n)})^{-1})^*(\Omega_{U_i^{(n)}})\lra
\Omega_{U_i\times{\bf Z}_n}\]
be the canonical morphism (cf. \ref{Can_sh}). Let
\[\rho_{ij} \ = \ H^0(\theta_i^{(n)})(\tau_{ij}) \ \in \ 
H^0(U_{ij}\times{\bf Z}_n,\Omega_{U_{ij}\times{\bf Z}_n}) \ . \]
We have
\[\rho_{ij}+\theta_{ij}^{(n)}(\rho_{jk}) \ = \ \rho_{ik} \ , \]
where \[\theta_{ij}^{(n)}:((\delta_{ij}^{(n)})^{-1})^*(\Omega_{U_{ij}\times{\bf 
Z}_n})\lra\Omega_{U_{ij}\times{\bf Z}_n}\]
is the morphism described in \ref{Can_sh}, corresponding to 
\m{\delta_{ij}^{(n)}}. Conversely, every family \m{(\rho_{ij})} satisfying 
the preceding relation defines an element of \m{H^1(\Omega_{X_n})}.
\end{subsub}

\sepprop

\begin{subsub}\label{ext_can} Extension to multiplicity \m{n+1} -- \rm
Suppose that \m{X_n} can be extended to \m{X_{n+1}}, of 
multiplicity \m{n+1}, and let
\[\delta_i^{(n+1)}:U_i^{(n+1)}\lra U_i\times{\bf Z}_{n+1}\]
be trivializations extending \m{\delta_i^{(n)}}.
\end{subsub}

We can view every element of \m{H^0(\ko_{U_{ij}\times{\bf Z}_n})} as 
an element of \m{H^0(\ko_{U_{ij}\times{\bf Z}_{n+1}})} (with coefficient 0 in 
degree $n$ with respect to $t$).

For every \m{\beta\in H^0(\ko_{U_{ij}\times{\bf Z}_{n+1}})}, 
\m{\beta=\sigg_{p=0}^n\beta_kt^k}, \m{\beta_k\in\ko_X(U_{ij})},
let \m{[\beta]_n=\sigg_{p=0}^{n-1}\beta_kt^k\in H^0(\ko_{U_{ij}\times{\bf 
Z}_n})}.

Similarly, we can view every element of \m{H^0(\Omega_{U_{ij}\times{\bf Z}_n})} 
as an element of \m{H^0(\Omega_{U_{ij}\times{\bf Z}_{n+1}})}. Let \ 
\m{\omega\in H^0(\Omega_{U_{ij}\times{\bf Z}_{n+1}})}. We can write
\begin{equation}\label{equ10}\omega \ = \ 
\sigg_{k=0}^nt^kb_k+\big(\sigg_{p=0}^{n-1}c_pt^p\big)dt \ , \end{equation}
with \m{b_k\in H^0(U_{ij},\Omega_X)}, \m{c_p\in\ko_X(U_{ij})}. Let
\[[\omega]_n \ = \ \sigg_{k=0}^{n-1}t^kb_k+\big(\sigg_{p=0}^{n-2}c_pt^p\big)dt 
\ , \]
Then we have \ \m{\omega\in H^0(\Omega_{U_{ij}\times{\bf Z}_n})} \ if and only 
if \m{\omega=[\omega]_n}, i.e. if and only if \m{b_n=0} and \m{c_{n-1}=0}.

Let \ 
\[\theta_i^{(n+1)}:((\delta_i^{(n+1)})^{-1})^*(\Omega_{U_i^{(n+1)}})\lra
\Omega_{U_i\times{\bf Z}_{n+1}} \ , \]
\[\theta_{ij}^{(n+1)}:((\delta_{ij}^{(n+1)})^{-1})^*(\Omega_{U_{ij}\times{\bf 
Z}_{n+1}})\lra\Omega_{U_{ij}\times{\bf Z}_{n+1}}\]
be the canonical morphisms. We have
\[\theta_{ij}^{(n+1)}(\omega) \ = \ 
\theta_{ij}^{(n)}([\omega]_n)+t^n\Upsilon_0(\omega)+\Upsilon_1(\omega)t^{n-1}dt 
\ , \]
with \m{\Upsilon_0^{ij}(\omega)\in H^0(U_{ij},\Omega_X)} and 
\m{\Upsilon_1^{ij}(\omega)\in
\ko_X(U_{ij})}. If \m{R(\theta_{ij}^{(n+1)}(\omega))} is the restriction of 
\m{\theta_{ij}^{(n+1)}(\omega)} to \m{U_{ij}\times{\bf Z}_n}, we have
\[R(\theta_{ij}^{(n+1)}(\omega)) \ = \ 
\theta_{ij}^{(n)}([\omega]_n)+\Upsilon_1^{ij}(\omega)t^{n-1}dt \ . \]
Suppose that there exist \m{x_1,\ldots,x_d\in\ko_X(U_{ij})}, \m{d=\dim(X)}, 
such that \Nligne \m{\Omega_X(U_{ij})=\C dx_1\oplus\cdots\oplus\C dx_d}.

\sepprop


\begin{subsub}\label{triv_xn} The case of a trivial \m{X_n} -- \rm We suppose 
that \m{X_n} is the trivial primitive multiple scheme associated to $L$ and 
that \m{n\geq 2}. We have then \ \m{(\delta_{ij}^{(n)})^*(\alpha)=\alpha}  \ 
for every \m{\alpha\in\ko_X(U_{ij})}, and \ 
\m{(\delta_{ij}^{(n)})^*(t)=\alpha_{ij}t}, where \ 
\m{\alpha_{ij}\in\ko_X(U_i)^*} (and \m{(\alpha_{ij})} is a cocycle which 
defines $L$). We have then, for every \m{\alpha\in\ko_X(U_{ij})},
\[(\delta_{ij}^{(n+1)})^*(\alpha) \ = \ \alpha+\eta_{ij}(\alpha)t^{n} \ , \]
where \m{\eta_{ij}} is a derivation of \m{\ko_X(U_{ij})} and \ 
\m{(\delta_{ij}^{(n+1)})^*(t)=\alpha_{ij}(1+\epsilon_{ij}t^{n-1})t}, with 
\m{\epsilon_{ij}\in\ko_X(U_{ij})}. Now let \m{\omega\in 
H^0(\Omega_{U_{ij}\times{\bf Z}_n})}, written as in \ref{equ10}, with \m{b_n=0} 
and \m{c_{n-1}=0}. We have then
\[\theta_{ij}^{(n)}(\omega) \ = \ \sigg_{k=0}^{n-1}\alpha_{ij}^kt^kb_k+
\sigg_{p=0}^{n-2}c_p\alpha_{ij}^pt^p(td\alpha_{ij}+\alpha_{ij}dt) \ , \]
\begin{eqnarray*}\theta_{ij}^{(n+1)}(\omega) & = & \theta_{ij}^{(n+1)}(b_0)+
\alpha_{ij}(1+\epsilon_{ij}t^{n-1})t\theta_{ij}^{(n+1)}(b_1)+
\sigg_{k=2}^{n-1}\alpha_{ij}^kt^k\theta_{ij}^{(n+1)}(b_k)\\
& & +\big[c_0+\eta_{ij}(c_0)t^n+c_1\alpha_{ij}(1+\epsilon_{ij}t^{n-1})t+
\sigg_{p=2}^{n-2}c_p\alpha_{ij}^pt^p\big]d\big((\delta_{ij}^{(n+1)})^*(t),
\end{eqnarray*}
hence
\[R(\theta_{ij}^{(n+1)}(\omega)) \ = \
\sigg_{k=0}^{n-1}\alpha_{ij}^kt^k\theta_{ij}^{(n+1)}(b_k)+
\big(\sigg_{p=0}^{n-2}c_p\alpha_{ij}^pt^p\big)\big(td\alpha_{ij}+(\alpha_{ij}+
n\alpha_{ij}\epsilon_{ij}t^{n-1})dt\big) \ .
\]
Suppose that \ \m{b_k=\sigg_{q=1}^d\gamma_q^{(k)}dx_q}. Then
\[t^k\theta_{ij}^{(n+1)}(b_k) \ = \ t^k\sigg_{q=1}^d
(\gamma_q^{(k)}+\mu_{ij}(\gamma_q^{(k)})t^n)(dx_q+n\mu_{ij}(x_q)t^{n-1}dt+
d(\mu_{ij}(x_q))t^n) \ . \]
Restricting to \m{U_{ij}\times{\bf Z}_n} we have
\[\theta_{ij}^{(n+1)}(b_0)_{|U_{ij}\times{\bf Z}_n} \ = \ 
b_0+n\big(\sigg_{q=1}^d\gamma_q^{(0)}\mu_{ij}(x_q)\big)t^{n-1}dt \ , \]
and \ \m{t^k\theta_{ij}^{(n+1)}(b_k)_{|U_{ij}\times{\bf Z}_n}=t^kb_k} \ if 
\m{k\geq 1}. Finally
\[R(\theta_{ij}^{(n+1)}(\omega)) \ = \ \theta_{ij}^{(n)}(\omega)+n\Big(
\big(\sigg_{q=1}^d\gamma_q^{(0)}\mu_{ij}(x_q)\big)+
c_0\alpha_{ij}\epsilon_{ij}\Big)t^{n-1}dt \ , \]
whence
\[\Upsilon_1^{ij}(\omega) \ = \ 
n\Big(\big(\sigg_{q=1}^d\gamma_q^{(0)}\mu_{ij}(x_q)\big)+
c_0\alpha_{ij}\epsilon_{ij}\Big) \ = \ 
n(\pline{\mu_{ij},b_0}+c_0\alpha_{ij}\epsilon_{ij}) \ . \]
In particular, if \m{n=1} we have \m{\omega=b_0} and \ 
\m{\Upsilon_1^{ij}(\omega)=\pline{\mu_{ij},\omega}}.
\end{subsub}

\sepprop

Let \ \m{\delta:H^1(\Omega_{X_n})\to H^2(L^n)} \ be the map deduced from the 
exact sequence \Nligne \m{0\to L^n\to\Omega_{X_{n+1}|X_n}\to\Omega_{X_n}\to 0}.

\sepprop

\begin{subsub}\label{desc_delta} Description of $\delta$ -- \rm Let \ 
\m{\tau\in H^0(\Omega_{X_n})} \ as in \ref{can_sh2}. Let \m{\ov{\rho}_{ij}} be 
\m{\rho_{ij}} viewed as an element og \m{H^0(\Omega_{U_{ij}\times{\bf 
Z}_{n+1}})}, and
\[\tau'_{ij} \ = \ H^0((\theta_i{(n+1)})^{-1})(\ov{\rho}_{ij}) \ \in \
H^0(U_{ij}^{(n+1)},\Omega_{X_{n+1}}) \ . \]
Let \ \m{\ov{\tau}_{ij}=\tau'_{ij|U_{ij}^{(n)}}\in 
H^0(U_{ij}^{(n)},\Omega_{X_{n+1}|X_n})}, which is over \ \m{\tau_{ij}\in 
H^0(U_{ij}^{(n)},\Omega_{X_n})}. Then
\[\omega_{ijk} \ = \ \ov{\tau}_{ij}+\ov{\tau}_{jk}-\ov{\tau}_{ik} \ \in \
H^0(U_{ijk},L^n) \ , \]
and \m{\omega_{ijk})} is a cocycle which represents \m{\delta(\tau)}. It is 
also represented (as in \ref{cech1}) by \ \m{(\theta_i^{(n+1)}(\omega'_{ijk})
_{|U_{ijk}\times{\bf Z}_n})}, where \ 
\m{\omega'_{ijk}=\tau'_{ij}+\tau'_{jk}-\tau'_{ik}}. We have
\begin{eqnarray*}\theta_i^{(n+1)}(\omega'_{ijk})_{|U_{ijk}\times{\bf Z}_n} & = &
\big(\ov{\rho}_{ij}+\theta_{ij}^{n+1}(\ov{\rho}_{jk})-
\ov{\rho}_{ik}\big)_{|U_{ijk}\times{\bf Z}_n}\\
& = & \Upsilon^{ij}_1(\rho_{jk})t^{n-1}dt \ .
\end{eqnarray*}
It follows that \m{\delta(\tau)} is represented (in the sense of \ref{cech1}) by
\m{(\Upsilon^{ij}_1(\rho_{jk}))}.
\end{subsub}

\sepprop

\begin{subsub}\label{desc_delta3} Description of $\delta$ in the case of a 
trivial \m{X_n} -- \rm We suppose also that:
\begin{enumerate}
\item[--] (i) $\tau$ is the canonical class of the line bundle 
$\ki_{X,X_{n+1}}$ on $X_n$: $\tau=\nabla_0(\ki_{X,X_{n+1}})$.
\item[--] (ii) $\ki_{X,X_{n+1}}$ is the same as the ideal sheaf of $X$ in the 
trivial extension of $X_n$ in multiplicity $n+1$.
\end{enumerate}
For example, condition (ii) will be satisfied if \m{h^1(L^{n-1})=0}. With the 
notations of \ref{triv_xn}, we have 
\m{b_0=\dsp\frac{d\alpha_{ij}}{\alpha_{ij}}} and \m{c_0=0}. Hence 
\m{\delta(\tau)} is represented by the cocycle \ 
\m{\Big(\dsp\frac{\mu_{ij}(\alpha_{jk})}{\alpha_{jk}}\Big)}.

\end{subsub}

\end{sub}

\sepsec

\section{Automorphisms of primitive multiple schemes}\label{aut_prim}

\Ssect{Derivations and automorphisms}{Deriv_1}

Let $U$ be a smooth affine variety, \m{U=\spec(A)}, of dimension \m{d>0}, such 
that the vector bundle \m{\Omega_U} is trivial, generated by 
\m{dx_1,\ldots,dx_d}, with \m{x_1,\ldots,x_d\in A}. Then for every closed point 
\m{P\in U}, the images of \m{x_i-x_i(P)}, \m{1\leq i\leq d}, form a basis of 
\m{m_P/m_P^2}.

Let \m{n\geq 2} be an integer, and \ \m{R_n=A[t]/(t^n)}.

\sepprop

\begin{subsub}\label{deriv_Rn} Derivations of \m{R_n} -- \rm They are of 
the form
\begin{equation}\label{equ5} D \ = \ \sigg_{i=1}^da_i\frac{\partial}{\partial 
x_i}+b 
t\frac{\partial}{\partial t} \ , \end{equation}
with \ \m{a_1,\ldots,a_d,b\in R_n}.
If \ \m{\alpha=\sigg_{p=0}^{n-1}\alpha_pt^p}, with \ 
\m{\alpha_0,\ldots,\alpha_{n-1}\in A}, then
\[D(\alpha) \ = \ \sigg_{i=1}^da_i\big(\sigg_{p=0}^{n-1}\frac{\partial\alpha_p}
{\partial x_i}t^p\big)+b\sigg_{p=1}^{n-1}p\alpha_pt^p \ . \]
Let \m{\text{Der}_0(R_n)} be the vector space of derivations $D$ of \m{R_n} 
such that
\[D(R_n) \ \subset \ (t) \qquad \text{and} \qquad D((t)) \ \subset \ (t^2) \ ,
\]
i.e. derivations $(\ref{equ5})$ with \m{a_1,\ldots,a_d,b} multiple of $t$. In 
other words, a derivation $D$ belongs to \m{\text{Der}_0(R_n)} if and only if 
there exists another derivation \m{D_0} of \m{R_n} such that \m{D=tD_0}.

For every \ \m{D\in\text{Der}_0(R_n)}, and every integer \m{k\geq 1}, we have \ 
\m{D^k(R_n)\subset(t^k)} (where \m{D^k} denotes the $k$ times composition of 
$D$).
\end{subsub}

\sepprop

\begin{subsub}\label{auto_Rn} Automorphisms from derivations -- \rm
Let \ \m{D\in\text{Der}_0(R_n)}. Then we have \m{D^k=0} if \m{k\geq n}, so we 
can define the map
\begin{equation}\label{equ6}
\chi_D \ = \ \sigg_{k\geq 0}\frac{1}{k!}D^k:R_n\lra R_n \ ,
\end{equation}
such that for every \m{\alpha\in R_n}, the terms of degree 0 of $\alpha$ and 
\m{\chi_D(\alpha)}, with respect to $t$, are the same.
\end{subsub}

Let \m{D_0} be a derivation of \m{R_n} and \m{D\in\text{Der}_0(R_n)}. Then \
\m{D+t^{n-1}D_0\in\text{Der}_0(R_n)}, and we have \ \m{\chi_{D+t^{n-1}D_0}=
\chi_D+t^{n-1}D_0}.

\sepprop

\begin{subsub}\label{prop11}{\bf Proposition: } {\bf 1 -- } The map \m{\chi_D} 
is an automorphism of $\C$-algebras.

{\bf 2 -- } If \m{D,D'\in\text{Der}_0(R_n)} commute, then \ 
\m{\chi_{D+D'}=\chi_D\circ\chi_{D'}}.
\end{subsub}
\begin{proof} The fact that \m{\chi_D} is a morphism of $\C$-algebras and {\bf 
2-} are immediate. From {\bf 2-} we see that \m{\chi_{-D}} is the inverse of 
\m{\chi_D}, so \m{\chi_D} is an automorphism.
\end{proof}

\sepprop

\begin{subsub}{\bf Examples:} \rm
If \m{n=2} and \m{D,D'\in\text{Der}_0(R_2)}, then \m{DD'=D'D=0} (multiple of 
\m{t^2}).

If \m{n=3} and \m{D,D'\in\text{Der}_0(R_3)}, then we have 
\[\chi_D\circ\chi_{D'} \ = \ \chi_{D+D'}+\frac{1}{2}(D'\circ D-D\circ D') \
= \ \chi_{D+D'+\frac{1}{2}(D'\circ D-D\circ D')} \]
(\m{D'\circ D-D\circ D'} is a multiple of \m{t^2}).
\end{subsub}

\sepprop

For every \m{D\in\text{Der}_0(R_n)} have \ \m{\chi_D\in\kg_n(U)} (cf. 
\ref{PMV_cons}). Let \m{\kg_n^0(U)} denote the subgroup of \m{\kg_n(U)} of 
automorphisms $\gamma$ such that \ \m{\gamma(t)=(1+t\epsilon)t} \ for some 
\m{\epsilon\in R_n}. It is also the set of \m{\gamma\in\kg_n(U)} which can be 
written as \ \m{\gamma=I+t\psi}, for some map \m{\psi:R_n\to R_{n-1}}.

We have an obvious ``division by $t$'' map \ \m{(t)\to R_{n-1}}, that we can 
denote by multiplication by \m{\dsp\frac{1}{t}}. With the notations of 
\ref{str_aut} we have \ \m{\chi_D=\phi_{\eta,\mu}}, with
\[\eta \ = \ \frac{1}{t}\sigg_{k\geq 1}\frac{1}{k!}D^k \ , \qquad
\mu \ = \ 1+\frac{1}{t}\sigg_{k\geq 1}\frac{1}{k!}D^k(t) \ . \]

\sepprop

\begin{subsub}\label{prop12}{\bf Theorem: } The map 
\[\xymatrix@R=5pt{\text{Der}_0(R_n)\ar[r] & \kg_n^0(U)\\ D\fmaps[r] & \chi_D}\]
is a bijection.
\end{subsub}
\begin{proof} It similar to the proof of theorem 4.2.5 of \cite{dr1}. One can 
also define the inverse using the logarithm
\xmat{I+t\psi\fmaps[r] & 
\dsp\sigg_{k=1}^{n-1}(-1)^{k-1}\frac{1}{k}(t\psi)^k 
.}
\end{proof}

\sepprop

\begin{subsub}{\bf Remark: } \rm If we allow in $(\ref{equ6})$ derivations $D$ 
which are not multiple of $t$ we could still get a convergent series but not in 
\m{R_n} (in the ring of holomorphic functions on $U$ instead), or the terms of 
degree 0 of \m{\in R_n} and \m{\chi_D(\alpha)} can be distinct. For example in 
this way, with \m{A=\C[x]}, if we take \
\m{D=\dsp\frac{\partial}{\partial x}}, then for every \m{\alpha\in\C[x]} we 
have \ \m{\chi_D(\alpha)=\alpha(x+1)}.
\end{subsub}

\end{sub}

\sepsub

\Ssect{Sheaves of differentials}{sh_diff}

Let $X$ be a smooth connected variety of dimension \m{d>0} and $n$ an integer 
such that \m{n\geq 2}. Let \m{X_n} be a primitive multiple scheme, of 
multiplicity $n$, such that \ \m{(X_n)_{red}=X}. We give here some properties 
of the sheaf of differentials \m{\Omega_{X_n}} and of its dual \ 
\m{\kt_n=(\Omega_{X_n})^*}, that will be used in 8 and 9 . Recall that we have 
a canonical filtration by primitive multiple schemes
\[X_1=X\subset X_2\subset\cdots\subset X_{n-1}\subset X_n\]
(cf. \ref{PMV_def}).

The sheaves \m{\Omega_{X_n}} and \m{\kt_n} are {\em quasi locally free} (cf. 
\cite{dr2}, \cite{dr4}, i.e. they are locally isomorphic to a direct sum of 
sheaves \m{\ko_{X_i}}, \m{1\leq i\leq n}). They are locally isomorphic to \ 
\m{d.\ko_{X_n}\oplus\ko_{X_{n-1}}}: if \m{P\in X} and \m{x_1,\ldots,x_d\in 
m_{{X_n}_P}} generate \m{m_{X,P}/m_{X,P}^2}, and \m{z\in\ko_{{X_n}_P}} is a 
generator of the ideal of $X$, then
\[\Omega_{{X_n}_P} \ \simeq \ \ko_{{X_n}_P}dx_1\oplus\cdots\oplus
\ko_{{X_n}_P}dx_d\oplus\ko_{{X_n}_P}dz \ , \]
and we have \ \m{z^{n-1}dz=\frac{1}{n}d(z^n)=0} 
(\m{\ko_{{X_n}_P}dz\simeq\ko_{X_{n-1,P}}}).

Let $\ke$ be a coherent sheaf on \m{X_n}. For \m{0\leq i\leq n}, let \ 
\m{\ke_i=\ki_X^i\ke\subset\ke}, and \m{\ke^{(i)}\subset\ke} the subsheaf of 
elements annihilated by \m{\ki^i_X}. The filtration
\[0=\ke_n\subset\ke_{n-1}\subset\cdots\subset\ke_1\subset\ke=\ke_0\]
is called the {\em first canonical filtration of $\ke$}, and
\[0=\ke^{(0)}\subset\ke^{(1)}\subset\cdots\subset\ke^{(n-1)}\subset\ke^{(n)}=\ke
\]
the {\em second canonical filtration} of $\ke$.

For \m{0\leq i<n}, let
\[G_i(\ke) \ = \ \ke_i/\ke_{i+1} \ , \qquad G^{(i+1)}(\ke)=\ke^{(i+1)}/\ke^{(i)}
\ . \]
If \m{1\leq i<n} we have \ \m{(\ke^*)^{(i)}=(\ke/\ke_i)^*}, and it follows 
easily that if $\ke$ is quasi locally free, then 
\[G^{(i+1)}(\ke^*) \ = \ G_i(\ke)^*\ot L^{n-1} \ . \]
In the simplest case, i.e. if
$\ke$ is locally free, the two filtrations are the same, and we have, if 
\m{E=\ke_{|X}}, \m{G_i(\ke)\simeq E\ot L^i} \ for \m{0\leq i<n}.
\sepprop

\begin{subsub} The case of \m{\Omega_{X_n}} and \m{\kt_n} -- \rm We have 
\m{\Omega_{X_n|X}\simeq\Omega_{X_2|X}}, and an exact sequence
\[0\lra L\lra\Omega_{X_2|X}\lra\Omega_X\lra 0 \ , \]
associated to \ \m{\sigma\in Ext^1_{\ko_X}(\Omega_X,L)=H^1(T_X\ot L)}. We have 
seen in \ref{g2} that \m{\sigma=0} if and only \m{X_2} is trivial, and that if 
\m{\sigma\not=0}, then \m{\C\sigma\in\P(H^1(T_X\ot L))} defines completely 
\m{X_2}. We have
\[G_i(\Omega_{X_n}) \ \simeq \ \Omega_{X_2|X}\ot L^i \qquad \text{for} \ 1\leq 
i<n-1 \ , \]
and \ \m{G_{n-1}(\Omega_{X_n})\simeq\Omega_X\ot L^{n-1}}.
\end{subsub}
For every \ \m{\nabla\in H^1(\Omega_X)=\Ext^1_{\ko_X}(\ko_X,\Omega_X)}, let \ 
\m{0\lra\Omega_X\lra E_\nabla\lra\ko_X\lra 0} \ be the exact sequence 
associated to $\nabla$, where \m{E_\nabla} is a rank \m{d+1} vector bundle on 
$X$. Recall that \m{\nabla_0(L)\in H^1(\Omega_X)} denotes the canonical class 
of $L$ (cf. \ref{can_cl}).

\sepprop

\begin{subsub}\label{prop20}{\bf Theorem: } We have \ 
\m{\Omega_{X_n}^{(1)}\simeq E_{\nabla_0(L)}\ot L^{n-1}}.
\end{subsub}
\begin{proof} Suppose that \m{X_n} is constructed as in \ref{PMV_cons}, using 
the open cover \m{(U_i)_{i\in I}} of $X$. We suppose also that for every 
\m{i\in I}, there are \m{x_1,\ldots,x_d\in\ko_X(U_i)} such that 
\m{\Omega_{X|U_i}} is generated by \m{dx_1,\ldots,dx_n}.
We will first construct \m{ \Omega_{X_n}} by the method and notations of 
\ref{const_sh}. We take \ \m{\ke_i=\Omega_{U_i\times {\bf Z}_n}}. We have 
canonical isomorphisms induced by \m{\delta_{ij}}
\[\Psi_{ij}:\delta_{ij}^*(\Omega_{U_{ij}\times {\bf Z}_n})\lra
\Omega_{U_{ij}\times {\bf Z}_n} \ . \]
Suppose that \ 
\m{\delta_{ij}=\phi_{\eta^{ij},\mu^{ij}}}, cf. \ref{str_aut}), and let \ 
\m{\lambda_k^{ij}=\eta^{ij}(x_k)} \ for \m{1\leq k\leq d}.
We have then
\begin{eqnarray*} \Psi_{ij}(dx_k) & = & d(x_k+\lambda_k^{ij}t)\\ & = & 
dx_k+d(\lambda_k^{ij})t+\lambda_k^{ij}dt \ , \\
\Psi_{ij}(dt) & = & d(\mu^{ij}t)\\ & = & d(\mu^{ij})t+\mu^{ij}dt \ .
\end{eqnarray*}
Let
\[\Psi_{ij}^{(1)}:\delta_{ij}^*(\Omega_{U_{ij}\times {\bf Z}_n})^{(1)}=
(\Omega_{U_{ij}\times {\bf Z}_n})^{(1)}\lra
(\Omega_{U_{ij}\times {\bf Z}_n})^{(1)}\]
be the isomorphism induced by \m{\Psi_{ij}}. We have \ \m{(\Omega_{U_{ij}\times
{\bf Z}_n})^{(1)}\simeq\Omega_{X|U_{ij}}\oplus\ko_{U_{ij}}}, and the sheaf 
\m{(\Omega_{U_{ij}\times {\bf Z}_n})^{(1)}} is generated by 
\m{t^{n-1}dx_1,\ldots,t^{n-1}dx_d} and \m{t^{n-2}dt}. We have
\begin{eqnarray*}\Psi_{ij}^{(1)}(t^{n-1}dx_k) & = & 
(\mu^{ij}_0)^{n-1}t^{n-1}dx_k \ ,\\
\Psi_{ij}^{(1)}(t^{n-2}dt) & = & (\mu^{ij}_0)^{n-1}(\frac{d\mu^{ij}_0}
{\mu^{ij}_0}t^{n-1}+t^{n-2}dt) \ .
\end{eqnarray*}
In other words, \m{\Psi_{ij}^{(1)}:\Omega_{X|U_{ij}}\oplus\ko_{U_{ij}}
\to\Omega_{X|U_{ij}}\oplus\ko_{U_{ij}}} is defined by the matrix \Nligne 
\m{\dsp(\mu^{ij}_0)^{n-1}\begin{pmatrix}
I & \frac{d\mu^{ij}_0}{\mu^{ij}_0}\\ 0 & 1
\end{pmatrix}}. Proposition \ref{prop20} follows then from \ref{cech3} and 
\ref{I_X}.
\end{proof}

\sepprop

It follows easily that we have \ \m{G^{(i)}(\Omega_{X_n})\simeq 
E_{\nabla_0(L)}\ot L^{n-i}} \ for \m{1\leq i<n} and \ 
\m{G^{(n)}(\Omega_{X_n})\simeq\Omega_X}.

\sepprop

\begin{subsub}\label{coro2}{\bf Corollary: }{\bf 1 -- } We have \ 
\m{G_i(\kt_n)\simeq (E_{\nabla_0(L)})^*\ot L^i} \ for \ \m{0\leq i<n-1} and 
\Nligne \m{G_{n-1}(\kt_n)\simeq T_X\ot L^{n-1}}.

{\bf 2 -- } We have \ \m{G^{(i)}(\kt_n)\simeq(\Omega_{X_2|X})^*\ot L^{n-i}} \ 
for \ \m{1\leq i\leq n-1} and \ \m{G^{(n)}(\kt_n)\simeq T_X}.
\end{subsub}

\end{sub}

\sepsub

\Ssect{Automorphisms of primitive multiple schemes}{aut_PMS}

We use the notations of \ref{sh_diff}. Suppose that \m{X_n} corresponds to 
\m{\gamma\in H^1(X,\kg_n)}. Let \m{\AAut(X_n)} be the sheaf of groups of 
automorphisms of \m{X_n} leaving $X$ invariant, which can be identified with 
\m{\kg_n^\gamma} (cf. \ref{coh_gr}). Suppose that \m{n\geq 3}. From lemma 
\ref{gn_lem2} we have an exact sequence of sheaves of groups on $X$
\[0\lra(\Omega_{X_2|X})^*\ot L^{n-1}\lra\AAut(X_n)\lra\AAut(X_{n-1})\lra 0 \ .
\]
Let \ \m{\Aut(X_n)=H^0(\AAut(X_n))} \ be the group of global automorphisms of 
\m{X_n} leaving $X$ invariant.

From \ref{I_X} we have a canonical morphism \ \m{\xi_n:\kg_n\to\ko_X^*}, 
inducing \ \m{\xi_n^\gamma:\kg_n^\gamma\to(\ko_X^*)^\gamma=\ko_X^*} and
\[H^0(\xi_n^\gamma):\Aut(X_n)\lra\C^* \ , \]
sending to an automorphism of \m{X_n} the induced automorphism of $L$. 
Let \Nligne \m{\Aut_0(X_n)=\ker(H^0(\xi_n^\gamma))}.

\sepsubsub

\begin{subsub}\label{Auto_der} Automorphisms from derivations 2 -- \rm
We have \ \m{\ki_X\kt_n\simeq\kt_{n-1}\ot\ki_X}: this follows easily from the 
fact that \m{\ki_X} is a line bundle on \m{X_{n-1}}. Let \ \m{D\in 
H^0(X_{n-1},\ki_X\kt_n)}. Then for every affine open subset $U$ of \m{X_n} such 
that \ \m{\ki_{X|U}\simeq\ko_{X_{n-1}|U}}, \m{D_{|U}} is a derivation of 
\m{\ko_{X_n}(U)} such that \ \m{\imm(D_{|U})\subset\ki_{X|U}} \ and \ 
\m{D_{|U}(\ki_{C|U})\subset\ki_{X|U}^2}, and we obtain from \ref{Deriv_1} an 
automorphism \m{\chi_{D_{|U}}} of $U$. These automorphisms glue to define an 
automorphism \m{\chi_D} of \m{X_n}. From theorem \ref{prop12} we deduce
\end{subsub}

\sepprop

\begin{subsub}\label{prop22}{\bf Theorem: } For every \m{\phi\in\Aut_0(X_n)}, 
there is an unique \ \m{D\in H^0(X_{n-1},\ki_X\kt_n)} \ such that \ 
\m{\phi=\chi_D}.
\end{subsub}

\sepprop

Of course, if \m{n>2} the map \m{D\to\chi_D} is not a morphism of groups, so we 
have \Nligne \m{\Aut_0(X_n)\simeq H^0(X_{n-1},\ki_X\kt_n)} \ as sets only. But 
if \m{n=2}, it is actually a morphism of groups, and \ \m{\ki_X\kt_n=T_X\ot L}.

\sepprop

\begin{subsub}\label{Auto_2} Automorphisms of double primitive schemes -- \rm
Suppose first that \m{X_2} is not trivial. Then the automorphisms of \m{X_2} 
leaving $X$ invariant are of the form \ \m{\Theta_D=I_{X_2}+D}, where $D$ is a 
sheaf of derivations with value in \m{\ki_X}. With the notations of 
\ref{sh_diff}, at \m{P\in X}, we have
\[\kt_{2,P} \ \simeq \ \ko_{X_2,P}\frac{\partial}{\partial x_1}\oplus\cdots
\oplus\ko_{X_2,P}\frac{\partial}{\partial 
x_1}\oplus\ko_{X_2,P}.z\frac{\partial}{\partial z} \ . \]
It follows that we must have \ \m{D\in 
H^0(X,G^{(1)}(\kt_2))=H^0(X,(\Omega_{X_2|X})^*\ot L)}. From the \Nligne exact 
sequence of lemma \ref{g2_lem2} and the fact that \m{X_2} is not trivial, we 
see that \Nligne \m{H^0(X,(\Omega_{X_2|X})^*\ot L)\simeq H^0(X,T_X\ot L)}. It 
follows that

\sepprop

\begin{subsub}{\bf Corollary: }\label{lem7} If \m{X_2} is not trivial, then
\ \m{\Aut(X_2)=\Aut_0(X_2)=H^0(T_X\ot L)}.
\end{subsub}

\sepprop

If \m{X_2} is trivial, see \ref{triv_sh}.
\end{subsub}

\sepprop

\begin{subsub} The general case -- \rm We suppose now that \m{n\geq 2}, and 
that \m{\Aut_0(X_n)} is trivial. This is the case for example if \ 
\m{h^0(T_X\ot L^p)=0} \ for \ \m{1\leq p\leq n-1} and \ \m{h^0(L^p)=0} \ for 
\ \m{1\leq p\leq n-2} (using corollary \ref{coro2} )
\end{subsub}

\sepprop

\begin{subsub}\label{prop21}{\bf Theorem: } If \ \m{\imm(H^0(\xi_n^\gamma))} 
contains a number $\lambda$ such that \m{\lambda^p\not=1} for every integer $p$ 
such that \m{1\leq p\leq n-1}, then \m{X_n} is trivial.
\end{subsub}
\begin{proof} By induction on $n$. The result is true for \m{n=2} by corollary 
\ref{lem7}.

Suppose that \m{n\geq 3} and that it is true for \m{n-1}. Suppose that 
\m{\imm(H^0(\xi_n^\gamma))} contains a number 
$\lambda$ such that \m{\lambda^i\not=1} for every integer $i$ such that 
\m{1\leq i\leq n-1}. Then, by the induction hypothesis, \m{X_{n-1}} is trivial,
it is the \m{(n-1)}th infinitesimal neighborhood of $X$ in \m{L^*}.

We now use the notations of \ref{PMV_cons} and \ref{str_aut}. Let 
\m{(\nu_{ij})} (\m{\nu_{ij}\in\C^*}) be a cocycle which defines $L$. Then 
\m{X_{n-1}} is defined by the cocycle \m{(\phi_{0,\nu_{ij}})} of \m{\kg_{n-1}}. 
According to proposition \ref{gn_lem1}, \m{X_n} is defined by the cocycle 
\m{\delta_{ij}^*}, with \ \m{\delta_{ij}^*\in H^0(U_{ij},T_X\oplus\ko_X)\subset 
H^0(U_{ij},\kg_{n})}. We can write
\[\delta_{ij}^* \ = \ \phi_{D_{ij}t^{n-2},\nu_{ij}(1+\rho_{ij}t^{n-2})} \ , \]
i.e.
\[\delta_{ij}^* \ = \ \phi_{D_{ij}t^{n-2},1+\rho_{ij}t^{n-2}}\circ
\phi_{0,\nu_{ij}} \ , \]
where \m{((D_{ij},\rho_{ij}))} is a cocyle which represents \m{\sigma\in
H^1(X,(\Omega_{X_2|X})^*\ot L^{n-1})} (cf. \ref{par_ext}), in the way of 
\ref{cech}.. To prove that \m{X_n} is trivial, we must show that \m{\sigma=0}, 
i.e. that there exist derivations \m{F_i} of \m{\ko_X(U_i)}, 
\m{f_i\in\ko_X(U_i)}, such that
\begin{equation}\label{equ7}
D_{ij} \ = \ F_i-\nu_{ij}^{n-1}F_j \ , \qquad
\rho_{ij} \ = \ f_i-\nu_{ij}^{n-2}f_j \ .
\end{equation}
Let $\chi$ be an automorphism of \m{X_n}, such that, if \m{\lambda^p\not=1} for 
every integer $p$ such that \m{1\leq p\leq n-1}. From theorem \ref{prop22} and 
the hypothesis, \m{\chi_{|X_{n-1}}} is induced by the multiplication by 
$\lambda$ in \m{L^*}. If \ \m{\chi_i=\delta_i\circ\chi\circ\delta_{i}^{-1}\in
\Aut(U_i\times{\bf Z}_n)}, we can write
\[\chi_i^* \ = \ \phi_{t^{n-2}E_i,\lambda(1+\epsilon_it^{n-2}) }\ , \]
where \m{E_i} is a derivation of \m{\ko_X(U_i)} and \m{\epsilon_i\in\ko_X(U_i)}.
On can then see easily that the equations $(\ref{equ7})$ follow from the 
equality \ \m{\delta_{ij}^*\circ\chi_j^*=\chi_i^*\circ\delta_{ij}^*}, with
\[F_i \ = \ \frac{1}{1-\lambda_{n-1}}E_i \ , \qquad
f_i\ = \ \frac{1}{1-\lambda_{n-2}}\epsilon_i \ . \]
\end{proof}

\sepprop

\begin{subsub}\label{coro3}{\bf Corollary: } Suppose that \ 
\m{h^0(T_X\ot L^p)=0} \ for \ \m{1\leq p\leq n-1} and \ \m{h^0(L^p)=0} \ for 
\ \m{1\leq p\leq n-2}. If \m{X_n} is non trivial, then \m{\Aut(X_n)} is finite.
\end{subsub}

\sepprop

\begin{subsub}\label{triv_sh} Automorphisms of trivial primitive multiple 
schemes -- \rm Suppose that \m{X_n} is trivial. then \m{\Aut_0(X_n)} is a 
normal subgroup of \m{\Aut(X_n)}. In this case \m{X_n} is the 
$n$th-infinitesimal neighborhood of \m{X\subset L^*} (via the zero section), 
and for every \m{\lambda\in\C^*} the multiplication \ \m{\times\lambda:L^*\to 
L^*} \ induces a corresponding automorphism of \m{X_n}. We can see in this way 
\m{\C^*} as a subgroup of \m{\Aut(X_n)}. We have \ 
\m{\Aut(X_n)=\Aut_0(X_n).\C^*} \ and \ \m{\Aut_0(X_n)\cap\C^*=\{I_{X_n}\}}. 
Hence {\em \m{\Aut(X_n)} is the semi-direct product of \m{\C^*} and 
\m{\Aut_0(X_n)}}. For example, if \m{n=2}, as a set, 
\m{\Aut(X_2)=\Aut_0(X_2)\times\C^*}, with the group 
law
\[\big((\sigma,\lambda),(\sigma',\lambda')\big)\mapsto(\sigma+\lambda\sigma',
\lambda\lambda') \ . \]
\end{subsub}

\end{sub}

\sepsub

\Ssect{The case of pairs \m{(X_n,\L)}}{pairs}

With the notations of \ref{prim_can}, let \m{g\in\kh_n}, corresponding to a 
pair \m{(X_n,\L)}, where \m{X_n} is a primitive multiple scheme such that \ 
\m{(X_n)_{red}=X} \ and $\L$ is a vector bundle on \m{X_n} such that \ 
\m{\L_{X_{n-1}}=\ki_{X,X_n}}. Then it is easy to see, using \ref{const_sh}, 
that \m{H^0(X,\kh_n^g)} is isomorphic to the group of automorphisms $\chi$ of 
\m{X_n} leaving $X$ invariant and such that \m{\chi^*(\L)\simeq\L}.

\end{sub}

\sepsub

\Ssect{Parametrizations of the extensions of trivial double schemes to 
multiplicity $3$}{hi_mu}

We use the notations of \ref{PMV_def} to \ref{gn}. Let \m{X_2} be the trivial 
primitive double scheme such that \m{(X_2)_{red}=X}, with associated line 
bundle \m{L\in\Pic(X)}. Then \m{X_2} can be extended to a primitive 
multiple scheme \m{X_3} of multiplicity $3$ (for example the trivial one). 
Recall that \ \m{\rho_3:\kg_3\to\kg_2} \ is the canonical morphism. Let 
\m{g_i\in H^1(X,\kg_i)}, \m{i=2,3}, correspond to \m{X_i}. We have \ 
\m{H^1(\rho_3)(g_3)=g_2}. From \ref{coh_gr3} and \ref{gn_lem2} we have a 
canonical surjective map
\[\lambda_{g_3}:H^1(X,(\Omega_{X_2|X})^*\ot L^2)\lra H^1(\rho_3)^{-1}(g_2)
\]
which sends $0$ to \m{g_3}. Note that we have \ \m{(\Omega_{X_2|X})^*\simeq 
L^{-1}\oplus T_X}.

There is an action of the group \m{\Aut(X_2)} on 
\m{H^1(X,(\Omega_{X_2|X})^*\ot L^2)} such that the fibers of \m{\lambda_{g_3}} 
are the orbits of this action, 
so that there is a bijection between the set of isomorphism classes of 
primitive schemes of multiplicity 3 extending \m{X_2} and the quotient \ 
\m{H^1(X,(\Omega_{X_2|X})^*\ot L^2)/\Aut(X_2)}. The action is as follows: let 
\m{(g_{ij})} be a cocycle of \m{\kg_3} (with respect to \m{(U_i)_{i\in I}}) 
representing \m{g_3}, \m{\gamma\in H^1(X,(\Omega_{X_2|X})^*\ot L^2)}, 
represented by \m{(\gamma_{ij})} (seen as a cocycle in \m{\kg_3}, cf. 
\ref{gn_lem1}), and \m{\chi\in\Aut(X_2)}. Let \ \m{\delta_i:U_i^{(3)}\to 
U_i\times{\bf Z}_3} \ be isomorphisms such that \ 
\m{g_{ij}=(\delta_j^*)^{-1}\delta_i^*}. The induced isomorphisms \ 
\m{U_i^{(2)}\to U_i\times{\bf Z}_2} \ will also be denoted by \m{\delta_i}. Let 
\m{\chi_i=(\delta_j^*)^{-1}\chi^*\delta_i^*\in\Aut(U_i\times{\bf 
Z}_2)=\kg_2(U_i)}. Suppose that \m{\chi_i} can be extended to \ 
\m{c_i\in\kg_3(U_i)}. Then we have \ \m{\chi\gamma=(\theta_{ij})}, with
\begin{equation}\label{equ8}
\theta_{ij} \ = \ c_i\gamma_{ij}g_{ij}c_j^{-1}g_{ij}^{-1} \ .
\end{equation}
We use the notations of \ref{nequ3}. Let \m{(\nu_{ij})} be a cocycle 
representing $L$, \m{\nu_{ij}\in\ko_X(U_{ij})^*}. Then we can take
\[g_{ij} \ = \ \Phi_{0,\nu_{ij},0} \ . \]
Let \ \m{\gamma=(\eta,\epsilon)}, with \m{\eta\in H^1(X,L)} and 
\m{\epsilon\in H^1(X,T_X\ot L^2)}, represented by cocycles \m{(\eta_{ij})}, 
\m{(\epsilon_{ij})}, \m{\eta_{ij}\in\ko_X(U_{ij})}, \m{\epsilon_{ij}\in 
H^0(U_{ij},T_X)}, with the cocycle relations
\[\eta_{ik} \ = \eta_{ij}+\nu_{ij}\eta_{jk} \ , \qquad \epsilon_{ik} \ = \ 
\epsilon_{ij}+\nu_{ij}^2\epsilon_{jk}\]
(cf. \ref{cech1}). Using \ref{gn_lem2}, we see that
\[\gamma_{ij} \ = \ \Phi_{0,1+\rho_{ij}t,\epsilon_{ij}} \ . \]
According to \ref{aut_PMS}, as a set we have \ \m{\Aut(X_2)=H^0(T_X\ot 
L)\times\C^*}.

\sepprop

\begin{subsub} The action of \m{H^0(T_X\ot L)} -- \rm We suppose that \ 
\m{\chi\in\Aut_0(X_2)\simeq H^0(T_X\ot L)}. In this case we can take \ 
\m{\chi_i=\Phi_{D_i,1}} (with the notations of \ref{nequ2}), where \m{(D_i)} 
represents $\chi$ (with the cocycle relation \ \m{D_i=\nu_{ij}D_j}). We can take
\[c_i \ = \ \Phi_{D_i,1,K_i} \ , \]
where \m{K_i} is a derivation of \m{\ko_X(U_i)}. A lengthy calculation of 
formula $(\ref{equ8})$ shows then that \m{\chi\gamma=(\eta',\epsilon')}, where 
\m{\eta'} (resp. \m{\epsilon'}) is represented by the cocycle \m{(\eta'_{ij})} 
(resp. \m{(\epsilon'_{ij})}), with
\begin{equation}\label{equ9}\eta'_{ij} \ = \ \eta_{ij}+D_j(\nu_{ij}) \ , \qquad
\epsilon'_{ij} \ = \ \epsilon_{ij}-\nu_{ij}\big(\eta_{ij}+
\frac{1}{2}D_j(\nu_{ij})\big)D_j+K_i-\nu_{ij}^2K_j \ . \end{equation}
The duality \m{\Omega_X\simeq T_X^*} induces a canonical bilinear map
\[\xymatrix@R=5pt{
H^0(T_X\ot L)\times H^1(\Omega_X)\ar[r] & H^1(L)\\
(\chi,\tau)\fmaps[r] & \pline{\chi,\tau} \ ,
}\]
and we have a canonical product \ \m{H^1(L)\times H^0(T_X\ot L)\to 
H^1(T_X\ot L^2)}. Recall that \m{\nabla_0(L)} denotes the canonical class of 
$L$ in \m{H^1(\Omega_X)} (cf. \ref{can_cl}).
\end{subsub}

\sepprop

\begin{subsub}\label{prop23}{\bf Proposition: } We have
\[\eta' \ = \ \eta+\pline{\chi,\nabla_0(L)} \ , \qquad
\epsilon' \ = \ \epsilon-\eta\chi-\frac{1}{2}\pline{\chi,\nabla_0(L)}\chi \ .
\]
\end{subsub}

The proof follows easily from formulas $(\ref{equ9})$, with the conventions of 
\ref{cech}. See also the similar more detailed proof of theorem \ref{prop4}.

\sepprop

{\bf Remarks: }{\bf 1 -- } The term ``\m{K_i-\nu_{ij}^2K_j}'' in 
\m{\epsilon'_{ij}} can be suppressed (as expected), since \Nligne
\m{(K_i-\nu_{ij}^2K_j)} is a boundary.

{\bf 2 -- } Let \m{\chi'\in\Aut_0(X_2)}. The formulas of proposition 
\ref{prop23} do not show directly that \ 
\m{\chi'(\chi(\eta,\epsilon))=(\chi'\chi)(\eta,\epsilon)}. We find that
\[\chi'(\chi(\eta,\epsilon))-(\chi'\chi)(\eta,\epsilon) \ = \
\big(0,\frac{1}{2}(\pline{\chi',\nabla_0(L)}\chi-
\pline{\chi,\nabla_0(L)}\chi')\big)
\ = \ (0,\beta) \ . \]
If \ \m{\chi_i=(\delta_j^*)^{-1}\chi^*\delta_i^*=\phi_{D'_i,1}} \ , then 
$\beta$ is represented by the cocycle \Nligne 
\m{\dsp\big(\frac{\nu_{ij}}{2}(D'_j(\nu_{ij})D_j-D_j(\nu_{ij})D'_j))\big)}.
But we have
\[\frac{\nu_{ij}}{2}(D'_j(\nu_{ij})D_j-D_j(\nu_{ij})D'_j) \ = \ 
K_i-\nu_{ij}^2K_j \ , \]
with \ \m{\dsp K_i=\frac{1}{2}(D'_iD_i-D_iD'_i)}, so actually \m{\beta=0} in 
\m{H^1(T_X\ot L^2)}.

\sepprop

\begin{subsub}The action of \m{\C^*} -- \rm Let \m{\lambda\in\C^*}. For the 
corresponding element $\chi$ of \m{\Aut(X_2)} we have \ 
\m{\chi_i=\Phi_{0,\lambda}} \ , and we can take \ \m{c_i=\Phi_{0,\lambda,0}}. 
Let \ \m{(\eta',\epsilon')=\chi\gamma}. Then it follows easily that
\end{subsub}

\sepprop

\begin{subsub}\label{prop25}{\bf Proposition: } We have \ \m{\eta'=\lambda\eta} 
\ and \ \m{\epsilon'=\lambda^2\epsilon} \ .
\end{subsub}

\end{sub}

\sepsub

\Ssect{Parametrizations of the extensions of trivial schemes of multiplicity 
$n>2$}{hi_mu2}

Let \m{X_n} be the trivial primitive multiple scheme of multiplicity $n$ such 
that \m{(X_n)_{red}=X}, with associated line bundle \m{L\in\Pic(X)}. Then 
\m{X_n} can be extended to a primitive multiple scheme \m{X_{n+1}} of 
multiplicity \m{n+1} (for example the trivial one). Let \ 
\m{\rho_{n+1}:\kg_{n+1}\to\kg_n} \ be the canonical morphism. Let \m{g_i\in 
H^1(X,\kg_i)}, \m{i=n,n+1}, correspond to \m{X_i}. We have \ 
\m{H^1(\rho_{n+1})(g_{n+1})=g_n}. From \ref{coh_gr3} and \ref{gn_lem2} we have 
a canonical surjective map
\[\lambda_{g_{n+1}}:H^1(X,(\Omega_{X_2|X})^*\ot L^n)\lra 
H^1(\rho_{n+1})^{-1}(g_n)\]
which sends $0$ to \m{g_{n+1}}. We have \ \m{(\Omega_{X_2|X})^*\simeq 
T_X\oplus L^{-1}}. There is an action of the group \m{\Aut(X_n)} on 
\m{H^1(X,(\Omega_{X_2|X})^*\ot L^n)} such that the fibers of 
\m{\lambda_{g_{n+1}}} are the orbits of this action, so that there is a 
bijection between the set of isomorphism classes of primitive schemes of 
multiplicity \m{n+1} extending \m{X_n} and the quotient \ 
\m{H^1(X,(\Omega_{X_2|X})^*\ot L^n)/\Aut(X_n)}.

We will consider only the action of the subgroup \ \m{\C^*\subset\Aut(X_n)} \ on
\m{H^1(X,(\Omega_{X_2|X})^*\ot L^n)} (cf. \ref{triv_sh}). We have
\[H^1(X,(\Omega_{X_2|X})^*\ot L^n) \ \simeq \ H^1(X,L^{n-1})\times
H^1(X,T_X\ot L^n) \ . \]
Then the following result is an easy generalization of proposition \ref{prop24}:

\sepprop

\begin{subsub}\label{prop24}{\bf Proposition: } For every \ 
\m{(\eta,\epsilon)\in H^1(X,L^{n-1})\times H^1(X,T_X\ot L^n)} \ and 
\m{\lambda\in\C^*}, we have \ 
\m{\lambda.(\eta,\epsilon)=(\lambda^{n-1}\eta,\lambda^n\epsilon)}.
\end{subsub}

\end{sub}

\sepsec

\section{Extension to higher multiplicity}\label{ext_hi}

Let \m{n\geq 1} be an integer, $X$ a smooth projective variety, and \m{X_n} a 
primitive multiple scheme of multiplicity $n$ such that \ \m{(X_n)_{red}=X}. 
Let $L$ be the associated line bundle on $X$.
We will see when it is possible to extend \m{X_n} in multiplicity \m{n+1}, i.e. 
to construct a primitive multiple scheme \m{X_{n+1}} of multiplicity \m{n+1} 
whose canonical subscheme of multiplicity $n$ is isomorphic to \m{X_n}. If 
such a \m{X_{n+1}} exists, we will also see when it is possible to extend a 
vector bundle on \m{X_n} to a vector bundle on \m{X_{n+1}}. The two problems 
are related: \m{\ki_{X,X_n}} is a line bundle on \m{X_{n-1}}, and a necessary 
condition for the existence of \m{X_{n+1}} is the extension of \m{\ki_{X,X_n}} 
to a line bundle on \m{X_n} (that will be \m{\ki_{X,X_{n+1}}}).

\sepsub

\Ssect{Extension of vector bundles}{ext_VB}

Let $r$ be a positive integer. If $Y$ is a primitive multiple scheme, let 
\m{\GL(r,\ko_Y)} denote the sheaf of groups on $Y$ of invertible \m{r\times 
r}-matrices with coefficients in \m{\ko_Y} (the group law being the 
multiplication of matrices), and if $\ke$ is a coherent sheaf on $Y$, let 
\m{M(r,\ke)} denote the sheaf of groups on $Y$ of \m{r\times r}-matrices with 
coefficients in $\ke$ (the group law being the addition of matrices). There is 
a canonical bijection between \m{H^1(Y,\GL(r,\ko_Y))} and the set of 
isomorphism classes of rank $r$ vector bundles on $Y$.

Suppose that \m{X_n}, of multiplicity $n$, can 
be extended to \m{X_{n+1}} of multiplicity \m{n+1}. For every integer $p$ such 
that \m{1\leq p\leq n+1}, we can view the sheaf of groups \m{\GL(r,\ko_{X_p})} 
on \m{X_p} as a sheaf of groups on \m{X_{n+1}}, and we have a canonical 
surjective morphism \ \m{\GL(r,\ko_{X_{n+1}})\to\GL(r,\ko_{X_n})}. Its kernel 
is isomorphic to \m{M(r,L^n)}. So we have an exact sequence of sheaves of 
groups on \m{X_{n+1}}
\[0\lra M(r,L^n)\lra\GL(r,\ko_{X_{n+1}})\hfl{p_n}{}\GL(r,\ko_{X_n})\lra 0 \ . \]

Let \m{g\in H^1(X_n,\GL(r,\ko_{X_n}))}, $\E$ the corresponding vector bundle on 
\m{X_n} and \ \m{E=\E_{|X}}.

\sepprop

\begin{subsub}\label{prop3}{\bf Proposition:} We have \ \m{M(r,L^n)^g\simeq 
E^*\ot E\ot L^n}.
\end{subsub}

\begin{proof} The action of \m{\GL(r,\ko_{X_n})} is the action by conjugation 
after restriction to $X$. The result follows then easily from \ref{cech2c}. 
\end{proof}

\sepprop

Recall that we have a canonical map
\[\Delta:H^1(X_{n},\GL(r,\ko_{X_n}))\lra H^2(E^*\ot E\ot L^n)\]
such that \ \m{g\in\imm(H^1(p_n))} \ if and only \ \m{\Delta(g)=0} (cf. 
\ref{coh_gr4}). Hence {\em \m{\Delta(g)} is the obstruction to the extension of 
$\E$ to a vector bundle on \m{X_{n+1}}}. In particular, if $\E$ is a line 
bundle then this obstruction lies in \m{H^2(L^n)}.

We have a canonical exact sequence of sheaves on \m{X_n}
\[0\lra L^n\lra\Omega_{X_{n+1}|X_n}\lra\Omega_{X_n}\lra 0 \ , \]
corresponding to \ \m{\sigma_n\in\Ext^1_{\ko_{X_n}}(\Omega_{X_n},L^n)}, 
inducing \ \m{\ov{\sigma}_n\in\Ext^1_{\ko_{X_n}}(\E\ot\Omega_{X_n},\E\ot L^n)}. 
In \ref{can_cl} we have defined \ 
\m{\nabla_0(\E)\in\Ext^1(\E,\E\ot\Omega_{X_n})}, the canonical class of $\E$. 
We have a canonical product
\[\Ext^1_{\ko_{X_n}}(\E\ot\Omega_{X_n},\E\ot L^n)\times
\Ext^1_{\ko_{X_n}}(\E,\E\ot\Omega_{X_n})\lra
\Ext^2_{\ko_{X_n}}(\E,\E\ot L^n)=H^2(X,E^*\ot E\ot L^n) \ . \]

\sepprop

\begin{subsub}\label{prop4}{\bf Theorem:} We have \ 
\m{\Delta(g)=\ov{\sigma}_n\nabla_0(\E)} \ .
\end{subsub}

\begin{proof} For every subset $U$ of $X$, let \m{U^{(n)}} (resp. 
\m{U^{(n+1)}}) denote the corresponding open subset of \m{X_n} (resp. 
\m{X_{n+1}}).  Let \m{(U_i)_{i\in I}} be an open cover of $X$ such that $L$ is 
trivial on every \m{U_i}, \m{U_i^{(n+1)} \simeq U_i\times {\bf Z}_{n+1}}, and 
that $\E$ is represented by a cocycle \m{(\theta_{ij})}, 
\m{\theta_{ij}\in\Aut(\ko_{U_{ij}}^{(n)}\ot\C^r)}.
Let \m{t_i:\ko_{U_i}\to L_{|U_i}} be an isomorphism. For every \m{i,j\in I} 
such that \m{i\not=j}, let \m{\ov{\theta}_{ij}} be an extension of 
\m{\theta_{ij}} to an automorphism of \m{\ko_{U_{ij}}^{(n+1)}\ot\C^r}. Using 
the isomorphisms \ \m{U_i^{(n+1)} \simeq U_i\times {\bf Z}_{n+1}} \ and 
\ref{ex_pol}, we can assume that \ \m{\ov{\theta}_{ji}=\ov{\theta}_{ij}^{-1}}. 
Then by \ref{coh_gr}, \m{\Delta(g)} is represented by the family 
\m{(\ov{\theta}_{ij}\ov{\theta}_{jk}\ov{\theta}_{ki})}. We can write 
\[\ov{\theta}_{ij}\ov{\theta}_{jk}\ov{\theta}_{ki} \ = \ 
I_{U_{ijk}^{(n+1)}\ot\C^r}+t_i^n\rho_{ijk} \ , \]
and \m{\Delta(g)} is also represented by the family \m{(t_i^n\rho_{ijk})}. More 
precisely, let \ \m{\theta_i:\E_{|U_i^{(n)}}\to\ko_{U_i^{(n)}}\ot\C^r} \ be 
isomorphisms such that \ \m{\theta_{ij}=\theta_i\theta_j^{-1}}. Then 
\m{\Delta(g)} is represented by the true cocycle 
\m{(t_i^n.\theta_i^{-1}\rho_{ijk}\theta_i)}.

On the other hand, consider the morphism
\[\delta:H^1(\E^*\ot\E\ot\Omega_{X_n})\lra H^2(\E^*\ot\E\ot L^n)\]
coming from the exact sequence
\[0\lra \E^*\ot\E\ot 
L^n\lra\E^*\ot\E\ot\Omega_{X_{n+1}|X_n}\lra\E^*\ot\E\ot\Omega_{X_n}\lra 0 \ . \]
Then we have
\[\ov{\sigma}_n.\nabla_0(\E) \ = \ \delta(\nabla_0(\E)) \ . \]
The class \m{\nabla_0(\E)} is 
represented by the cocycle \m{((d\theta_{ij})\theta_{ij}^{-1})} in the way of 
\ref{cech}. In the ordinary way (local sections of 
\m{\E^*\ot\E\ot\Omega_{X_n}}), it is represented by the cocycle 
\m{(\theta_i^{-1}(d\theta_{ij})\theta_j)}. To describe \m{\delta(\nabla_0(E))} 
we first take \m{\theta_i^{-1}(d\ov{\theta}_{ij})\theta_j\in H^0(U_{ij}^{(n)},
\E^*\ot\E\ot\Omega_{X_{n+1}|X_n})}, which is over 
\m{\theta_i^{-1}(d\theta_{ij})\theta_j}. Then we have
\[\tau_{ijk} \ = \ 
\theta_i^{-1}(d\ov{\theta}_{ij})\theta_j-\theta_i^{-1}(d\ov{\theta}_{ik})\theta_k
+\theta_j^{-1}(d\ov{\theta}_{jk})\theta_k \ \in \ H^0(U_{ijk}^{(n)},
\E^*\ot\E\ot L^n) \ , \]
and \m{(\tau_{ijk})} is a cocycle representing \m{\delta(\nabla_0(E))}. We have
\[\ov{\theta}_{ij}\ov{\theta}_{jk} \ = \ 
(I_{U_{ijk}^{(n+1)}\ot\C^r}+t_i^n\rho_{ijk})\ov{\theta}_{ik} \ , \]
hence
\[\ov{\theta}_{ik} \ = \ 
\ov{\theta}_{ij}\ov{\theta}_{jk}-t_i^n\rho_{ijk}\theta_{ij}^0\theta_{jk}^0
\]
(where the exponent $0$ means restriction to $X$). Restricting to \m{X_n} we 
have
\[d\ov{\theta}_{ik} \ = \ 
\theta_{ij}.d\ov{\theta}_{jk}+d\ov{\theta}_{ij}.\theta_{jk}-
nt_i^{n-1}\theta_{ij}^0\theta_{jk}^0\rho_{ijk}dt_i \ . \]
An easy computation shows then that 
\[\tau_{ijk} \ = \ n t_i^{n-1}dt_i.(\theta_i^{-1}\rho_{ijk}\theta_i) \ , \]
i.e. it is the image in \m{H^0(U_{ijk}^{(n)},\E^*\ot\E\ot\Omega_{X_{n+1}|X_n})} 
of \ \m{t_i^n\rho_{ijk}\in H^0(U_{ijk}^{(n)},\C^{r*}\ot\C^r\ot L^n)}. This 
proves theorem \ref{prop4}.
\end{proof}

\sepprop

\begin{subsub} Other construction of the obstruction -- \rm To get an extension 
of $\E$ to a vector bundle on \m{X_{n+1}} we can try to build an extension
\[0\lra E\ot L^n\lra\ke\lra\E\lra 0\]
on \m{X_{n+1}} such that $\ke$ is locally free. We have an exact sequence
\xmat{0\ar[r] & H^1(\HHom(\E,E\ot L^n))\fleq[d]\ar[r] & 
\Ext^1_{\ko_{X_{n+1}}}(\E,E\ot L^{n})\ar[r]^-\zeta & 
H^0(\EExt^1_{\ko_{X_{n+1}}}(\E,E\ot L^n))\fleq[d]\xrightarrow{\nabla'} \\
& H^1(E\ot E^*\ot L^n) & & \End(E)}
\[\qquad\qquad\qquad\qquad  \xrightarrow{\hspace*{2.5cm}}
H^2(\HHom(\E,E\ot L^n)) \ = \ H^2(E\ot E^*\ot L^n) \ . \]
Suppose that the preceding extension corresponds to \ 
\m{\sigma\in\Ext^1_{\ko_{X_{n+1}}}(\E,E\ot L^{n``})}. Then one can prove that 
$\ke$ is locally free if and only \m{\zeta(\sigma)} is an automorphism. So it 
can be possible to build a locally free $\ke$ if and only if \ 
\m{\nabla'(I_E)=0}. One can then show that \ \m{\nabla'(I_E)=\Delta(g)}.
\end{subsub}

\end{sub}

\sepsub

\Ssect{Extension of multiple schemes}{ext_MV}

We use the notations of \ref{PMV}. Let \ \m{g_n\in H^1(X,\kg_n)} and \m{X_n} be
the corresponding multiple scheme of multiplicity $n$. Recall that we have an 
exact sequence 
\xmat{0\ar[r] & T_X\oplus\ko_X\ar[r] & \kg_{n+1}\ar[rr]^-{\rho_{n+1}} & &
\kg_n\ar[r] & 0 \ , }
that \ \m{(T_X\oplus\ko_X)^{g_n}\simeq(\Omega_{X_2|X})^*\ot L^n}, and that we 
have a canonical map
\[\Delta_n:H^1(X,\kg_n)\lra H^2((\Omega_{X_2|X})^*\ot L^n)\]
such that \ \m{g_n\in\imm(H^1(\rho_{n+1}))} \ if and only \ \m{\Delta_n(g_n)=0} 
(cf. \ref{coh_gr4}). Hence {\em \m{\Delta_n(g_n)} is the obstruction to the 
extension of \m{X_n} to a primitive multiple scheme of multiplicity \m{n+1}}.

Suppose that  \ \m{\Delta_n(g_n)=0}. By \ref{coh_gr3}, given an extension of 
\m{X_n} in multiplicity \m{n+1}, corresponding to \ \m{g_{n+1}\in 
H^1(\kg_{n+1})}, there is a canonical surjective map
\[\lambda_{g_{n+1}}:H^1((\Omega_{X_2|X})^*\ot L^n)\lra 
H^1(\rho_{n+1})^{-1}(g_n) \]
such that \ \m{\lambda_{g_{n+1}}(0)=g_{n+1}}.

By \ref{ext_VB}, we have also a canonical map
\[\Delta'_n:H^1(X_{n-1},\ko^*_{X_{n-1}})=\Pic(X_{n-1})\lra H^2(L^{n-1})\]
such that a line bundle $\kl$ on \m{X_{n-1}} can be extended to a line bundle 
on \m{X_n} if and only if \ \m{\Delta'_n(\kl)=0}.

By \ref{g2} we have an exact sequence
\xmat{0\ar[r] & T_X\ar[r]^-\iota & (\Omega_{X_2|X})^*\ar[r]^-\upsilon  & 
L^{-1}\ar[r] & 0 \ .}
Let \ \m{\Upsilon_n:H^1(X,\kg_n)\to \Pic(X_{n-1})} \ be the map defined as 
follows: if \m{g_n\in H^1(X,\kg_n)} and \m{X_n} is the corresponding multiple 
scheme of multiplicity $n$, then \m{\Upsilon_n(g_n)=\ki_{X,X_n}}.

\sepprop

The following result can be proved easily using the methods of \ref{gn}:

\sepprop

\begin{subsub}\label{prop6}{\bf Proposition: } The following diagram
\xmat{H^1(X,\kg_n)\ar[rr]^-{\Delta_n}\ar[d]^{\Upsilon_n} & & 
H^2((\Omega_{X_2|X})^*\ot L^n)\ar[d]^{H^2(\upsilon\ot I_{L^n})}\\
\Pic(X_{n-1})\ar[rr]^-{\Delta'_n} & & H^2(L^{n-1})}
is commutative.
\end{subsub}

\sepprop

This result is not surprising: if \m{g_n\in H^1(X,\kg_n)} and \m{X_n} is the 
corresponding primitive multiple scheme, then \m{\Delta_n(g_n)} is the 
obstruction to its extension in multiplicity \m{n+1}, and 
\m{\Delta'_n(\Upsilon_n(g_n))} is the obstruction to the extension of 
\m{\ki_{X,X_n}} to a line bundle on \m{X_n}. It follows that \
\m{\Delta_n(g_n)=0} \ implies that \ \m{\Delta'_n(\Upsilon_n(g_n))=0}.

Suppose that \ \m{\Delta'_n(\Upsilon_n(g_n))=0}. It follows that there exists \ 
\m{\eta\in H^2(T_X\ot L^n)} \ such that \ \m{\Delta_n(g_n)=H^2(\iota\ot 
I_{L^n})(\eta)}. This can be seen more precisely using \ref{prim_can}. From 
$(\ref{equ2})$ we have a canonical map
\[\Delta''_n:H^1(X,\kh_n)\lra H^2(T_X\ot L^n)\]
such that for any pair \m{(X_n,\L)} (where $\L$ is an extension of 
\m{\ki_{X,X_n}} to a line bundle on \m{X_n}), corresponding to \ \m{h_n\in 
H^1(X,\kh_n)}, there exists an extension of \m{X_n} to a primitive multiple 
scheme \m{X_{n+1}} of multiplicity \m{n+1} such that \ \m{\ki_{X,X_{n+1}}=\L} \ 
if and only if \ \m{\Delta''_n(h_n)=0}. We have also a commutative diagram
\xmat{H^1(X,\kh_n)\ar[r]^-{\Delta''_n}\ar[d]^{H^1(\epsilon_n)} & H^2(T_X\ot 
L^n)\ar[d]^{H^2(\iota\ot I_{L^n})}\\
H^1(X,\kg_n)\ar[r]^-{\Delta_n} & H^2((\Omega_{Y|X})^*\ot L^n)}

Since \ \m{\Delta'_n(\Upsilon_n(g_n))=0}, \m{\ki_{X,X_n}} can be extended to a 
line bundle $\L$ on \m{X_n}, and \m{(X_n,\L)} corresponds to \ \m{h_n\in 
H^1(X,\kh_n)}. It follows that we can take \ \m{\eta=\Delta''_n(h_n)}. 

Let $\kx$ be the set of extensions of \m{X_n} to a primitive multiple scheme 
\m{X_{n+1}} of multiplicity \m{n+1} such that \ \m{\ki_{X,X_{n+1}}\simeq\L}. 
Then there exists a canonical surjective map
\[\lambda:H^1(T_X\ot L^n)\lra\kx\]
such that the fibers of $\lambda$ are the orbits of the group of automorphisms 
$\chi$ of \m{X_n} leaving $X$ invariant and such that \m{\chi^*(\L)\simeq\L}.

\sepprop

\begin{subsub}\label{par_ext} Parametrization of the extensions -- \rm We use 
the notations of \ref{PMV_cons} and \ref{str_aut}. Suppose that \ 
\m{g_n=H^1(\rho_{n+1})(g_{n+1})}, and that \m{g_{n+1}} is defined by the 
cocycle \m{(\delta_{ij}^*)}, \ \m{\delta_{ij}^*=\phi_{\eta_{ij},\mu_{ij}}}.
Then the elements of \m{H^1(\rho_{n+1})^{-1}(g_n)} are represented by the 
cocycles of the form \m{(\tau_{ij}\delta_{ij}^*)}, with
\[\tau_{ij} \ = \ \phi_{D_{ij}t^{n-1},1+\epsilon_{ij}t^{n-1}} \ , \]
where \m{D_{ij}} is a derivation of \m{\ko_X(U_{ij})} and \ 
\m{\epsilon_{ij}\in\ko_X(U_{ij})}. The cocyle relation satisfied by 
\m{(\tau_{ij}\delta_{ij}^*)} (to ensure that it defines an element of 
\m{H^1(X,\kg_{n+1})}) is equivalent to a cocyle relation satisfied by 
\m{(D_{ij},\epsilon_{ij})}, so that it defines an element of 
\m{H^1(X,(\Omega_{X_2|X})^*\ot L^n)} in the way described in \ref{cech}
(cf. \ref{coh_gr3}).

\end{subsub}

\sepprop

\begin{subsub}\label{weigh} Extensions of trivial primitive multiple schemes -- 
\rm Let \m{X_n^0} be trivial primitive multiple scheme of multiplicity \m{n\geq 
2}, with associated smooth variety $X$ and associated line bundle $L$. We 
suppose that \ \m{h^0(X_n,\ki_X\kt_n)=0}.
Since \ \m{\Aut_0(X^0_n)=\{I_{X_n^0}\}}, we have \m{\Aut(X_n^0)=\C^*}. It 
follows from proposition \ref{prop24} that the non trivial extensions of 
\m{X_n^0} to a primitive multiple scheme of multiplicity \m{n+1} can be 
identified with the quotient \ \m{\big(H^1(X,L^{n-1})\times H^1(X,T_X\ot 
L^n)\big)/\C^*} with the action
\[\xymatrix@R=5pt{\C^*\times\big(H^1(X,L^{n-1})\times H^1(X,T_X\ot L^n)\big)
\ar[r] & H^1(X,L^{n-1})\times H^1(X,T_X\ot L^n)\\
(\lambda,(\eta,\epsilon))\fmaps[r] & (\lambda^{n-1}\eta,\lambda^n\epsilon)
}\]
i.e. it is a {\em weighted projective space}.
\end{subsub}

\end{sub}

\sepsub

\Ssect{Extension of double schemes}{ext_DB}

We use the notations of \ref{PMV_cons} and \ref{g2}: \m{X_2} is a double 
scheme with underlying smooth variety $X$ and associated line bundle $L$ on 
$X$, \m{(U_i)} is an open affine cover of $X$, and \m{X_2} is obtained by 
gluing the varieties \m{U_i\times {\bf Z}_2} using the automorphisms \ 
\m{\delta_{ij}:U_{ij}\times {\bf Z}_2\to U_{ij}\times {\bf Z}_2}. We can write 
(using the notations of \ref{struc})
\[\delta_{ij}^* \ = \ \Phi_{D_{ij},\alpha_{ij}} \ , \]
where \m{D_{ij}} is a derivation of \m{\ko_X(U_{ij})} and 
\m{\alpha_{ij}\in\ko_X(U_{ij})} is invertible. The cocycle \m{(\alpha_{ij})} 
defines the line bundle \m{L}, and we have \ 
\m{D_{ik}=D_{ij}+\alpha_{ij}D_{jk}}.

Suppose that $L$ can be extended to a line bundle on \m{X_2}, and let 
\m{(\alpha_{ij}+\beta_{ij}t)} \Nligne
(\m{\alpha_{ij}+\beta_{ij}t\in\ko_X(U_{ij})[t]/(t^2)}) be a family defining 
\m{\L} (the cocycle defining \m{\L} on \m{X_2} is 
\m{(\delta_i^*(\alpha_{ij}+\beta_{ij}t))}). The pair \m{(X_2,\L)} corresponds 
to \ \m{\eta\in H^1(X,\kh_2)} (cf. \ref{prim_can}). The cocycle 
\m{((\Phi_{D_{ij},\alpha_{ij}},\alpha_{ij}+\beta_{ij}t))} 
represents $\eta$. We have a canonical map
\[\Delta''_2:H^1(X,\kh_2)\lra H^2(T_X\ot L^2)\]
such that there exists an extension of \m{X_2} to a primitive multiple scheme 
\m{X_3} of multiplicity 3 such that \ \m{\ki_{X,X_{3}}=\L} \ if and only if \ 
\m{\Delta''_2(\eta)=0}.

\sepprop

\begin{subsub}\label{prop8}{\bf Theorem: } The cocycle \m{(\rho_{ijk})}, 
with
\[\rho_{ijk} \ = \ 
\frac{1}{2}\Big(\beta_{ij}D_{jk}+\big(D_{jk}(\alpha_{ij})+\frac{\alpha_{ij}}
{\alpha_{jk}}D_{jk}(\alpha_{jk})-\frac{\alpha_{ij}}{\alpha_{jk}}\beta_{jk}\big)
D_{ij}+\alpha_{ij}(D_{ij}D_{jk}-D_{jk}D_{ij})\Big)\]
represents \m{\Delta''_2(\eta)}.
\end{subsub}
\begin{proof}
The construction of \m{\Delta''_2(\eta)} follows \ref{coh_gr4}, using the exact 
sequence of $(\ref{equ2})$
\xmat{0\ar[r] & T_X\ot L^2\ar[r] & \kg_3^\eta\ar[r] & \kh_2^\eta\ar[r] & 0 \ . }
We take \ \m{\Psi_{ij}=\Psi(D_{ij},\alpha_{ij}+\beta_{ij}t)\in\kg_3(U_{ij})} \ 
as a pull-back of \m{(\Phi_{D_{ij},\alpha_{ij}},\alpha_{ij}+\beta_{ij}t)\in 
\kh_2(U_{ij})}. The family \m{(\Psi_{ij})} has the correct property, i.e. \ 
\m{\Psi_{ij}^{-1}=\Psi_{ij}}, by lemma \ref{lem6}, and 
\m{(\Psi_{ij}\Psi_{jk}\Psi_{ki})} represents \m{\Delta''_2(\eta)}. We have
\[\Psi_{ij}\Psi_{jk}\Psi_{ki} \ = \ \Phi_{0,1,\rho_{ijk}} \ , \]
from proposition \ref{prop0}, and the result follows.
\end{proof}

\sepprop

\begin{subsub} The case \ \m{L=\ko_X} -- \rm In this case we take 
\m{\alpha_{ij}=1}, and \m{(D_{ij})} is a cocycle representing the element \ 
\m{\sigma\in H^1(T_X)} \ which defines \m{X_2} (cf. \ref{param_2}). On the 
other hand, \m{(\beta_{ij})} is a cocycle representing \m{\theta_\L\in 
H^1(\ko_X)} associated to $\L$. We have a canonical bilinear map ({\em Poisson 
bracket})
\[\xymatrix@R=5pt{H^1(T_X)\times H^1(T_X)\ar[r] & H^2(T_X)\\
(u,v)\fmaps[r] & [u,v]
}\]
and it follows easily from theorem \ref{prop8} that 
\[\Delta''_2(\eta) \ = \ \theta_\L\sigma+\frac{1}{2}[\sigma,\sigma] \ . \]

\end{subsub}

\end{sub}

\sepsub

\Ssect{The case of surfaces -- Double primitives schemes with canonical 
associated sheaf}{surf}

We suppose that $X$ is a surface and \m{L=\omega_X}. The non trivial primitive 
double schemes with support $X$ and associated line bundle \m{\omega_X} are 
parametrized by \ \m{\P(H^1(T_X\ot\omega_X))=\P(H^1(\Omega_X))}.

Suppose that \m{\sigma\in H^1(\Omega_X)}, \m{\sigma\not=0}, and let \m{X_2} be 
the primitive double scheme corresponding to \m{\C\sigma}. The first necessary 
condition to extend \m{X_2} to multiplicity 3 is that \m{\omega_X} can be 
extended to a line bundle on \m{X_2}. This is equivalent to
\[\sigma.\nabla_0(\omega_X) \ = \ 0\]
in \ \m{H^2(\omega_X)\simeq\C}. If this is true, a sufficient condition to 
extend \m{X_2} to multiplicity 3 is that
\[h^2(T_X\ot\omega_X^2) \ = \ h^0(T_X) \ = \ 0 \ . \]
In other words

\sepprop

\begin{subsub}\label{prop28}{\bf Proposition: } If \ \m{H^0(T_X)=\nsp}, then
\m{X_2} can be extended to a primitive multiple scheme of multiplicity 3 if and 
only if \ \m{\sigma.\nabla_0(\omega_X)=0}.
\end{subsub}

\end{sub}

\sepsub

\Ssect{Relations with the construction of D. Bayer and D. Eisenbud}{b_e}

D.~Bayer and D.~Eisenbud gave in \cite{ba_ei} the construction and 
parametrization of primitive double schemes, also called {\em ribbons}. These 
constructions can be generalized to higher multiplicities, as we will see 
(mostly without proofs).

Given \m{X_n}, if \m{X_{n+1}} exists, then we have an exact sequence of sheaves 
on \m{X_n}
\begin{equation}\label{equ3b}
\xymatrix{0\ar[r] & L^n\ar[r] & \Omega_{X_{n+1}|X_n}\ar[r]^-p & 
\Omega_{X_n}\ar[r] & 0 \ ,}
\end{equation}
and \m{\Omega_{X_{n+1}|X_n}} is locally free. Let
\[d_n:\ko_{X_n}\lra\Omega_{X_n} \ , \quad 
d_{n+1}:\ko_{X_{n+1}}\lra\Omega_{X_{n+1}}\]
be the canonical derivations, and \ 
\m{\rho:\Omega_{X_{n+1}}\to\Omega_{X_{n+1}|X_n}} \ the restriction morphism. 
Then we have a commutative diagram
\xmat{\ko_{X_{n+1}}\ar[r]\ar[d]^{\rho d_{n+1}} & \ko_{X_n}\ar[d]^{d_n}\\
\Omega_{X_{n+1}|X_n}\ar[r]^-p & \Omega_{X_n}}
inducing, for every \m{x\in X_n}, a bijection between \m{\ko_{X_{n+1},x}} and 
the set of pairs \Nligne 
\m{(\alpha,\beta)\in\ko_{X_n,x}\times \Omega_{X_{n+1}|X_n,x}} \ such that \ 
\m{d_{n,x}(\alpha)=\rho_xd_{n+1,x}(\beta)}. It follows that {\em given \Nligne
\m{\sigma\in\Ext^1_{\ko_{X_n}}(\Omega_{X_n},L^n)}, there exists at most one 
extension \m{X_{n+1}} such that $\sigma$ is associated to $(\ref{equ3b})$}.

Conversely, given an exact sequence of sheaves on \m{X_n}
\begin{equation}\label{equ3}
\xymatrix{0\ar[r] & L^n\ar[r] & \ke\ar[r]^-p & \Omega_{X_n}\ar[r] & 0 \ ,}
\end{equation}
where $\ke$ is locally free, we can  define a sheaf of abelian groups $\kc$ to 
be the pullback
\xmat{\kc\ar[r]\ar[d] & \ko_{X_n}\ar[d]^{d_n}\\ \ke\ar[r]^-p & \Omega_{X_n}}
and make $\kc$ into a sheaf of $\C$-algebras with an appropriate product law 
(cf. \cite{ba_ei}, proof of theorem 1.2). We can then define \m{X_{n+1}} as \
\m{X_{n+1}=\spec(\kc)}. It follows that {\em \m{X_n} can be extended to a 
primitive multiple scheme of multiplicity \m{n+1} if and only if there exists 
an extension $(\ref{equ3})$ such that $\ke$ is locally free}.

To study the extensions of \m{X_n} in multiplicity \m{n+1} we have to consider 
extensions $(\ref{equ3})$. We have from \cite{go}, 7.3, a canonical exact 
sequence
\begin{equation}\label{equ4}
0\lra H^1(\HHom(\Omega_{X_n},L^n))\lra \Ext^1_{\ko_{X_n}}(\Omega_{X_n},L^n)\lra
\qquad\qquad\qquad
\end{equation}
\xmat{\qquad\qquad\qquad H^0(\EExt^1_{\ko_{X_n}}(\Omega_{X_n},L^n))\ar[r] & 
H^2(\HHom(\Omega_{X_n},L^n)) \ .
}

\sepprop

\begin{subsub}{\bf Lemma:}\label{lem3} We have \ 
\m{\EExt^1_{\ko_{X_n}}(\Omega_{X_n},L^n)\simeq\ko_C} \ .
\end{subsub}

\begin{proof} We use the exact sequence
\[L^n\lra\Omega_{X_{n+1}|X_n}\lra\Omega_{X_n}\lra 0 \ . \]
For every \m{x\in X}, there is an open neighborhood $U$ of $x$ in \m{X_n} such 
that \m{\ki_{X,X_n|U}} can be extended to a line bundle $\L$ of $U$. We obtain 
a locally free resolution of \m{\Omega_{X_n|U}}:
\[\cdots\lra\L^{n+1}\lra\L^n\lra\Omega_{X_{n+1}|U}\lra\Omega_{X_n|U}\lra 0 \ . 
\]
It follows that we have an isomorphism \ 
\m{\EExt^1_{\ko_{U}}(\Omega_{X_n|U},L_{|U}^n)\simeq\ko_U}. It is easy to see 
that all these isomorphisms glue and thus define the isomorphism of lemma 
\ref{lem3}.
\end{proof}

\sepprop

\begin{subsub}{\bf Lemma:}\label{lem4} 
Let \m{x\in X} and \ 
\m{\alpha\in\Ext^1_{\ko_{X_n,x}}(\Omega_{X_n,x},L_x^n)\simeq\ko_{X,x}}. Let
\[0\lra L_x^n\lra M\lra\Omega_{X_n,x}\lra 0\]
be the corresponding extension of \m{\ko_{X_n,x}}-modules. Then $M$ is free if 
and only if $\alpha$ is invertible.
\end{subsub}

\begin{proof} We use the free resolution
\xmat{\cdots\ar[r] & 
\L^{n+1}_x\ar[r] & \L^n_x\ar[r]^-\beta & \Omega_{X_{n+1}|U,x}\ar[r] & 
\Omega_{X_n,x}\ar[r] & 0 \ , 
}
where $\beta$ induces \ \m{\ov{\beta}:L^n_x\to\Omega_{X_{n+1}|U,x}}. Then, by 
the usual construction of extensions, $M$ is isomorphic to the cokernel of
\[\alpha\oplus\ov{\beta}:L_x^n\lra L_x^n\oplus\Omega_{X_{n+1}|U,x} \ . \]
Lemma \ref{lem4} follows easily.
\end{proof}

\sepprop

We have \ \m{\HHom(\Omega_{X_n},L^n)\simeq(\Omega_{X_2|X})^*\ot L^n}. Hence 
we have an exact sequence
\xmat{0\ar[r] & H^1((\Omega_{X_2|X})^*\ot L^n)\ar[r]^-\nu & 
\Ext^1_{\ko_{X_n}}(\Omega_{X_n},L^n)\ar[r]^-\tau & \C\ar[r]^-\delta & 
H^2((\Omega_{X_2|X})^*\ot L^n) \ . }
{\em it follows that \m{X_n} can be extended to primitive multiple scheme of 
multiplicity \m{n+1} if and only if \ \m{\delta(1)=0}}, and that the set of 
isomorphism classes of extensions of \m{X_n} in multiplicity \m{n+1} can be 
identified with \m{\tau^{-1}(1)}.

Let \m{X_{n+1}} be an extension of \m{X_n} in multiplicity \m{n+1} and \ 
\m{\sigma\in\tau^{-1}(1)} the corresponding element. We have a bijection
\[\xymatrix@R=5pt{\ov{\nu}:H^1((\Omega_{X_2|X})^*\ot L^n)\ar[r] & 
\tau^{-1}(1)\\ u\fmaps[r] & \sigma+\nu(u) \ .} \]
With the notations of \ref{ext_MV}, \m{\tau^{-1}(1)} can be identified with 
\m{H^1(\rho_{n+1})^{-1}(g_n)}, and we have a canonical map
\[\lambda_{g_{n+1}}:H^1((\Omega_{X_2|X})^*\ot L^n)\lra 
H^1(\rho_{n+1})^{-1}(g_n) \ . \]

The following result is proved in \cite{dr6}, 3.1, when $X$ is a curve, but the 
proof is similar in any dimension:

\sepprop

\begin{subsub}{\bf Theorem:}\label{prop7} We have, for every \ \m{u\in 
H^1((\Omega_{X_2|X})^*\ot L^n)}, \ \m{\lambda_{g_{n+1}}(u)=\ov{\nu}(nu)}.
\end{subsub}

\end{sub}

\sepsec

\section{The case of projective spaces}\label{pro_sp}

We use the notations of \ref{PMV_def} to \ref{gn}.

Let $V$ be a complex vector space of dimension \m{m+1}, \m{m\geq 2}. We suppose 
that \m{X=\P(V)=\P_m} (the set of lines in $V$). Let $k$ be an integer, and 
\m{L=\ko_{\P_m}(k)}. We will prove

\sepprop

\begin{subsub} {\bf Theorem: }\label{theo4} {\bf 1 -- } If \m{m>2}, all 
primitive multiple schemes $Y$ such that \m{Y_{red}=\P_m} are trivial.

{\bf 2 -- } There are only two non trivial primitive multiple schemes $Y$ such 
that \m{Y_{red}=\P_2}, one in multiplicity $2$, with \m{L=\ko_{\P_2}(-3)}, and 
the other in multiplicity $4$, with \m{L=\ko_{\P_2}(-1)}.
\end{subsub}
\begin{proof}
There exists a non trivial primitive double scheme \m{X_2} with 
underlying smooth variety $X$ and associated line bundle $L$ if and only if \ 
\m{h^1(T_X\ot L)\not=0}. This is true only if \m{m=2} and \m{k=-3}.

Suppose that \m{m>2}. According to \ref{ext_MV}, there can be non trivial 
extensions of a trivial primitive multiple scheme \m{X_n} to multiplicity 
\m{n+1} only if \ \m{h^1(\P_m,(\Omega_{X_2|X})^*\ot L^n)\not=0}. This never 
happens if \m{m>2}. So {\bf 1} is proved.

Suppose now that \m{m=2} and \m{k=-3}. In this 
case, according to \ref{g2}, the non trivial primitive double schemes are 
parametrized by \m{\P(H^1(T_{\P_2}(-3)))}. Since \m{h^1(T_{\P_2}(-3))=1}, there 
is just only one non trivial primitive double scheme, that we will denote by 
\m{{\bf X}_2}. We will prove

\sepprop

\begin{subsub}{\bf Theorem: }\label{theo2} The only line bundle on \m{{\bf 
X}_2} is the trivial bundle \m{\ko_{{\bf X}_2}}.
\end{subsub}

\sepprop

It follows that \m{{\bf X}_2} is not quasi-projective and cannot be extended to 
a primitive scheme of multiplicity 3. The scheme \m{{\bf X}_2} appears also in 
\cite{b_m_r} (it is an example of {\em non projective K3 
carpet}).

Hence any non trivial primitive multiple scheme with underlying smooth variety 
\m{X=\P_2} and associated line bundle $L$, and of multiplicity \m{n+1>2} is an 
extension of a trivial multiple scheme of multiplicity \m{n\geq 2}. This can 
only happen if \ \m{h^1(T_X\ot L^n)\not=0} (cf. \ref{par_ext}). In this case we 
have \m{L^n=\ko_{\P_2}(-3)}, hence \m{n=3} and \m{L=\ko_{\P_2}(-1)}. Starting 
with the trivial primitive multiple scheme \m{X_3} of multiplicity 3 (with 
\m{L=\ko_{\P_2}(-1)}), we get a non trivial primitive multiple scheme \m{{\bf 
X}_4} of multiplicity 4. Since \m{h^1(T_X\ot L^3)=1}, proposition \ref{prop25} 
implies that \m{{\bf X}_4} is unique. We will prove

\sepprop

\begin{subsub}{\bf Theorem: }\label{theo3} The only line bundle on \m{{\bf 
X}_4} is the trivial bundle \m{\ko_{{\bf X}_4}}.
\end{subsub}

\sepprop

It follows that \m{{\bf X}_4} is not quasi-projective and cannot be extended to 
a primitive scheme of multiplicity 5. Hence {\bf 2} is proved.\end{proof}

\sepsub

\Ssect{Description of some sheaves on the projective plane}{can_pr}


We assume that \m{m=2}. Let \m{x_0,x_1,x_2\in V^*} be coordinates in \m{\P_2}, 
and for i=0,1,2, let
\[U_i \ = \ \{\C(u_0,u_1,u_2)\in\P_2;u_i\not=0\} \ . \]
We will take indices in \m{\Z/3\Z}, so that \m{x_3=x_0,x_4=x_1,\cdots}.

We have \ \m{\ko_{\P_2}(U_i)=\C[\frac{x_{i+1}}{x_i},\frac{x_{i+2}}{x_i}]}.

\sepprop

\begin{subsub}\label{can_pr1} Line bundles -- \rm For every \m{x\in\P_2}, we 
have \m{\ko_{\P_2}(-1)_x=x}. For \m{i=0,1,2} we have a trivialization of 
\m{\ko_{\P_2}(-1)} on \m{U_i}, 
\m{\alpha_i:\ko_{\P_2}(-1)_{|U_i}\to\ko_{U_i}\ot\C}, such that for every 
\m{x=\C(u_0,u_1,u_2)} in \m{U_i}, we have \ \m{\alpha_{i,x}(u_0,u_1,u_2)=u_i}.
Hence \m{\ko_{\P_2}(-1)} is defined by the cocycle \m{(\alpha_{ij})}, where \ 
\m{\alpha_{ij}:U_{ij}\to\C^*}, and for every \m{x=\C(u_0,u_1,u_2)} 
in \m{U_{ij}}, we have \ \m{\dsp\alpha_{ij,x}(x)=\frac{u_i}{u_j}}.

Similarly, let \m{k\in\Z}. Then \m{\ko_{\P_2}(k)} is defined by the cocycle 
\m{(\alpha_{ij}^{(k)})}, where \ 
\m{\alpha_{ij}^{(k)}:U_{ij}\to\C^*}, and for every 
\m{x=\C(u_0,u_1,u_2)} in \m{U_{ij}}, we have \ 
\m{\dsp\alpha_{ij,x}^{(k)}(x)=\big(\frac{u_i}{u_j}\big)^{-k}}.

Let \m{U\subset\P_2} be a nonempty open subset. Then \m{H^0(U,\ko_{\P_2}(k))} 
can be identified with the set of rational functions \m{\dsp\frac{P}{Q}}, 
\m{P,Q\in\C[X_0,X_1,X_2]}, such that $P$, $Q$ are homogeneous, $Q$ does not 
vanish on $U$ and \ \m{\deg(P)-\deg(Q)=k}. If \m{k<0}, then for every \m{x\in 
U}, \m{x=\C.u}, \m{\ko_{\P_2}(k)_x=x^{-k}\in S^{-k}V}, and \ 
\m{\dsp\frac{P}{Q}(x)=\frac{P(u)}{Q(u)}.u^{-k}}.

For example, if \m{k=-3}, these rational functions \m{\dsp\frac{P}{Q}} such 
that $Q$ does not vanish on \m{U_{012}} are the linear combinations of 
\m{\dsp\frac{1}{X_0X_1X_2}}, \m{\dsp\frac{X_k^a}{X_i^bX_j^c}}, 
\m{\dsp\frac{X_i^dX_j^e}{X_k^f}}, \m{0\leq i<j\leq 2}, \m{0\leq k\leq 2}, 
\m{k\not=i,j} \m{a,b,c,d,e,f\geq 0}, \m{a-b-c=3}, \m{d+e-f=3}. We have \ 
\m{H^2(\ko_{\P_2}(-3))\simeq\C}, generated by the element 
corresponding to the cocycle \m{\dsp\Big(\frac{1}{X_0X_1X_2}\Big)}.
\end{subsub}

\sepprop

\begin{subsub}\label{can_pr2} Canonical bundle and tangent bundle -- \rm We 
have \ \m{H^1(\Omega_{\P_2})\simeq\C}, and it is generated by the canonical 
class of \m{\ko_{\P_2}(-1)}, which is represented by the cocycle 
\m{(\rho_{ij})}, with
\[\rho_{ij} \ = \ \frac{x_j}{x_i}d\Big(\frac{x_i}{x_j}\Big) \ . \]

The \m{\ko_{\P_2}(U_i)}-module \m{H^0(T_{\P_2}\ot\ko_{\P_2}(-3))} is generated 
by the derivations \m{\dsp\frac{\partial}{\partial(x_{i+1}/x_i)}}, 
\m{\dsp\frac{\partial}{\partial(x_{i+2}/x_i)}}.
We have \ \m{\Omega_{\P_2}\simeq T_{\P_2}\ot\ko_{\P_2}(-3)}, and \ 
\m{H^1(T_{\P_2}\ot\ko_{\P_2}(-3))\simeq\C} \ is generated by the element 
\m{\rho'} corresponding to the cocycle \m{(\rho'_{ij})} with
\[\rho'_{i,i+1} \ = \ 
\frac{x_i}{x_{i+1}}\frac{\partial}{\partial(x_{i+2}/x_i)}\ot
\frac{1}{x_i^3}\]
(where \m{\dsp\frac{1}{x_i^3}} corresponds to a section of \m{\ko_{\P_2}(-3)} 
on \m{U_i}, cf. \ref{can_pr1}).
\end{subsub}

\end{sub}

\sepsub

\Ssect{The primitive scheme ${\bf X}_2$ -- Proof of theorem \ref{theo2}}{X2}

It is the primitive double scheme defined by \ \m{\rho'\in 
H^1(T_{\P_2}\ot\ko_{\P_2}(-3))} (cf. \ref{can_pr2}). It is obtained from the 
cocycle \m{(\delta_{ij}^*)}, with \ \m{\delta_{ij}^*=\phi_{D_{ij},\nu_{ij}}}, 
where
\[D_{i,i+1} \ = \ \frac{x_i}{x_{i+1}}\frac{\partial}{\partial(x_{i+2}/x_i)} \ , 
\qquad \nu_{ij} \ = \ \Big(\frac{x_i}{x_j}\Big)^3 \]
(with the notations of \ref{can_sh2}).

Let $p$ be an integer, \m{\not=0}. The line bundle \m{\ko_{\P_2}(p)} is defined 
by the cocycle \m{\Big(\Big(\dsp\frac{x_i}{x_j}\Big)^{-p}\Big)}.

Consider the exact sequence \ \m{0\to\ko_{\P_2}(-3)\to\Omega_{{\bf 
X}_2|\P_2}\to\Omega_{\P_2}\to 0}, and the associated map \ 
\m{\delta:H^1(\Omega_{\P_2})\to H^2(\ko_{P_2}(-3))}. According to theorem 
\ref{prop4}, \m{\ko_{\P_2}(p)} can be extended to a line bundle on \m{{\bf 
X}_2} if and only if \ \m{\delta(\nabla_0(\ko_{\P_2}(p)))\not=0}.

The class \m{\nabla_0(\ko_{\P_2}(p))} is represented by the cocycle 
\m{\dsp\Big(\Big(-p\frac{x_j}{x_i}d(x_i/x_j)\Big)\Big)}.
From \ref{desc_delta3},\Nligne \m{\delta(\nabla_0(\ko_{\P_2}(p)))} is 
represented by the cocycle \m{(\lambda_{012})}, with
\[\lambda_{012} \ = \ \frac{1}{x_0^3}\frac{x_0}{x_1}\begin{pmatrix}\dsp
\frac{\frac{\partial}{\partial{x_2/x_0}}\Big(\Big(\frac{x_1}{x_2}\Big)^{-p}\Big)}
{\Big(\Big(\frac{x_1}{x_2}\Big)^{-p}\Big)} 
\end{pmatrix} \ = \ \frac{p}{x_0x_1x_2} \ . \]
It follows from \ref{can_pr1} that \ 
\m{\delta(\nabla_0(\ko_{\P_2}(p)))\not=0}. This proves theorem \ref{theo2}.

\end{sub}

\sepsub

\Ssect{The primitive scheme ${\bf X}_4$ -- Proof of theorem \ref{theo3}}{X4}

The computations are exactly the same as in \ref{X2}.

\end{sub}

\sepsec

\section{The case of projective bundles over curves}\label{pro_c}

Let $C$ be a smooth irreducible projective curve of genus $g$
and $E$ a rank 2 vector bundle 
on $C$. Let \ \m{X=\P(E)} \ be the associated projective bundle and \ 
\m{\pi:\P(E)\to C} \ the canonical projection. For every closed point \m{x\in 
C}, \m{\pi^{-1}(x)} is the projective line \m{\P(E_x)} of lines in \m{E_x}.

\sepsub

\Ssect{Preliminaries}{can_proj}

The line bundles on $X$ are of the form
\[\kl \ = \ \pi^*(D)\ot\ko_{\P(E)}(k) \ , \]
for some uniquely determined \m{D\in\Pic(C)} and integer $k$.
We have then
\begin{eqnarray*}
\pi_*(\kl) & = & (S^kE^*)\ot D \quad \text{if} \quad k\geq 0 ,\\
& = & 0 \quad \text{if} \quad k<0 ,\\
R^1\pi_*(\kl) & = & (S^{-2-k}E)\ot\det(E)\ot D \quad \text{if} \quad k\leq -2 
,\\
& = & 0 \quad \text{if} \quad k>-2 .
\end{eqnarray*}
it follows that
\begin{eqnarray*}
H^0(\kl) & = & H^0((S^kE^*)\ot D) \quad \text{if} \quad k\geq 0 ,\\
& = & \nsp \quad \text{if} \quad k<0 ,\\
H^1(\kl) & = &  H^1((S^kE^*)\ot D) \quad \text{if} \quad k\geq 0 ,\\
& = & H^0((S^{-2-k}E)\ot\det(E)\ot D) \quad \text{if} \quad k\leq -2 ,\\
& = & \nsp \quad \text{if} \quad k=-1 ,\\
H^2(\kl) & = &  H^1((S^{-2-k}E)\ot\det(E)\ot D) \quad \text{if} \quad k\leq -2 
,\\
& = & \nsp \quad \text{if} \quad k>-2 .
\end{eqnarray*}
We have canonical exact sequences
\begin{equation}\label{equ1b}
0\lra\Omega_{\P(E)/C}\lra\pi^*(E^*)\ot\ko_{\P(E)}(-1)\lra\ko_{\P(E)}\lra 0 \ 
, \end{equation}
\[0\lra\pi^*(\omega_C)\lra\Omega_{\P(E)}\lra \Omega_{\P(E)/C}\lra 0 \ . \]

\sepprop

\begin{subsub}\label{prop29}{\bf Proposition:} If $E$ is simple and \m{g\geq 
2}, we have \ \m{H^0(T_{\P(E)})=\nsp}.
\end{subsub}
\begin{proof} Since \m{g\geq 2}, we have \ \m{H^0(\omega_C^*)=\nsp}, and 
\m{H^0(T_{\P(E)})\simeq H^0((\Omega_{\P(E)/C})^*)} (by the second exact 
sequence). Applying \m{\pi_*} to the dual of the first exact sequence we find 
\[0\lra\ko_C\lra E\ot 
E^*\lra\pi_*((\Omega_{\P(E)/C})^*)\simeq\text{Ad}(E)\lra 0 \ , \]
hence \ \m{H^0((\Omega_{\P(E)/C})^*)=H^0(\text{Ad}(E))=\nsp} \ because $E$ is 
simple.
\end{proof}

\sepprop

Using proposition \ref{prop28} we get

\sepprop

\begin{subsub}\label{coro4}{\bf Corollary: } Suppose that $E$ is simple. Let 
\m{X_2} be a non trivial primitive double scheme with support \m{\P(E)} and 
associated line bundle \m{\omega_{\P(E)}}. Let \m{\C\sigma} be the element of 
\m{\P(H^1(\Omega_{\P(E)}))} associated do \m{X_2}. Then \m{X_2} can be extended 
to a primitive multiple scheme of multiplicity 3 if and only if \ 
\m{\sigma.\nabla_0(\omega_{\P(E)})=0} in \m{H^2(\omega_{\P(E)})}.
\end{subsub}

\end{sub}

\sepsub

\Ssect{Annulation of obstructions}{ann_lem}

Let \m{\Gamma,F\in\Pic(C)}, \m{\gamma=\deg(\Gamma)}, $n$, $k$, $p$ integers, 
with \m{n,k>0}. We will consider 3 cases:
\begin{enumerate}
\item[(i)] $E$ is semi-stable and \m{\deg(E)=0},
\item[(ii)] $E$ is semi-stable and \m{\deg(E)=-1},
\item[(iii)] $E$ is not semi-stable of degree $\epsilon=0$ or $-1$. Let 
\m{L_1\subset E} be the Harder-Narasimhan filtration of $E$, where \m{L_1} is a 
line bundle of degree \m{\epsilon_1>\epsilon}.
\end{enumerate}

\sepprop

\begin{subsub}\label{lem8}{\bf Lemma: } We have \ \m{h^0((S^{kn+p}E)\ot 
F\ot\Gamma^n)=0} if \ \m{\gamma<\gamma_0}, and
\ \m{\dsp\gamma_0=-\frac{\deg(F)}{n}} \ in case (i), \
\m{\dsp\gamma_0=\frac{k}{2}-\frac{\deg(F)}{n}+\frac{p}{2n}} \ in case (ii),
\ \m{\dsp\gamma_0=-\frac{\deg(F)}{n}-\frac{p\epsilon_1}{n}-k\epsilon_1} \ in 
case (iii).
\end{subsub}
\begin{proof} In case (i), \m{S^{kn+p}E} is semi-stable of degree 0, so 
\m{h^0((S^{kn+p}E)\ot F\ot\Gamma^n)=0} if\Nligne \m{\deg(F\ot\Gamma^n)<0}, 
which 
is equivalent to \m{\dsp\gamma<-\frac{\deg(F)}{n}}. The proof in case (ii) is 
similar. To prove case (iii) we use the fact that there is a filtration of 
\m{S^{kn+p}E} with graduates \m{L_1^a\ot(\ko_C/L_1)^{kn+p-a}}, \m{0\leq a\leq 
kn+p}.\end{proof}
 
\sepprop

\begin{subsub}\label{lem9}{\bf Lemma: } We have \ \m{h^0((S^{kn+p}E)\ot E\ot 
F\ot\Gamma^n)=0} if \ \m{\gamma<\gamma_0}, and
\ \m{\dsp\gamma_0=-\frac{\deg(F)}{n}} \ in case (i), \
\m{\dsp\gamma_0=\frac{k}{2}-\frac{\deg(F)}{n}+\frac{p+1}{2n}} \ in case (ii),
\ \m{\dsp\gamma_0=-\frac{\deg(F)}{n}-\frac{(p+1)\epsilon_1}{n}-k\epsilon_1} \ 
in case (iii).
\end{subsub}
\begin{proof} Similar to that of lemma \ref{lem8}. \end{proof}

\sepprop

Let $n$, $k$ be positive integers, and \ \m{L=\pi^*(D)\ot\ko_{\P(E)}(-k)} \ a 
line bundle on \m{\P(E)}. If we want to study the primitive multiple schemes 
with associated smooth variety \m{\P(E)} and associated line bundle $L$, we 
need to consider the cohomology groups \m{H^2(L^n)}, \m{H^2(T_X\ot L^n)}, 
\m{n>0}.

\sepprop

\begin{subsub}\label{obs1} Obstructions to the extension of line bundles in 
higher multiplicity -- \rm We have \Nligne 
\m{H^0(L^n)=\nsp}, \m{H^1(L^n)\simeq H^0((S^{kn-2}E)\ot\det(E)\ot D^n)},\Nligne
\m{H^2(L^n)\simeq H^1((S^{kn-2}E)\ot\det(E)\ot D^n)}, and by Serre duality
\[h^1((S^{kn-2}E)\ot\det(E)\ot D^n) \ = \ h^0((S^{kn-2}E^*)\ot\det(E^*)\ot 
D^{-n}\ot\omega_C) \ = \ h^0((S^{kn-2}E)\ot F\ot\Gamma^n) \ , \]
with \ \m{F=\det(E)\ot\omega_C}, \m{\Gamma=D^*\ot\det(E)^{-k}}. Hence, by lemma 
\ref{lem8}, \m{h^2(L^n)=0} \ if \Nligne
\m{\deg(D)>\delta_0}, with \ \m{\dsp\delta_0=\frac{2g-2}{n}} in case (i), 
\m{\delta_0=\dsp\frac{k}{2}+\frac{2g-2}{n}} in case (ii), and \Nligne 
\m{\delta_0=\dsp\frac{2g-2+\deg(E)-2\epsilon_1}{n}+k\epsilon_1-k\deg(E)} \ in 
case (iii).
\end{subsub}

\sepprop

\begin{subsub} Obstructions to the extension of the schemes in higher 
multiplicity -- \rm We have exact sequences
\begin{equation}\label{equ11}
0\lra L^n\lra\pi^*(E)\ot\ko_{\P(E)}(1)\ot L^n\lra\Omega_{\P(E)/C}^*\ot 
L^n\lra 0 \ ,\end{equation}
\[0\lra\Omega_{\P(E)/C}^*\ot L^n\lra T_{\P(E)}\ot L^n\lra\pi^*(\omega_C^*)\ot 
L^n\lra 0 \ . \]
Hence we have \ \m{H^2(T_{\P(E)}\ot L^n)=\nsp} \ whenever
\[h^2(\pi^*(E)\ot\ko_{\P(E)}(1)\ot L^n) \ = \ h^2(\pi^*(\omega_C^*)\ot L^n)
\ = \ 0 . \]
We have
\[h^2(\pi^*(E)\ot\ko_{\P(E)}(1)\ot L^n) \ = \ h^1((S^{kn-3}E)\ot E\ot\det(E)
\ot D^n) \ , \]
and by Serre duality
\begin{eqnarray*} h^1((S^{kn-3}E)\ot E\ot\det(E)\ot D^n) & = &
h^0((S^{kn-3}E^*)\ot E^*\ot\det(E^*)\ot D^{-n}\ot\omega_C)\\
& = & h^0((S^{kn-3}E)\ot E\ot F\ot\Gamma^n) \ ,
\end{eqnarray*}
with \ \m{F=\det(E)\ot\omega_C}, \m{\Gamma=D^*\ot\det(E)^{-k}}. Hence, by lemma 
\ref{lem9},\Nligne \m{h^2(\pi^*(E)\ot\ko_{\P(E)}(1)\ot L^n)=0} \ if the same 
conditions as in \ref{obs1} are satisfied.

We have \ \m{H^2(\pi^*(\omega_C^*)\ot L^n)\simeq H^1((S^{kn-2}E)\ot\det(E)\ot
\omega_C^*\ot D^n)}, and
\[h^1((S^{kn-2}E)\ot\det(E)\ot\omega_C^*\ot D^n) \ = \
h^0((S^{kn-2}E)\ot F\ot\Gamma^n) \ , \]
with \ \m{F=\det(E)\ot\omega_C^2}, \m{\Gamma=D^*\ot\det(E)^{-k}}. Hence, by 
lemma \ref{lem8}, \m{h^2(\pi^*(\omega_C^*)\ot L^n)=0} \ if
\m{\deg(D)>\delta_0}, with \ \m{\dsp\delta_0=\frac{4g-4}{n}} in case (i), 
\m{\delta_0=\dsp\frac{k}{2}+\frac{4g-4}{n}} in case (ii), and \Nligne 
\m{\delta_0=\dsp\frac{4g-4+\deg(E)-2\epsilon_1}{n}+k\epsilon_1-k\deg(E)} \ in 
case (iii).
\end{subsub}

\end{sub}

\sepsub

\Ssect{Construction of primitive multiple schemes}{cnst}

Suppose that \ \m{L=\pi^*(D)\ot\ko_{\P(E)}(-k)}, with \m{k\geq 3}. Let \ 
\m{d=\deg(D)}. The following is a consequence of \ref{ann_lem}:

\sepprop

\begin{subsub}\label{lem10}{\bf Lemma: } We have \ \m{H^2(L^n)=\nsp} \ for 
every \m{n\geq 1} and \ \m{H^2(T_{\P(E)}\ot L^n)=\nsp} \ for every \m{n\geq 2} 
if

In case (i):
\begin{enumerate}
\item[--] If $g=0$, \ $d\geq 0$.
\item[--] If $g=1$, \ $d>0$.
\item[--] If $g\geq 2$, \ $d>2g-2$. 
\end{enumerate}

In case (ii):
\begin{enumerate}
\item[--] If $g=0$, \ $d\geq\dsp\frac{k}{2}$.
\item[--] If $g=1$, \ $d>\dsp\frac{k}{2}$.
\item[--] If $g\geq 2$, \ $d>\dsp\frac{k}{2}+2g-2$.
\end{enumerate}

In case (iii):
\begin{enumerate}
\item[--] If $g=0$ or $1$, \ $d\geq k(\epsilon_1-\deg(E))$.
\item[--] If $g\geq 2$, \ $d\geq k(\epsilon_1-\deg(E))+2g-2$.
\end{enumerate}
\end{subsub}

\sepprop

Recall that to extend a primitive scheme \m{X_n} of multiplicity \m{n\geq 2} to 
one of multiplicity \m{n+1}, we must first extend the ideal sheaf of \m{\P(E)} 
in \m{X_n} to a line bundle $\L$ on \m{X_n} (with an obstruction in 
\m{H^2(L^{n-1})}), and then to extend \m{X_n} to \m{X_{n+1}} such that the 
ideal sheaf of \m{\P(E)} in \m{X_{n+1}} is $\L$, we have another obstruction in 
\m{H^2(T_{\P(E)}\ot L^n)}. Hence if the conditions of lemma \ref{lem10} are 
satisfied, the obstructions vanish, and it is possible to extend a primitive 
multiple scheme of multiplicity $n$ to one of multiplicity \m{n+1}.

Similarly we have

\sepprop

\begin{subsub}\label{lem12}{\bf Lemma: } We have \ \m{H^2(T_{\P(E)}\ot L)=\nsp} 
if

In case (i):
\begin{enumerate}
\item[--] If $g=0$, \ $d>-2$.
\item[--] If $g=1$, \ $d>0$.
\item[--] If $g\geq 2$, \ $d>4g-4$. 
\end{enumerate}

In case (ii):
\begin{enumerate}
\item[--] If $g=0$, \ $d>\dsp\frac{k}{2}-2$
\item[--] If $g=1$, \ $d>\dsp\frac{k}{2}$
\item[--] If $g\geq 2$, \ $d>\dsp\frac{k}{2}+4g-4$
\end{enumerate}

In case (iii):
\begin{enumerate}
\item[--] If $g=0$, \ $d>k(\epsilon_1-\deg(E))+\deg(E)-2-2\epsilon_1$.
\item[--] If $g=1$, \ $d>k(\epsilon_1-\deg(E))+\deg(E)-2\epsilon_1$.
\item[--] If $g\geq 2$, \ $d>k(\epsilon_1-\deg(E))+\deg(E)-2\epsilon_1+4g-4$.
\end{enumerate}
\end{subsub}

\sepprop

\begin{subsub}\label{lem11}{\bf Lemma: } Suppose that \m{n\geq 1}. Then the map
\ \m{H^1(L^n)\to H^1(\pi^*(E)\ot\ko_{\P(E)}(1)\ot L^n)} \ induced by 
$(\ref{equ11})$ is injective.
\end{subsub}
\begin{proof}
This map is the canonical one
\[H^0(D^n\ot(S^{kn-2}E)\ot\det(E))\lra H^0(D^n\ot E\ot (S^{kn-3}E)\ot\det(E)) \ 
, \]
which is clearly injective.
\end{proof}

\sepprop

Suppose that the conditions of lemmas \ref{lem10} and \ref{lem12} are 
satisfied and \m{n\geq 1}. We have \Nligne \m{H^0(\Omega_{\P(E)/C}^*\ot 
L^n)=\nsp}, because for every \m{x\in C}, \m{\Omega_{\P(E)/C}^*\ot 
L^n_{|\pi^{-1}(x)}\simeq\ko_{\P(E_x)}(2-kn)}. Hence
\[h^1(\Omega_{\P(E)/C}^*\ot L^n) \ = \ h^1(\pi^*(E)\ot\ko_{\P(E)}(1)\ot L^n)-
h^1(L^n) \ . \]
We have then
\begin{eqnarray*}
h^1(T_{\P(E)}\ot L^n) & = & h^1(\Omega_{\P(E)/C}^*\ot L^n)+
h^1(\pi^*(\omega_C^*)\ot L^n) \\
& = & h^1(\pi^*(E)\ot\ko_{\P(E)}(1)\ot L^n)-h^1(L^n)+
h^1(\pi^*(\omega_C^*)\ot L^n)\\
& = & h^0(D^n\ot E\ot (S^{kn-3}E)\ot\det(E))-h^0(D^n\ot(S^{kn-2}E)\ot\det(E))\\
& & +h^0(D^n\ot\omega_C^*\ot(S^{kn-2}E)\ot\det(E))\\
& = & \chi(D^n\ot E\ot 
(S^{kn-3}E)\ot\det(E))-\chi(D^n\ot(S^{kn-2}E)\ot\det(E))\\
& & +\chi(D^n\ot\omega_C^*\ot(S^{kn-2}E)\ot\det(E))\\
& = & n(kn-2)(2d+k\deg(E))+2(1-g)(2kn-3) \ . 
\end{eqnarray*}
It is easily verified that \ \m{h^1(T_{\P(E)}\ot L^n)>0}. It follows that {\em 
there exist infinite sequences
\[X_1=\P(E)\subset X_2\subset\cdots\subset X_n\subset X_{n+1}\subset\cdots \ ,
\]
where \m{X_n} is a non trivial projective primitive multiple scheme with 
associated smooth variety \m{\P(E)} and associated line bundle $L$}.

We have also
\begin{eqnarray*}h^1(L^n) & = & \chi(D^n\ot(S^{kn-2}E)\ot\det(E))\\
& = & (kn-1)\big(\frac{kn}{2}\deg(E)+nd+1-g\big) \ .
\end{eqnarray*}

\sepprop

The space \m{\P(H^1(T_{\P(E)}\ot L))} parametrizes the primitive double schemes 
(cf. \ref{param_2}).

Let \m{X_n} be a non trivial primitive multiple scheme of multiplicity $n$. Let 
\m{\kx} be the set of extensions of \m{X_n} to a primitive multiple scheme 
\m{X_{n+1}} of multiplicity \m{n+1}. Then from \ref{ext_MV}, there is a 
canonical surjective map
\[H^1((\Omega_{X_2|\P(E)})^*\ot L^n)\lra\kx\]
whose fibers are the orbits of an action of \m{\Aut(X_n)}. As a set, we have 
\Nligne \m{\Aut_0(X_n)\simeq H^0(\kt_{n-1}\ot\ki_{\P(E),X_n})} (cf. theorem 
\ref{prop22}). Using corollary \ref{coro2} and the exact sequence \ \m{0\to 
T_{\P(E)}\to(\Omega_{X_2|\P(E)})^*\to L^{-1}\to 0}, it is easy to see that \ 
\m{H^0(\kt_{n-1}\ot\ki_{\P(E),X_n})=\nsp}. From \ref{coro3}, \m{\Aut(X_n)} is 
finite. Hence $\kx$ is parametrized, up to the action of a finite group, by 
\m{H^1((\Omega_{X_2|\P(E)})^*\ot L^n)}. We have
\begin{eqnarray*}
h^1((\Omega_{X_2|\P(E)})^*\ot L^n) & = & h^1(T_{\P(E)}\ot L^n)+h^1(L^{n-1})\\
& = & (3kn^2-5n-2kn+k+1)\big(\frac{k}{2}\deg(E)+d\big)+2(1-g)(5kn-7-k) \ .
\end{eqnarray*}

\end{sub}

\sepsub

\Ssect{The case of \m{\P_1\times\P_1}}{arg2}

We have here \m{C=\P_1}, \m{E=\ko_{\P_1}\ot\C^2}, $\pi$ is the projection
\m{\P_1\times\P_1\to\P_1} \ on the first factor, and, with the usual notations, 
\m{\ko_{\P(E)(1)}=\ko(0,1)}, \m{\pi^*(\ko_{\P_1}(1))=\ko(1,0)}.

Suppose that \m{k\geq 3}, and \m{L=\ko(d,-k)}. The conditions of lemmas 
\ref{lem10} and \ref{lem12} are: \m{d>0}.

There exist infinite sequences
\[X_1=\P_1\times\P_1\subset X_2\subset\cdots\subset X_n\subset 
X_{n+1}\subset\cdots \ , \]
where \m{X_n} is a non trivial projective primitive multiple scheme with 
associated smooth variety \m{\P_1\times\P_1} and associated line bundle $L$. 

The extensions of \m{X_n} to a primitive multiple scheme of 
multiplicity \m{n+1} form a family of dimension \ 
\m{d(3kn^2-5n-2kn+k+1)+5kn-7-k}.

\end{sub}

\vskip 4cm


\vskip 1.5cm

\begin{thebibliography}{99} 
\bibitem{b_g_g} Bangere, P., Gallego, F.J., González, M. {\em Deformations of 
hyperelliptic and generalized hyperelliptic polarized varieties.} 
https://arxiv.org/abs/2005.00342 (2020).
\bibitem{b_m_r} Bangere, P., Mukherjee, J., Raychaudhry, D. {\em K3 carpets on 
minimal rational surfaces and their smoothing.}  Internat. J. Math. 32 (2021), 
no. 6, Paper No. 2150032
\bibitem{ba_fo} B\u anic\u a, C., Forster, O. {\em Multiple structures on space
curves.} In: Sundararaman , D. (Ed.) Proc. of Lefschetz Centennial Conf.
(10-14 Dec. Mexico), Contemporary Mathematics 58, AMS, 1986, 47-64.
\bibitem{ba_ei}Bayer, D., Eisenbud, D. {\em Ribbons and their canonical
embeddings.} Trans. of the Amer. Math. Soc., 1995, 347-3, 719-756.
\bibitem{ch-ka} Chen, D., Kass, J. L. {\em Moduli of generalized line bundles 
on a ribbon.} J. of Pure and Applied Algebra 220 (2016), 822-844.
\bibitem{dr1b}Dr\'ezet, J.-M. {\em D\'eformations des extensions larges de
faisceaux}. Pacific Journ. of Math. 220, 2 (2005), 201-297.
\bibitem{dr2}Dr\'ezet, J.-M. {\em Faisceaux coh\'erents sur les courbes
multiples} . Collect. Math. 2006, 57-2, 121-171.
\bibitem{dr1}Dr\'ezet, J.-M. {\em Param\'etrisation des courbes multiples
primitives} Adv. in Geom. 2007, 7, 559-612.
\bibitem{dr4}Dr\'ezet, J.-M. {\em Faisceaux sans torsion et faisceaux quasi
localement libres sur les courbes multiples primitives.} Mathematische
Nachrichten, 2009, 282-7, 919-952.
\bibitem{dr5}Dr\'ezet, J.-M. {\em Sur les conditions d'existence des faisceaux
semi-stables sur les courbes multiples primitives.} Pacific Journ. of Math.
2011, 249-2, 291-319.
\bibitem{dr6}Dr\'ezet, J.-M. {\em Courbes multiples primitives et
d\'eformations de courbes lisses.} Annales de la Facult\'e des Sciences de
Toulouse 22, 1 (2013), 133-154.
\bibitem{dr7}Dr\'ezet, J.-M. {\em Fragmented deformations of primitive
multiple curves.} Central European Journal of Mathematics 11, n${}^o$ 12
(2013), 2106-2137.
\bibitem{dr8}Dr\'ezet, J.-M. {\em Reducible deformations and
smoothing of primitive multiple curves.} Manuscripta Mathematica 148 (2015), 
447-469.
\bibitem{dr9}Dr\'ezet, J.-M. {\em Reachable sheaves on ribbons
and deformations of moduli spaces of sheaves.} Intern. J. of Math. 28,12 
(2017), 1750086.
\bibitem{ei_gr}Eisenbud, D., Green, M. {\em Clifford indices of ribbons.}
Trans. of the Amer. Math. Soc., 1995, 347-3, 757-765.
scientifiques et industrielles 1252, Hermann, Paris (1964). 
\bibitem{fo} Fong, L-Y. {\em Rational ribbons and deformation of hyperelliptic 
curves.} Journ. Alg. Geom. 2, n${}^o$ 2 (1993), 295-307.
\bibitem{fe}Ferrand, D. {\em Courbes gauches et fibr\'es de rang 2.} C.R. Acad. 
Sci. Paris, 281 (1977), 345-347.
\bibitem{fr}Frenkel, J. {\em Cohomologie non ab\'elienne et espaces fibr\'es.}
Bull. de la Soc. Math. de France 85 (1957), 135-220. 
\bibitem{ga_go_pu0} Gallego, F.J., González, M., Purnaprajna, B.P. {\em
Deformation of finite morphisms and smoothing of ropes.} Compos. Math. 144 
n${}^o$ 3 (2008), 673-688.
\bibitem{ga_go_pu1} Gallego, F.J., González, M., Purnaprajna, B.P. {\em An 
infinitesimal condition to smooth ropes.} Rev. Mat. Complut. 26 n${}^o$ 1, 
(2013), 253-269.
\bibitem{ga-go-pu}Gallego, F.J., Gonz\'alez, M., Purnaprajna, B.P. {\em K3 
double structures on Enriques surfaces and their smoothings.} J. of Pure and 
Applied Algebra 212 (2008), 981-993.
\bibitem{go} Godement, R. {\em Topologie alg\'ebrique et th\'eorie des 
faisceaux.} Actualit\'es scientifiques et industrielles 1252, Hermann, Paris 
(1964).
\bibitem{gonz1} Gonz\'alez, M. {\em Smoothing of ribbons over curves.} J. reine 
angew. Math. 591 (2006), 201-235.
\bibitem{sga1}Grothendieck, A. et al. SGA1. {\em Rev\^etements Etales et
Groupe Fondamental.} SGA1. Lect. Notes in Math. 224. Springer-Verlag (1971).
\bibitem{ha2} Hartshorne, R. {\em Ample vector bundles.} Publ. Math. IHES , 29 
(1966), 319-350.
\bibitem{ha} Hartshorne, R. {\em Algebraic geometry.} Grad. Texts in Math.,
Vol. 52, Springer (1977).
\bibitem{sa1} Savarese, M. {\em On the irreducible components of the 
compactified Jacobian of a ribbon.}\Nligne https://arxiv.org/abs/1803.05360 .
\bibitem{sa2} Savarese, M. {\em Coherent Sheaves on Ribbons and their Moduli.} 
https://arxiv.org/abs/1902.08510 .
\bibitem{sa3} Savarese, M. {\em Generalized line bundles on primitive multiple 
curves and their moduli.}\Nligne https://arxiv.org/abs/1902.09463 .
\end{thebibliography}
\end{document}